# Bethe仮説へのいざない

坂本玲峰（東京理科大学）

## 目 次

## 1 はじめに—量子可積分系とは

**はじめに**　本稿は量子可積分系における代表的な理論の一つである Bethe 仮説 (Bethe ansatz) 法についての入門的内容及び最近の進展についての解説です。読者としては物理学者、数学者の両方を念頭に置いて執筆しました。Bethe 仮説法には従来不透明な面が多くあり、数学者には接近しがたい印象を持たれていた方も多いようですが、最近の進展によりだいぶ色々と様子が分かってきたようです。本稿の内容は、本質的には最高ウェイトベクトルの具体的構成ですので、多くの数学者の方々にご興味をお持ちいただければと思います。内容は 2016 年 9 月 7 日に第 61 回代数学シンポジウム（日本数学会、佐賀大学）で行った筆者の講演に基づきます。講演で使用した資料は下記 URL

https://sites.google.com/site/affinecrystal/resources

にて入手可能です（本稿では色つき文字にはハイパーリンクがついています）。

Bethe 仮説法の分野ではすでに膨大な数の論文（Google Scholar で “Bethe ansatz” と検索すると約 2 万件表示される）が出版されています。本稿では Bethe の原論文に沿った内容を紹介しますが、他にも多くの切り口（例えば Gaudin 模型と幾何学的ラングランズ対応 [F95] など）が存在しますので、本稿の内容に限らず、ご専門、ご関心の方向に応じて興味をお持ちいただければ幸いです。

**古典力学における積分可能な系**　さて、表題の量子可積分系とはいかなる分野なのか、簡単に紹介するところから始めよう。歴史的経緯に基づきまずは解析力学について振り返ることから始めよう（詳細は、例えば優れた教科書 [A] を参照）。解析力学は古典力学の理論的研究であり、Euler, Lagrange, Jacobi, Poincaré 等をはじめとする多くの研究者によって研究され、力学のみならず物理学、数学の広範な領域にわたって多大な影響を与えた。20 世紀に入り、原子、分子などの微視的領域では古典力学は成立せず、量子力学と呼ばれる全く異なる力学によって記述されることが発見された。古典力学と量子力学は大きく性格が異なる理論であるから、量子力学の場合にも解析力学と同様に「深い」理論的研究を進めようと考えるのは自然な事と思われる。

ただしここで「深い」理解と述べたのは、何か特定の結果を想定している訳ではない。むしろ折角量子力学という新しい力学がある以上、古典力学の結果にこだわりすぎることなく、新しく面白い数学が見つかれば十分有益であると思う。そうは言っても我々にはすでに古典力学における数世紀にわたる蓄積があるのであるから、その過程、とくに研究者の問題意識と、その結果として得られた成果は大変参考になる。ここでは Jacobi[1] の研究について見てみよう（詳細は [L] の XVI 章や [H] を参照）。Jacobi の有名な力学の教科書 [J] では彼の研究の方針について次のような言明がある。

> この講義は、運動方程式の積分に際して、方程式の特別な形から得られる手法に関するものである。ラグランジュの『解析力学』には、定式化や微分方程式の変換の問題に関するあらゆることを見出せるが、それらの積分に関してはほとんど議論されていない。. . .

結果として Jacobi の研究が挙げた重要な成果を鑑みれば、彼の問題意識は大変興味深いものと思われる。その後の展開を見れば、可積分な力学系の性質が深く理解されるとともに、逆にその他大多数の積分不能な系について正確に認識できるようになった点も見逃せない（実際には、ほとんどの古典力学系は、数値計算すら原理的に不可能で、人間の論理のはるかに及ばない領域に存在するのであった[2]）。

[1] Carl Gustav Jacob Jacobi (1804–1851). 若いころから Poisson の研究などに興味を持っていたようであるが、Hamilton の力学の研究 (1835 年) に触発されて本格的に力学の研究に参入し、1843 年に過労で倒れるまで精力的に研究を行った。生前に研究を完成させることはできなかったが、他界後に出版された原稿は膨大なものであった。これらの力学の研究は文字通り Jacobi 先生が命を削って行ったものだと思うと思わず襟を正したくなる。Poisson 可換な保存量の存在と系の可積分性に関する有名な Liouville の定理も、元々は Jacobi によって得られた定理を自身の講義 (1842–43 年) において紹介したものがかなり広い範囲に流布したものだった様である。Jacobi の研究は力学の発展に重大な影響をもたらしたのみならず、Sophus Lie による連続群の研究に直接の影響を及ぼす（いまでも Lie 代数の定義式に名前が残っている）など意義深いものであった。

[2] 蛇足ながらもう少し説明を加えておくと、複雑な古典力学系では大抵の場合初期値の誤差が指数関数的に増大して伝播する。従って有限精度の数値計算では長時間の時間発展は計算不能であるし、無限精度の数値計算を行うことは原理的に不可能である。もちろん厳密解など得られるはずもないのであるから、結局大抵の古典力学系は人間の論理の力では決して理解することができない。平たく言えば、例えば地球大気が平均して地球を数周回る程度の未来について天気予報のような形で予言することは原理的に不可能である。これらは

**量子力学における積分可能な系**　量子力学は、古典力学と異なり、出来てからまだ百年にも満たない若い分野であるから、まだいろいろと未開拓の領域も残っているかもしれない。そこで本稿では量子力学における積分可能な系を考えたい。まず量子力学において解くべき方程式など、量子力学の基礎事項で必要となる部分から始めたい（量子力学の入門書としては、例えば [M] は洞察に富んでいて興味深い）。量子力学における基礎方程式は Schrödinger 方程式とよばれ、次のような固有値方程式の形をしている。

$$H\psi = E\psi \tag{1}$$

各項の意味づけは以下の通り。

- $H$: 演算子（Hamiltonian と呼ばれ物理系を表す）
- $\psi$: 固有ベクトル（波動関数、系の運動の状態—$|\psi|^2$ が粒子の確率分布となる）
- $E$: 固有値（観測される系のエネルギー）

物理量を通常のスカラーではなく作用素と考える点が古典力学と大きく異なる点である。その固有値が実際に実験的に観測される物理量となる。Hamiltonian は古典力学における全エネルギーに対応する演算子で、通常微分演算子となるが、本稿では有限サイズの行列となる場合を主に扱う。

古典力学における複雑さから容易に想像されるように、Schrödinger 方程式は解くことが非常に難しい。特に相互作用する多数の粒子が存在するような場合は一般には非常に大規模な数値計算に頼らざるを得ない。しかし時として驚くほど深く数学的な解析が可能な量子系が存在する。本稿の主題である Heisenberg 模型はその代表的な例であり、その解析手段である Bethe 仮説法は数学的にみても大変興味深い内容である。

**一つの例**　考え方の方向性を見るために、量子可積分系について本稿で関心を寄せる側面と通じるものが現れる他の典型的な例について触れておくのも興味深いかもしれない。周の長さが $L$ の円周上に $N$ 個の量子力学的な粒子を配置し、それぞれが距離の 2 乗分の 1 の相互作用ポテンシャルを持つ時（相互作用そのものは距離の 3 乗分の 1 に比例する）Calogero–Sutherland 模型と呼ばれる。参考までに、ハミルトニアンは

$$H = \sum_{i=1}^{N} \frac{1}{2} p_i^2 + \beta(\beta - 1) \sum_{i<j} \frac{(\frac{\pi}{L})^2}{\sin^2 \frac{\pi}{L}(q_i - q_j)},$$

ここで $q_i\ (0 \leq q_i \leq L)$ は各粒子の円周上の座標であり、$p_i = -\sqrt{-1}\frac{\partial}{\partial q_i}$ は運動量に対応する微分演算子である（$\beta$ は複素パラメータ）。

この模型は驚くほどすっきりと解く事（全ての波動関数 $\psi$ を求める事）ができる [Su71, Su72]。類似の積分可能な系と比べても際立って性質の良い模型である。その様な場合背後に深い数学的構造が潜んでいると考えるのが自然である。実際無限次元 Lie 代数の著名な例である Virasoro 代数（例えば [KR] 参照）の表現の特異ベクトルが Calogero–Sutherland 模型の波動関数（Jack 多項式と呼ばれる）の特殊な場合と一致することが知られていた [MY]。実はこの関係はより一般的なものであり、Jack 多項式を Virasoro 代数の表現の基底として

---

Poincaré 以来の研究によりすでに確立された事実である。しばしば、特に量子力学との比較の文脈で、「古典力学は決定論的な系である」などと耳にするが、この様な事実を踏まえればその様な言明は誤りである。

選ぶと、Virasoro 代数の作用が深い数学的背景があることを窺わせる形（組み合わせ論的意味がつく正整数係数一次式に因数分解される有理関数の和）で具体的に表示されるという極めて非自明な予想がある [SSAFR]。この予想は特殊な場合には解決しているが [CJ]、一般的な解決は容易ではなさそうである。ともあれ、以上の結果は Virasoro 代数の表す無限次元対称性の核心部に Jack 多項式が深くかかわっていることを示しており、未知の深い構造の存在を示している。

**本稿の内容**　§2 では本稿で主に対象とする Heisenberg 模型を定義する。§3 では、Faddeev 等の定式化に沿って Bethe 仮説法の天下り的な紹介をする。§4 では、理論全体のかなめとなる Bethe 仮説方程式（9 ページ (5) 式）の解について必要な事項を解説する。Bethe ベクトルとは Bethe 仮説方程式の解によって記述される最高ウェイトベクトルであり、Heisenberg 模型の波動関数となっている。§5 では、Bethe ベクトルを求める際に障害となる特異解について最近の進展を紹介する。§6 では、Bethe 仮説方程式の解の構造を理解する上で長年障害となっていたストリングの概念の問題に対して、新しい量子数を定義して解決する提案について紹介する。§7 は量子群に関係する場合の文献の紹介である。

最後の §8 では、以上のような理論から生み出された数学の例として、結晶基底に関わる一連の研究についてかなり詳細に紹介する。上記ストリングの概念に触発されて発展した理論であり、関連する対称性は著名な無限次元代数の一つアフィン量子群である。本節の内容は単独で見ても代数的組み合わせ論における重要な進展となっており興味深い。なお §8 の内容は本稿の他の部分とは独立して読むことができる。本節の概要については §8.1 および §8.2 を参照されたい。

**話題**　本稿執筆中に 2016 年度ノーベル物理学賞の発表があった。受賞した Haldane, Kosterlitz, Thouless のお三方とも本稿の内容と関連のある分野での研究に対する受賞であった。特に Haldane 氏は本稿の主題である一次元 Heisenberg 模型を基盤とする研究（1983 年）に対しての受賞であった（氏は上記 Calogero 系に関する研究でも有名である）。なお、他のお二方の研究は二次元の古典的 Heisenberg 模型に関するものである。

## 2　Heisenberg 模型

本稿では主に Heisenberg[3] 模型と呼ばれる磁性体のモデルを扱う。これは Schroödinger[4] 方程式 (1) で Hamiltonian $H$ が行列となる場合である。具体的には、電子が $N$ 個一次元格子状に並んだ系を考える。本稿の内容を理解するにあたって物理的な詳細を理解する必要は無いが、記号の導入もかねて簡単に背景にも触れておく。日常的に目にする磁石の起源である電子は、それぞれが小さな磁石としての性質を持っている。ただし量子的な効果により上向きまたは下向きの二つの状態しか持たないという一風変わった性質を持ち、一般にはそれら二状態の重ね合わせとなる。数学的に表示するために、上向きと下向きの状態をそれぞれ

$$v_+ = \begin{pmatrix} 1 \\ 0 \end{pmatrix}, \quad v_- = \begin{pmatrix} 0 \\ 1 \end{pmatrix}$$

[3] Werner Heisenberg (1901–1976). 量子力学の建設に大きな功績を残した。1932 年ノーベル物理学賞受賞。

[4] Erwin Schrödinger (1887–1961). 1926 年に量子力学を完成させたが、一連の論文は半年程度の間に集中して全て単著で著されており、その後の量子力学の教科書の入門的事項がほぼ達成されているという驚異的ペースだった。それらの論文では解析力学の深い造詣も披露されている。1933 年ノーベル物理学賞受賞。

とベクトル表示し、電子のスピン状態を二次元の状態空間 $V \simeq \mathbb{C}^2$ の元とみなす。今電子が $N$ 個一次元状に並んでいるのであるから、考える波動関数 $\psi$ は空間

$$\mathfrak{H}_N = \bigotimes_{j=1}^{N} V_j, \quad V_j \simeq \mathbb{C}^2.$$

の元であると考える。

さて、以前 Schrödinger 方程式を導入した時に、量子力学の基本的枠組みとして、物理量を作用素として考えることを紹介した。現在の場合重要となる物理量は角運動量であり、空間の $x, y, z$ 軸方向に応じて以下の Pauli 行列で与えられる（$sl_2$ の生成子である）

$$\sigma^1 = \begin{pmatrix} 0 & 1 \\ 1 & 0 \end{pmatrix}, \quad \sigma^2 = \begin{pmatrix} 0 & -i \\ i & 0 \end{pmatrix}, \quad \sigma^3 = \begin{pmatrix} 1 & 0 \\ 0 & -1 \end{pmatrix}.$$

上で定義された $v_+$ と $v_-$ は角運動量の $z$ 軸成分 $\sigma^3$ のそれぞれ固有値 1 および $-1$ の固有ベクトル（物理的状態）となっている。

**ウォームアップ**　Heisenberg 模型を定義するために、基本となる長さ $N = 2$ の場合を考えてみよう。記号として空間 $\mathfrak{H}_N$ の $k$ 番目だけに非自明に作用する作用素を

$$\sigma_k^a = I \otimes \cdots \otimes \underbrace{\sigma^a}_{k} \otimes \cdots \otimes I$$

などと表すことにする。さて $\mathfrak{H}_2$ 上の次の作用素 $h$ を考えてみよう。

$$h := \sum_{a=1}^{3} \sigma_1^a \sigma_2^a = \begin{pmatrix} 1 & 0 & 0 & 0 \\ 0 & -1 & 2 & 0 \\ 0 & 2 & -1 & 0 \\ 0 & 0 & 0 & 1 \end{pmatrix}.$$

この行列を対角化すると以下の結果を得る。

| 固有値 | 固有ベクトル |
|---|---|
| 1 | $(1, 0, 0, 0)$ |
| 1 | $(0, 1, 1, 0)$ |
| 1 | $(0, 0, 0, 1)$ |
| $-3$ | $(0, 1, -1, 0)$ |

固有ベクトルをベクトル $v_+$ と $v_-$ を用いて具体的に表示してみると（記号については必要ならば §3 参照)、作用素 $h$ は $\mathbb{C}^2 \otimes \mathbb{C}^2$ 内の $sl_2$ の 3 次元表現

$$\{v_+ \otimes v_+, v_+ \otimes v_- + v_- \otimes v_+, v_- \otimes v_-\}$$

の上で固有値 1、1 次元表現

$$\{v_+ \otimes v_- - v_- \otimes v_+\}$$

の上で固有値 $-3$ を持つことが分かった。

さて、量子力学の基本的思想に立ち帰り、作用素 $h$ を隣接する二つの電子間に作用する相互作用を表す作用素であると考えよう。その場合観測される物理量（$h$ の固有値—電子間の相互作用エネルギー）は同じ対称性（表現）に属する電子の組み合わせの間で同じ値を持つことになるので、相互作用として自然であろうと考えられる。

**Heisenberg 模型の定義**　ここで物理的な仮定（単純化）として、隣り合う二つの電子間には $h$ と同じ相互作用が生じ、それより遠方の電子とは何らの相互作用も持たないとしよう。具体的には次のハミルトニアンを考察する。

$$\mathcal{H}_N = \frac{J}{4}\sum_{k=1}^{N}(\sigma_k^1\sigma_{k+1}^1 + \sigma_k^2\sigma_{k+1}^2 + \sigma_k^3\sigma_{k+1}^3 - \mathbb{I}_N) \tag{2}$$

ただし $\mathbb{I}_N$ は $2^N$ 次元単位行列で $J$ はパラメータである（これらの部分は、さしあたってはおまけのような部分である）。この系を Heisenberg 模型 ([H], 1928 年) とよぶ。ただし上記のような $\mathcal{H}_N$ の定式化は Dirac[5]による ([D], 1929 年)。

どのような条件の物質に対して上述の仮定が正しいか、は物理として難しい問題であるが、実際に Heisenberg 模型で記述されると考えられている物質が存在する。物理学における文献については [SRFB] を参照。

# 3　Bethe 仮説法

相互作用する多数の量子力学的な粒子が存在するような系に対し Schrödinger 方程式 (1) を厳密に解くことは一般には不可能である。しかしながら一次元系に限れば、物理的に大変興味深い系も含め、ある程度系統的に解を求める方法が存在する。この様な方法は Heisenberg 模型 (2) に対して Bethe[6]によって 1931 年に導入された方法 [B] に基づく。この方法に対してつけられた Bethe 仮説法という名称は、原論文 [B] において波動関数が導入される際に "Wir machen den Ansatz" と述べられている事に基づき C. N. Yang[7] ら [YY] によって命名されたのが始まりのようである。

なお、Mermin–Wagner 理論 [MM] によれば一次元系では量子揺らぎが強すぎて相転移が起こらないとされている。よって二次元以上の系が解けるほうが興味深いのではあるが、現状では散在的な結果 [K06] は知られているものの一般性のある方法は絶望的なようである。

Bethe 仮説法はその後多くの拡張等が行われ、膨大な数の論文が出版されているが、本稿では Bethe の原論文と全く同じ設定の下で解説する。具体的には Hamiltonian (2) を周期境界条件 $\sigma_{N+1}^a = \sigma_1^a$ の下で考察することとする。

重要な性質として、Hamiltonian $\mathcal{H}_N$ は状態空間 $\mathfrak{H}_N$ 上の $sl_2$ 作用と可換になる。従って解の構成では $sl_2$ 作用に対する最高ウェイトベクトルのみを具体的に構成できれば、後は $sl_2$ の下降演算子を用いて残りの波動関数が求まることになる。結果として最高ウェイトベクトルの具体的構成という数学的にも興味深い問題を考察することになる。

---

[5] Paul Adrien Maurice Dirac (1902–1984). 量子力学の開拓者の一人で、電子の相対論的な波動方程式の導出などの業績で 1933 年ノーベル物理学賞受賞。解析力学に関する造詣が深く、有名な量子力学における経路積分法の導入 (1933 年) も正準変換の深い考察に基づく。

[6] Hans Albrecht Bethe (1906–2005). 特に量子力学の各方面への応用について多くの業績を挙げた。本稿で扱う Bethe 仮説法も同様であり、彼の最も有名な業績の一つ。1967 年ノーベル物理学賞受賞。

[7] Chen-Ning Yang (b. 1922). 素粒子標準理論の屋台骨を成す非可換ゲージ理論—Yang–Mills 理論の導入 (1954 年) や、量子可積分系で重要な Yang–Baxter 方程式などで有名。パリティの非保存（自然界では右と左は等価ではない）の発見により 1957 年ノーベル物理学賞受賞。彼の業績を調べてみると、Yang 先生は数理物理学者であると言って良いのではないかと思われるほどである（因みに最初の論文は Bull. Amer. Math. Soc. から出版）。実際 Faddeev 氏のエッセイ [F00] によれば、Yang–Mills 理論の発見において数学的直観が重要であったとご自身がよく述べておられたそうである。20 世紀最大級の物理学者の一人といって良いであろう。

**代数的 Bethe 仮説法**　以下では Faddeev[8] らによる定式化に基づいて天下り的に紹介する。詳細は [F96, TS, LR] などを参照して頂きたい。

代数的 Bethe 仮説法では、元々の状態空間 $\mathfrak{H}_N$ の代わりに空間 $\mathbb{C}^2 \otimes \mathfrak{H}_N$ を考える。ここで左側につけた $\mathbb{C}^2$ は補助空間と呼ばれ、特に強調したいときには 0 番目の成分と考え $\mathbb{C}^2_0$ と書く。以下の議論では $\mathbb{C}^2 \otimes \mathfrak{H}_N$ 上の作用素を $\mathbb{C}^2_0$ に関する 2 行 2 列の行列の形で表すと便利な事が多い。その場合各成分は $\mathfrak{H}_N$ 上の作用素 (つまり $2^N$ 次元行列) となる。

簡単な例で準備体操しておこう。二つの行列

$$A = \begin{pmatrix} a_{11} & a_{12} \\ a_{21} & a_{22} \end{pmatrix}, \qquad B = \begin{pmatrix} b_{11} & b_{12} \\ b_{21} & b_{22} \end{pmatrix},$$

を考える。するといわゆる Kronecker 積 $A \otimes B$ は $\mathbb{C}^2 \otimes \mathfrak{H}_1 = (\mathbb{C}^2)^{\otimes 2}$ に作用する作用素となる。具体的には

$$A \otimes B = \begin{pmatrix} a_{11}B & a_{12}B \\ a_{21}B & a_{22}B \end{pmatrix} = \begin{pmatrix} a_{11}b_{11} & a_{11}b_{12} & a_{12}b_{11} & a_{12}b_{12} \\ a_{11}b_{21} & a_{11}b_{22} & a_{12}b_{21} & a_{12}b_{22} \\ a_{21}b_{11} & a_{21}b_{12} & a_{12}b_{11} & a_{22}b_{12} \\ a_{21}b_{21} & a_{21}b_{22} & a_{12}b_{21} & a_{22}b_{22} \end{pmatrix}.$$

で与えられる。この式の真ん中のような表式を「$\mathbb{C}^2_0$ に関する 2 行 2 列の行列」と呼んでいる。この Kronecker 積の定義では $(\mathbb{C}^2)^{\otimes 2}$ の基底を

$$\{v_+ \otimes v_+,\ v_+ \otimes v_-,\ v_- \otimes v_+,\ v_- \otimes v_-\}.$$

と選んでいる事になる。重要な性質として $(C \otimes D)(A \otimes B) = CA \otimes DB$ が成り立っている。実際、

$$\begin{aligned} (C \otimes D)(A \otimes B) &= \begin{pmatrix} c_{11}D & c_{12}D \\ c_{21}D & c_{22}D \end{pmatrix} \begin{pmatrix} a_{11}B & a_{12}B \\ a_{21}B & a_{22}B \end{pmatrix} \\ &= \begin{pmatrix} (c_{11}a_{11} + c_{12}a_{21})DB & (c_{11}a_{12} + c_{12}a_{22})DB \\ (c_{21}a_{11} + c_{22}a_{21})DB & (c_{21}a_{12} + c_{22}a_{22})DB \end{pmatrix} = CA \otimes DB. \end{aligned}$$

参考までに $(\mathbb{C}^2)^{\otimes 3}$ では基底を

$$\begin{aligned} \{&v_+ \otimes v_+ \otimes v_+,\ v_+ \otimes v_+ \otimes v_-,\ v_+ \otimes v_- \otimes v_+,\ v_+ \otimes v_- \otimes v_-, \\ &v_- \otimes v_+ \otimes v_+,\ v_- \otimes v_+ \otimes v_-,\ v_- \otimes v_- \otimes v_+,\ v_- \otimes v_- \otimes v_-\}. \end{aligned}$$

と選ぶ。

さて、本題に戻り、$\mathbb{C}^2 \otimes \mathfrak{H}_N$ 上に以下の Lax[9] 作用素を定義する（右側の等号は直接計算で示される)。

$$L_k(\lambda) = \begin{pmatrix} \lambda \mathbb{I}_N + \frac{i}{2}\sigma_k^3 & \frac{i}{2}\sigma_k^- \\ \frac{i}{2}\sigma_k^+ & \lambda \mathbb{I}_N - \frac{i}{2}\sigma_k^3 \end{pmatrix} = \lambda I \otimes \mathbb{I}_N + \frac{i}{2}\sum_a^3 \sigma^a \otimes \sigma_k^a. \tag{3}$$

[8] Ludvig Faddeev (b. 1934). 非可換ゲージ理論における Feynman 則の導出や可積分系の研究などで著名な数理物理学者。

[9] Peter David Lax (b. 1926). 言わずと知れた解析学の大家だが、可積分系に関する業績を数多く残している。2005 年アーベル賞受賞。

ここで、

$$\sigma_k^+ = \sigma_k^1 + i\sigma_k^2 = 2\begin{pmatrix} 0 & 1 \\ 0 & 0 \end{pmatrix}_k, \quad \sigma_k^- = \sigma_k^1 - i\sigma_k^2 = 2\begin{pmatrix} 0 & 0 \\ 1 & 0 \end{pmatrix}_k$$

であり、$\lambda$ はある複素パラメータ（スペクトルパラメータと呼ばれる）である。$L_k(\lambda)$ は $\mathbb{C}_0^2 \otimes \mathfrak{H}_N$ の 0 番目と $k$ 番目の $\mathbb{C}^2$ にのみ非自明に作用する。

代数的 Bethe 仮説法では、Lax 行列から以下のようにして定義される**転送行列**が重要である。

$$T_N(\lambda) = L_N(\lambda)L_{N-1}(\lambda)\cdots L_1(\lambda). \tag{4}$$

転送行列を $\mathbb{C}_0^2$ に関する 2 行 2 列の行列で表示したときの各成分には通常以下の名前が付けられる。

$$T_N(\lambda) = \begin{pmatrix} A_N(\lambda) & B_N(\lambda) \\ C_N(\lambda) & D_N(\lambda) \end{pmatrix}.$$

以上の定式化と、もともと考えていた Heisenberg 模型との関係は以下の重要な定理によって理解される。

**定理**

$\tau_N(\lambda) = \mathrm{tr}\,|_{\mathbb{C}_0^2} T_N(\lambda) = A_N(\lambda) + D_N(\lambda)$ とするとき

$$\mathcal{H}_N = \frac{iJ}{2}\frac{d}{d\lambda}\log\tau_N(\lambda)\Big|_{\lambda=\frac{i}{2}} - \frac{NJ}{2}\mathbb{I}_N.$$

この定理の証明の鍵は $L_k(\frac{i}{2}) = iP_{0,k}$、ここで $P_{0,k}$ は $\mathbb{C}_0^2 \otimes \mathfrak{H}_N$ 上 0 番目と $k$ 番目を置換する作用素、となることに基づく。

以上まとめると、Hamiltonian $\mathcal{H}_N$ の対角化（Schrödinger 方程式の解を求める事）は $\tau_N(\lambda)$ の対角化の問題に帰着される。この様な定式化をすると、$\tau_N(\lambda)$ をパラメータ $\lambda$ により展開することにより、Hamiltonian を含む可換な作用素の族を構成することができて、Heisenberg 模型の量子可積分性が証明できるというご利益がある。

**Bethe ベクトル**　以上の準備の下で固有ベクトルの構成をする。次の特別なベクトルから出発する。

$$|0\rangle_N = v_+ \otimes \cdots \otimes v_+ \in \mathfrak{H}_N.$$

その時 Bethe ベクトルは

$$\Psi_N(\lambda_1, \cdots, \lambda_\ell) := B_N(\lambda_1)\cdots B_N(\lambda_\ell)|0\rangle_N \in \mathfrak{H}_N$$

で定義される。ここで $B_N(\lambda)$ どうしの可換性（下記定理に続く脚注を参照）

$$[B_N(\lambda), B_N(\mu)] = 0$$

が示されるため、Bethe ベクトルを定義する際のパラメータ $\lambda_1, \cdots, \lambda_\ell$ の順序は重要ではない。この時基本的な性質は以下の通り。

定理

Bethe ベクトル $\Psi_N(\lambda_1,\cdots,\lambda_\ell)$ はパラメータが Bethe 仮説方程式

$$\left(\frac{\lambda_k+\frac{i}{2}}{\lambda_k-\frac{i}{2}}\right)^N=\prod_{\substack{j=1\\ j\neq k}}^{\ell}\frac{\lambda_k-\lambda_j+i}{\lambda_k-\lambda_j-i},\quad (k=1,\cdots,\ell) \tag{5}$$

を満たす時 $\tau_N(\lambda)$ の固有ベクトルとなる。

定理の証明では Yang-Baxter 関係式から従う関係式[10] を用いて直接作用

$$\begin{aligned}&\{A_N(\lambda)+D_N(\lambda)\}B_N(\lambda_1)\cdots B_N(\lambda_\ell)|0\rangle_N\\ &=\Lambda(\lambda;\lambda_1,\cdots,\lambda_\ell)\prod_{j=1}^{\ell}B(\lambda_j)|0\rangle_N+\sum_{k=1}^{\ell}\left\{\Lambda_k(\lambda;\lambda_1,\cdots,\lambda_\ell)B(\lambda)\prod_{\substack{j=1\\ j\neq k}}^{\ell}B(\lambda_j)|0\rangle_N\right\}\end{aligned}$$

を計算し[11]、右辺第二項が消える条件を書き下すと Bethe 仮説方程式が得られる。

重要な性質として、パラメータが Bethe 仮説方程式を満たす時 Bethe ベクトルは $\mathfrak{H}_N$ 上の $sl_2$ 作用に対する最高ウェイトベクトルとなる (表現の次元は $N-2\ell+1$) ことが示される。

以上の構成は既に数十年前までに十分確立された内容であり、数学的にも特に問題となる点はないはずである。

**長さ 3 の例**　少し具体的な例（$N=3$）を見てみることにしよう。直接計算（例えば次の項目参照）により $B_3(\lambda)$ は以下の形になる。

$$\begin{pmatrix}0&0&0&0&0&0&0&0\\ i\left(\lambda-\frac{i}{2}\right)^2&0&0&0&0&0&0&0\\ i\left(\lambda^2+\frac14\right)&0&0&0&0&0&0&0\\ 0&i\left(\lambda-\frac{i}{2}\right)^2&i\left(\lambda^2+\frac14\right)&0&0&0&0&0\\ i\left(\lambda+\frac{i}{2}\right)^2&0&0&0&0&0&0&0\\ 0&i\left(\lambda^2+\frac14\right)&-i&0&i\left(\lambda^2+\frac14\right)&0&0&0\\ 0&0&i\left(\lambda^2+\frac14\right)&0&i\left(\lambda+\frac{i}{2}\right)^2&0&0&0\\ 0&0&0&i\left(\lambda-\frac{i}{2}\right)^2&0&i\left(\lambda^2+\frac14\right)&i\left(\lambda+\frac{i}{2}\right)^2&0\end{pmatrix}.$$

表現論の知識からすると $\mathfrak{H}_3=(\mathbb{C}^2)^{\otimes 3}$ は $sl_2$ の 4 次元表現一つと 2 次元表現二つに分解するはずである。それぞれの表現を Bethe 仮説法により構成してみよう。まず $\ell=0$ の時は $|0\rangle_3=v_+\otimes v_+\otimes v_+$ が確かに 4 次元表現の最高ウェイトベクトルとなる。

---

[10] $A_N(\lambda)B_N(\mu)=\frac{\lambda-\mu-i}{\lambda-\mu}B_N(\mu)A_N(\lambda)+\frac{i}{\lambda-\mu}B_N(\lambda)A_N(\mu)$ および $D_N(\lambda)B_N(\mu)=\frac{\lambda-\mu+i}{\lambda-\mu}B_N(\mu)D_N(\lambda)-\frac{i}{\lambda-\mu}B_N(\lambda)D_N(\mu)$。この様な $A_N,B_N,C_N,D_N$ 作用素間の交換関係の導出の概略は、$R$-行列

$$R(\lambda)=\begin{pmatrix}1&0&0&0\\0&b(\lambda)&c(\lambda)&0\\0&c(\lambda)&b(\lambda)&0\\0&0&0&1\end{pmatrix}=\frac{1}{\lambda+i}\left(\left(\frac{\lambda}{2}+i\right)I\otimes I+\frac{\lambda}{2}\sum_{a=1}^{3}\sigma^a\otimes\sigma^a\right),\ b(\lambda)=\frac{i}{\lambda+i},\ c(\lambda)=\frac{\lambda}{\lambda+i}$$

を用いて $\mathbb{C}_0^2\otimes\mathbb{C}_0^2$ 上の関係式 $R(\lambda-\mu)\left(L_k(\lambda)\otimes L_k(\mu)\right)=\left(L_k(\mu)\otimes L_k(\lambda)\right)R(\lambda-\mu)$ を直接示し、その積を取って $R(\lambda-\mu)\left(T_N(\lambda)\otimes T_N(\mu)\right)=\left(T_N(\mu)\otimes T_N(\lambda)\right)R(\lambda-\mu)$ を示すことによる。

[11] 最後に関係式 $A_N(\lambda)|0\rangle_N=\left(\lambda+\frac{i}{2}\right)^N|0\rangle_N$ および $D_N(\lambda)|0\rangle_N=\left(\lambda-\frac{i}{2}\right)^N|0\rangle_N$ を用いる。これらは転送行列の定義から直接示される。

次に $\ell=1$ の時、Bethe 仮説方程式は、右辺が自明となって、

$$\left(\frac{\lambda_1+\frac{i}{2}}{\lambda_1-\frac{i}{2}}\right)^3=1,$$

であり、解は $\lambda_1=\pm\frac{1}{\sqrt{12}}$ である。そのとき、先ほどの $B_3(\lambda)$ を用いると、

$$B_3\left(\frac{1}{\sqrt{12}}\right)|0\rangle_3=\left(0,-\frac{i-\sqrt{3}}{6},\frac{i}{3},0,-\frac{i+\sqrt{3}}{6},0,0,0\right)^t$$
$$B_3\left(-\frac{1}{\sqrt{12}}\right)|0\rangle_3=\left(0,-\frac{i+\sqrt{3}}{6},\frac{i}{3},0,-\frac{i-\sqrt{3}}{6},0,0,0\right)^t$$

が得られる。これらは互いに異なる二つの 2 次元表現の最高ウェイトベクトルとなっている。

表現論からの結果を見れば、以上で全ての最高ウェイトベクトルを尽くしているはずである。しかし、興味深いので、更に計算を進めることにする。$\ell=2$ の場合、Bethe 仮説方程式は

$$\left(\frac{\lambda_1+\frac{i}{2}}{\lambda_1-\frac{i}{2}}\right)^3=\frac{\lambda_1-\lambda_2+i}{\lambda_1-\lambda_2-i},\quad \left(\frac{\lambda_2+\frac{i}{2}}{\lambda_2-\frac{i}{2}}\right)^3=\frac{\lambda_2-\lambda_1+i}{\lambda_2-\lambda_1-i}$$

となる（ようやく右辺が非自明になる）。この解は

$$\{\lambda_1,\lambda_2\}=\{0,0\},\{\frac{\sqrt{3}}{2},\frac{\sqrt{3}}{2}\},\{-\frac{\sqrt{3}}{2},-\frac{\sqrt{3}}{2}\},\{\frac{i}{2},-\frac{i}{2}\},\{-\frac{i}{2},\frac{i}{2}\}$$

と求まるが、直接計算すると $B_3(\lambda_1)B_3(\lambda_2)=0$ となることが確かめられる。更に $\ell=3$ の場合には、1 次元の解が 4 つと 0 次元の解が 12 個求まりそのすべてについて $B_3(\lambda_1)B_3(\lambda_2)B(\lambda_3)=0$ が確かめられる。

**Mathematica による実装**　Mathematica をご利用にならない方でも、実際の計算がどのような調子であるのかを知ることは有益であると思うので、ここにコードの例を与える。以下本節では “n” はチェインの長さ $N$ を表し、“l” はスペクトルパラメータ $\lambda$ を表すこととする。なお以下の内容は “Mathematica Summer School on Theoretical Physics—4th edition (2012) Integrability and Super Yang–Mills” における下記の内容

http://msstp.org/sites/default/files/Exercise_byJoao.nb

を参考にして、現在の定義と合うように改変した。

まず Pauli 行列

```
s=PauliMatrix
```

を定義する。この定義より Pauli 行列は `s[a]` $=\sigma^a$ （左辺の様に Mathematica に入力する）として得られる。次に $2^n$ 次元の単位行列

```
id[n_]:=IdentityMatrix[2^n, SparseArray]
```

と定義する。一旦 “SparseArray” と宣言しておくと、以下の計算結果は全て疎な配列として扱われ、大幅に高速化する。そのとき Hamiltonian (2) は

```
H[n_,J_]:=(J/4)(Sum[KroneckerProduct[s[a],id[n-2],s[a]],{a,3}]
   +Sum[KroneckerProduct[id[k-1],s[a],s[a],id[n-k-1]],{k,n-1},{a,3}]
   -n*id[n])
```

となる。ただし和は $k = N$ と $k = 1, 2, \dots, N-1$ の二つの部分に分けてある。式 (3) 右端の表示式による Lax 作用素 $L_k(\lambda)$ は

```
L[n_,k_,l_]:=l*id[n+1]+
   (I/2)Sum[KroneckerProduct[s[a],id[k-1],s[a],id[n-k]],{a,3}]
```

また式 (4) の転送行列は、積を取ることにより

```
tm[n_,l_]:=Dot@@Table[L[n,k,l],{k,n,1,-1}]
```

作用素 $B_N(\lambda)$ は、転送行列の部分集合を取ることにより、

```
b[n_,l_]:=tm[n,l][[1;;2^n,2^n+1;;2^(n+1)]]
```

最後にベクトル $|0\rangle_N$ は

```
hw[n_]:=UnitVector[2^n,1]
```

例えばベクトル $B_4(l_1)B_4(l_2)|0\rangle_4$ は

```
Factor[b[4,l1].b[4,l2].hw[4]]
%//MatrixForm
```

として得られる。

# 4 Bethe 仮説方程式を解く

ここまでの議論より Bethe 仮説方程式の解が求まれば固有ベクトルが求まることが分かった。しかし Bethe 仮説方程式は複雑な方程式であるし、むしろ量子多体系の解の構成という難しい問題の重要な一部分を一組の代数方程式の形に昇華したかのような印象を持つ。従って Bethe 仮説方程式そのものの理解を深めることは本質的に重要な課題であると考えている。

もちろん Bethe 仮説方程式そのものの研究も行われてきた。例えば Langlands[12]–Saint-Aubin [LS94, LS97] の研究は最近 ASEP と呼ばれる系の研究に応用された [BDS]。しかしながら Bethe 仮説方程式にはなお研究すべき余地が多く存在すると考えられる。本稿では Bethe 仮説方程式の研究で古くから用いられている艤装配位（rigged configuration）と呼ばれる対象を軸として解説する。

艤装配位はもともと Bethe の原論文（1931 年）において Bethe 仮説方程式の解の性質を研究する目的で導入された組み合わせ論的対象である。その後 1986 年に Kerov, Kirillov, Reshetikhin ら [KKR, KR] によって大幅に拡張され、また艤装配位の組み合わせ論的研究が開始された。詳細は §8 で述べるが、後者の研究では Bethe 仮説方程式と Bethe ベクトルの対応の組み合わせ論的類似を構成することに主眼があり、それ自身極めて深い数学的内容を持つことが明らかになっている。

[12] Robert Pbelan Langlands (b. 1936). 整数論における業績とは別に数理物理学でも多くの業績を挙げている。ここでふれた Bethe 仮説に関する研究のほか Virasoro 代数や共形場の理論の研究なども有名である。

**艤装配位の定義**　ここでは現在の設定を若干拡張して $sl_2$ の同じ次元（2 次元とは限らない）の表現がテンソル積で並んでいる状況で定義を解説する。必要なデータとして二つの分割 $\mu$ と $\nu$ を準備する。以下では整数の分割と Young 図は自由に同一視して考えることとする。一つ目のデータ $\mu$ は状態空間の形状を表す。例えば本稿で主に扱っている状態空間 $\mathfrak{H}_N = (\mathbb{C}^2)^{\otimes N}$ では $\mu = (1^N)$ とする。一方状態空間 $(\mathbb{C}^3)^{\otimes N}$ に対する一般化された Heisenberg 模型の場合であれば $\mu = (2^N)$ とする、といった具合である。

更に進むためにいくつか重要な用語を定義する。分割 $\lambda$ に対し $Q_k(\lambda)$ を Young 図形 $\lambda$ の左側 $k$ 列分のます目の数とする。その時 **vacancy number** を

$$P_k(\mu,\nu) = Q_k(\mu) - 2Q_k(\nu)$$

で定義する（係数の 1, $-2$ という数列は実は $A$ 型の Cartan データの一部である）。その時艤装配位 $(\mu,(\nu,J))$ は次のように定義される（以下しばしば $\mu$ を略記して $(\nu,J)$ と表す）。$(\nu,J)$ を具体的に書くと、Young 図 $\nu$ の各行 $\nu_i$ に対して整数 $J_i$ を対応させた

$$(\nu,J) = \{(\nu_1,J_1),(\nu_2,J_2),\ldots,(\nu_l,J_l)\}$$

なる集合であり、ここに以下の二つの条件を課す（整数 $J_i$ を **rigging** と呼ぶ）。

- 分割 $\nu$ の各行は $0 \le P_{\nu_i}(\mu,\nu)$ を満たす。
- 整数 $J_i$ は $0 \le J_i \le P_{\nu_i}(\mu,\nu)$ を満たす。

この時 $(\nu,J)$ を（最高ウェイトベクトルの場合の）艤装配位と呼ぶ。以下に述べる理由により組 $(\nu_i,J_i)$ の事をしばしば**ストリング**と呼ぶ。なお、艤装配位では同じ長さの複数の行に対する riggings の順番は重要ではない。

様子を見るために簡単な例を考えてみよう。$sl_2$ の 3 次元表現が二つ並んでいる $(\mathbb{C}^3)^{\otimes 2}$ の場合、$\mu = (2,2)$ となる。その時以下の三つの艤装配位が存在する。

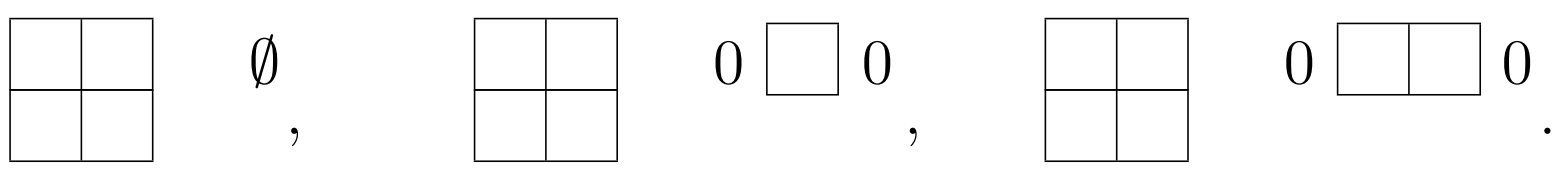

上図では左側に $\mu$ を、右側に $\nu$ を配置し、$\nu$ の各行の左側には対応する vacancy number を、右側には rigging を書いた。右側二つの艤装配位の場合 vacancy numbers が 0 なので riggings としては 0 のみが許される。例えば右端の艤装配位の場合 $\nu_1 = 2$ なので、対応する vacancy number は $P_2((2,2),(2)) = 4 - 2\cdot 2 = 0$ となる。一方 $\nu = (1,1)$ とすると、対応する vacancy number は $P_1((2,2),(1,1)) = 2 - 2\cdot 2 < 0$ となって艤装配位の条件を満たさない。

**Bethe 仮説法との関係**　では、この様にして定義される艤装配位と Bethe 仮説法との関係は何であろうか。それは、もともとは Bethe 自身が発見したように、Bethe 仮説方程式の解には非常に特徴的なパターンがみられるという観察である。雰囲気をつかむために、Bethe 仮説方程式の解の典型的な例を複素平面上で図示したものを示そう（$N = 12$, $\ell = 6$ の例）。

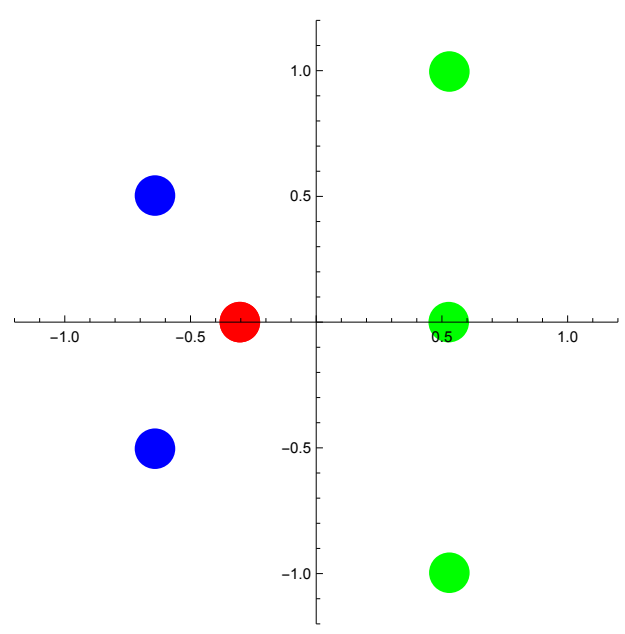

この例で観察されるように、根はいくつかのグループ（ストリングと呼ぶ）に分けられ、それぞれのストリングに属する根は、実部が「大体」同じであり、根どうしの間隔は「大体」$\sqrt{-1}$ である、というものである。そのような観察に基づき、ストリングの長さが $n$（$n$-ストリングと呼ぼう）であれば Young 図の長さ $n$ の行を対応させたい。例えば上で掲げた例であれば 3-ストリング（緑）、2-ストリング（青）、1-ストリング（赤）の三種類がそれぞれ一つずつあるので、以下の Young 図と対応させたい。

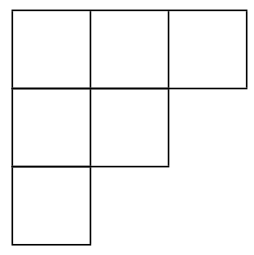

もちろんここでの記述から明らかなように、以上の観察はかなり大雑把なものである。その点をどのように科学的にするのかは §6 で改めて議論したい。しかしながら、結果的には Bethe 仮説方程式の "物理的な解" は rigged configurations と一対一に対応するであろうと考えている。その際 rigging は各ストリングの実部の位置を表すこととなる。

**$N = 12$ の場合の実例**　艤装配位と Bethe 仮説方程式の根の対応について具体的な状況を見るために長さ $N = 12$ で $\ell = 5$ の場合の例を考えてみよう。これらは次節で述べる「物理的特異解」となっている（実は以下の 5 つの解で $N = 12$、$\ell = 5$ の場合の物理的特異解を全て尽くしている—詳細は [Sa15, Appendix B] 参照）。この場合の Bethe 仮説方程式

$$\left(\frac{\lambda_k + \frac{i}{2}}{\lambda_k - \frac{i}{2}}\right)^{12} = \prod_{\substack{j=1 \\ j \neq k}}^{5} \frac{\lambda_k - \lambda_j + i}{\lambda_k - \lambda_j - i}, \quad (k = 1, \cdots, 5)$$

において $\lambda_1 = 0, \lambda_2 = i/2, \lambda_3 = -i/2$ となる解を求めると、重根となる場合を除いて $\lambda_4$ および $\lambda_5$ は以下の 5 次方程式

$$5120\xi^5 + 11520\xi^4 - 4992\xi^3 - 9312\xi^2 + 2020\xi - 55 = 0 \tag{6}$$

の五つの実数解から $\lambda_4 = \sqrt{\xi}$ および $\lambda_5 = -\sqrt{\xi}$ として得られる。方程式 (6) は、少なくとも Mathematica では厳密解が求められず、数値解としては

$$\xi_1 = -2.29679,\ \xi_2 = -0.999662,\ \xi_3 = 0.0320332,\ \xi_4 = 0.173735,\ \xi_5 = 0.840679$$

が得られ、それぞれの平方根は $\pm 1.51551i,\ \pm 0.999831i,\ \pm 0.178978,\ \pm 0.416816,\ \pm 0.916886$ となる。こうして得られる五つの解を複素平面上に図示すると図 1 のようになる。黄色のドットが $\lambda_4$ および $\lambda_5$ を表す。対応する艤装配位は順に以下の通り（riggings は §6 で説明する方法で決定した）。

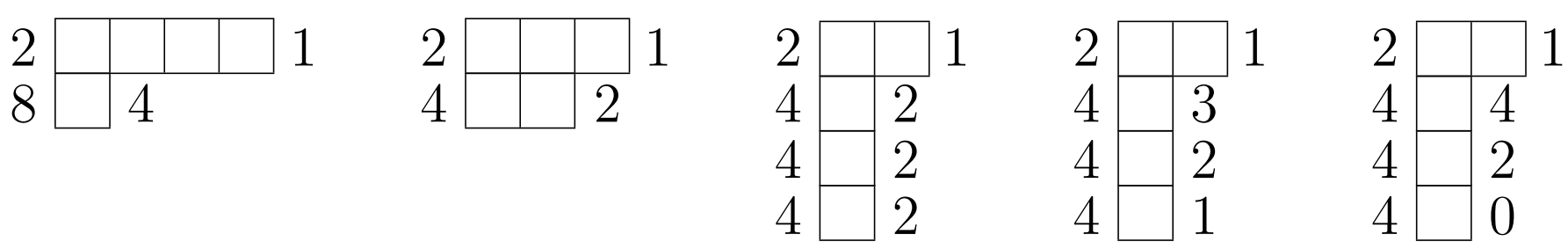

こうしてみると方程式 (6) の根が実に絶妙な位置に存在していることが分かる。

なお文献 [KS15] の補助ファイルとして $N = 12$ の Bethe 仮説方程式の根を複素平面上にプロットしたものが多数含まれているのでご利用頂きたい。

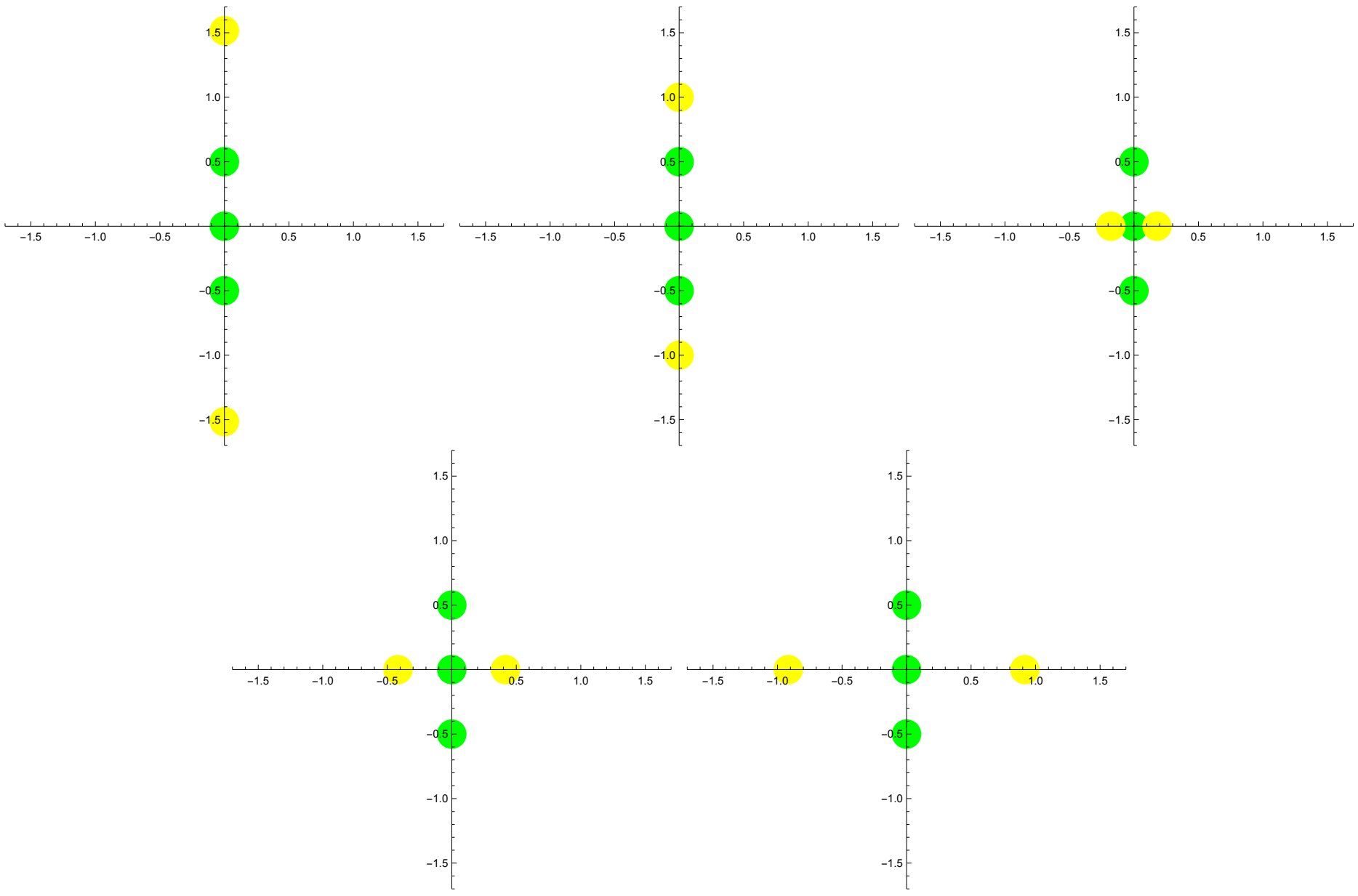

図 1: $N=12,\ \ell=5$ の物理的特異解

## 5 Bethe 仮説方程式の特異解

更に進んで Bethe 仮説方程式の解を解析しようとすると、いわゆる特異解の問題を避けて通ることはできない。今までの議論をおさらいすると、我々に関心のある Schrödinger 方程式とは固有値問題のことであり、ここまでに構成した固有ベクトルについてエネルギー固有値を計算すると、結局以下のような解を得たことになる（本質的には Bethe [B] による結果）。

$$\mathcal{H}_N \Psi_N(\lambda_1,\ldots,\lambda_\ell) = \mathcal{E}_{\lambda_1,\ldots,\lambda_\ell} \Psi_N(\lambda_1,\ldots,\lambda_\ell),$$
$$\mathcal{E}_{\lambda_1,\ldots,\lambda_\ell} := -\frac{J}{2}\sum_{j=1}^{\ell} \frac{1}{\lambda_j^2+\frac{1}{4}}.$$

さてこの結果を観察すると、Bethe 仮説方程式の解 $\{\lambda_1,\ldots,\lambda_\ell\}$ の中に $i/2$ または $-i/2$ が含まれると固有値 $\mathcal{E}_{\lambda_1,\ldots,\lambda_\ell}$ が発散することが分かる。そこでこの様な解を**特異解**とよぶ。

以下では Bethe 仮説方程式の解 $\{\lambda_1,\ldots,\lambda_\ell\}$ が互いに相異なり、かつ $i/2$ および $-i/2$ を含まないとき**正則な解**と呼ぶことにする。ここで厄介になるのは、正則な解のみを考えるのでは固有ベクトルが不足し（わずか $N=4$ の場合からこの問題が発生する）、一方 Bethe 仮説方程式の互いに相異なる解を全て持ってきてしまうと固有ベクトルの数を超過してしまうことである。従って特異解のうち物理的な意味づけを持つ適切な部分集合（**物理的特異解**と呼ぶ）を選択する必要がある。

参考までに、Bethe 仮説方程式の相異なる解の総数と物理的特異解の総数は二項係数を用いて記述されると予想されている [KS14a]。例えば与えられた共に偶数の $N$ および $\ell$ に対して物理的特異解の総数は

$$\binom{\frac{N-2}{2}}{\frac{\ell-2}{2}}$$

であると予想されている（関連する情報について、[DG14] および [Sa15, Appendix B] も参照)。

先に進む前に、一旦議論を中断して、Bethe 仮説について証明はまだ出来ていないのではないかと思われるが、おそらく正しいと思われる性質をまとめておくことにしよう。

- 正則な解と以下で述べる物理的特異解から構成される Bethe ベクトルは 1 次独立であると思われる。
- Bethe 仮説方程式の正則な解以外に対する Bethe ベクトルは 0 になると思われる。例えば 1 次元以上の解も存在するが、具体例で計算すると対応する Bethe ベクトルは 0 になる。
- 現在の $\mathfrak{H}_N = (\mathbb{C}^2)^{\otimes N}$ の場合は Bethe 仮説方程式の互いに相異なる解だけを考えれば十分であると考えられているが、より一般の場合には重根も考察しなければならないことが数値計算の結果により示唆されている [HNS14]。

さて特異解の問題に戻ることにしよう。Bethe 仮説法の思想に基づけば、ある条件を満たす特異解から出発して、何らかの正則化の手続きに基づき Bethe ベクトルを構成できるような方法を構築することが望ましい。この様な方向の研究は、最初 Avdeev と Vladimirov [AV] によって $N = 4$ の場合に提唱された（1986 年)。彼らの提案した方法はその後 Beisert ら [BMSZ] によってより一般的な例に対して拡張され（2003 年)、最終的に Nepomechie と Wang [NW] によって一般的な状況で計算された（2013 年)。ここでは [NW] に従って最終的な結果を紹介することとする（証明については、例えば [Sa15, Appendix A] を参照)。

一般的な特異解を

$$\left\{\frac{i}{2}, -\frac{i}{2}, \lambda_3, \ldots, \lambda_\ell\right\} \tag{7}$$

と表すことにする。この特異解が以下の Nepomechie–Wang の条件

$$\left(-\prod_{j=3}^{\ell} \frac{\lambda_j + \frac{i}{2}}{\lambda_j - \frac{i}{2}}\right)^N = 1, \tag{8}$$

を満たす場合**物理的特異解**とよばれ、次のような方法で対応する Bethe ベクトルを構成できる。定数 $c$ を

$$c = 2i^{N+1} \prod_{j=3}^{\ell} \frac{\lambda_j + \frac{3i}{2}}{\lambda_j - \frac{i}{2}}. \tag{9}$$

と定める。

Nepomechie–Wang の定理

条件 (8) を満たす物理的特異解 (7) に対し次の極限は 0 でないベクトルに収束する [NW]。

$$\lim_{\epsilon\to 0}\frac{1}{\epsilon^N}B_N\left(\frac{i}{2}+\epsilon+c\,\epsilon^N\right)B_N\left(-\frac{i}{2}+\epsilon\right)B_N(\lambda_3)\cdots B_N(\lambda_\ell)|0\rangle_N.$$

このベクトルは固有値

$$\mathcal{E}_{\frac{i}{2},-\frac{i}{2},\lambda_3,\ldots,\lambda_\ell}=-J-\frac{J}{2}\sum_{j=3}^{\ell}\frac{1}{\lambda_j^2+\frac{1}{4}}.$$

に対する固有ベクトルとなる [KS14b]。

Hao–Nepomehie–Sommese [HNS13] は Bethe 仮説方程式に対する大規模な数値計算を行い、Bethe 仮説方程式の重根を含まない解の総数および正則な解と物理的特異解の総数を求めた。その結果に基づき彼らは以下の予想を定式化した。

$$\begin{aligned}&(\mathfrak{H}_N\text{の最高ウェイトベクトルの総数})\\ &=(\text{正則な解の個数})+(\text{物理的特異解の個数}).\end{aligned}$$

この予想は同論文において $N=14$ の場合まで数値的に確認された。この予想と、前述した対応する Bethe ベクトルが一次独立であるという予想が共に証明されれば、Bethe 仮説法によって全ての固有ベクトルが求められ、Heisenberg 模型に対する Schrödinger 方程式の完全な解が得られたことになる。しかしながら、筆者の知る限り、これら二つの予想は未解決であると思われる。

艤装配位を用いて Bethe 仮説方程式の解を解析すると、物理的特異解についてより立ち入った情報を得ることができる。艤装配位の上にフリップ

$$\kappa:(\nu_i,J_i)\longmapsto(\nu_i,P_{\nu_i}(\mu,\nu)-J_i) \tag{10}$$

と呼ばれる involution を定義する（全ての rigging を対応する vacancy number と rigging の差に置き換える）。§6 で見るように、この操作は解の上では、全ての根を $(-1)$ 倍する

$$\{\lambda_1,\ldots,\lambda_\ell\}\to\{-\lambda_1,\ldots,-\lambda_\ell\}$$

という操作に対応する。

さて、論文 [KS14a] において、Bethe 仮説方程式の解を艤装配位と対応付けるとき、系の長さ $N$ が偶数の場合の物理的特異解は以下の性質を持つ艤装配位 $(\nu,J)$ に対応する解として特徴づけられる事が予想された（[Sa15, Apendix B] も参照）。

(i) フリップ不変であり、かつ

(ii) 長さ $\nu_i$ が偶数の行の総数が奇数である艤装配位。

すなわち該当する艤装配位を書き下して物理的特異解の大体の形を決定することができる。長さ $N=12$ かつ $\ell=5$ の場合の例が §4 に与えられている。この観察は艤装配位を使った Bethe 仮説方程式の解の解析の自然さを表す例であると考えている。

なお、上記特徴づけでは長さ $N$ が奇数の場合には一見物理的特異解が存在しないように見えるが、数は少ないながらも散在的に物理的特異解が見出される場合がある。この点については [KS14a, §4.1] 等を参照されたい。

# 6　Betheの量子数

Bethe 仮説方程式の解を艤装配位 $(\nu, J)$ を用いて解析すると、同じ形 $\nu$ に属するいくつかの解に正しくラベル $J$ を割り当てる必要がある。量子力学では、固有状態に対する「良い」ラベル（厳密な定義があるわけではない）のことを通常量子数と呼ぶ。従って個々の解に対して rigging $J$ という量子数を定義することが問題となる。

Heisenberg 模型に対しては、Takahashi の量子数 ([T]、1971 年) と呼ばれるものがあり、近似的には大変うまく行くことで有名であった。例えば熱力学的量の導出など、物理として意義深い帰結をもたらした。Takahashi の量子数に基づく解析は、文献ではしばしば「ストリング仮説」などとして言及される。

ただし Takahashi の量子数の定義では、前提として各ストリングが厳密に

$$a + bi,\ a + (b-1)i,\ a + (b-2)i,\ \ldots,\ a - bi, \qquad (a \in \mathbb{R},\ b \in \mathbb{Z}_{\geq 0}/2)$$

という形をとると仮定している。しかしながら、具体例を検討してみれば、長さ 3 以上のストリングであればほとんど常に各根は異なる実部を持っているし、また隣り合う根どうしの間隔も $\sqrt{-1}$ とは異なる値を持つ。従って Takahashi の量子数を具体的な解について計算しようとすれば通常一意的には定まらないし、また数値も整数（または半整数）ではない複雑な数となってしまう。このようないわゆるストリング仮説をめぐる問題は Bethe 仮説の研究において少なからぬ混乱をもたらしたようであり、しばしばストリング仮説の破れに関する議論を耳にしたり、場合によっては Bethe 仮説そのものの数学的正当性を疑うかのような見解を耳にする事すらあった。

筆者としては、この様な混乱は、Heisenberg 模型や Bethe 仮説そのものの欠陥によるものではなく、むしろそれらは数学的に正当化されるべきものであると考えたい。そのため、以下では Kirillov 氏と共同で行った研究 [KS15] に基づき、数学的に正当化可能であると考えられる方法で量子数を定義する提案について解説したい。

出発点として Bethe が原論文 [B] で述べている次のような考察から始めよう（詳細は [KS15, §3] を参照)。Arctan $z$ を実軸上の arctan $x$ ($\arctan 0 = 0$) を複素平面上 $(i, +i\infty)$ と $(-i, -i\infty)$ に branch cuts を入れて一価に拡張した関数とする。そのとき関係式

$$\log \frac{z-i}{z+i} = 2i\,\mathrm{Arctan}(z) + (2n+1)\pi i \qquad (n \in \mathbb{Z}). \tag{11}$$

が、両辺を微分することにより示される。ここで左辺の log は多価性を持ったフルバージョンの対数関数である。この式を用いて Bethe 仮説方程式（の右辺を左辺に移したもの）の対数

$$\log\left[\left(\frac{\lambda_k + \frac{i}{2}}{\lambda_k - \frac{i}{2}}\right)^N \prod_{\substack{j=1 \\ j\neq k}}^{\ell} \frac{\lambda_k - \lambda_j - i}{\lambda_k - \lambda_j + i}\right], \quad (k = 1, \ldots, \ell)$$

を計算してみよう。この式は Bethe 仮説方程式により $\log 1$ に等しいはずであるが、ここから複素平面上の位相因子に由来する整数が得られる。具体的に計算を実行してみれば、次のような整数 (半整数)

$$J_k = \frac{N}{2\pi}\left(2\,\mathrm{Arctan}(2\lambda_k) - \frac{2}{N}\sum_{\substack{j=1 \\ j\neq k}}^{\ell} \mathrm{Arctan}(\lambda_k - \lambda_j)\right) \tag{12}$$

を導出することができる。この数をBetheの量子数と呼ぶことにする。不思議なことにBetheの量子数はBetheの原論文に登場しているにもかかわらず、その後80年以上にわたり十分な注目を浴びてこなかったようである。

以下で例を見るようにBetheの量子数はそのままでは扱いづらい。そこで、新しい量子数として、$S$を一つのストリングとするときそれに属する各根に対するBetheの量子数の和

$$J(S) = \sum_{\eta \in S} J_\eta \tag{13}$$

を考察することが重要となる。この様に定義すると本来は“ストリング”という概念そのものに由来するあいまいさが存在する。例えば隣り合うストリング同士が混ざってしまって分離できない場合などがあり、その様なあいまいさの処理方法はよく分かっていなかった。しかしこの定義(13)を$N = 12$の全ての解について詳しく検討すると、様々なタイプの不定性に対して、そのあいまいさの処理の仕方によらず一意的な結論が得られる。これは状態空間の次元が$2^{12} = 4096$であるからかなり非自明な例である。

重要な観察として、$J(S)$は極めて安定なふるまいをする。例として、$N = 12$、$\ell = 6$で5-stringがひとつ、1-stringがひとつ含まれるような場合を考える。対応するrigged configurationsは、$r = 0, 1, \ldots, 8$に対して、次の9つ

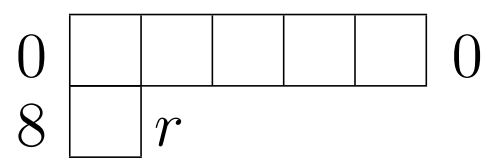

が存在する。以下それらの解を複素平面上にプロットしたものを掲げるが、経験則として解は最も大きなストリング$S$の実部$\Re(S)$によって並べると構造が分かりやすい。ここで$\Re(S)$とは$S$の両端の根の実部によって定義する（Vladimirovの結果[V]より両者は一致する）。

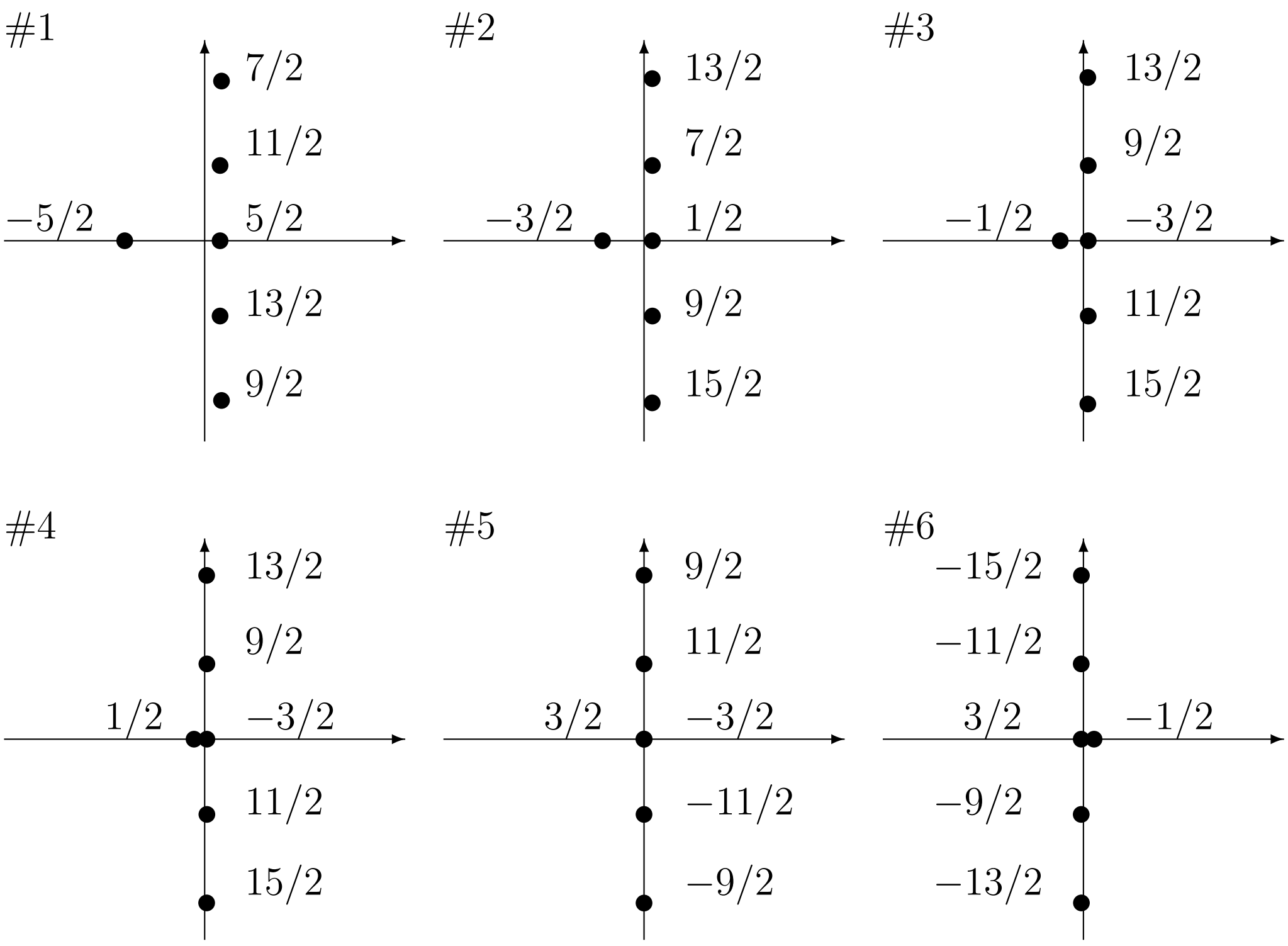

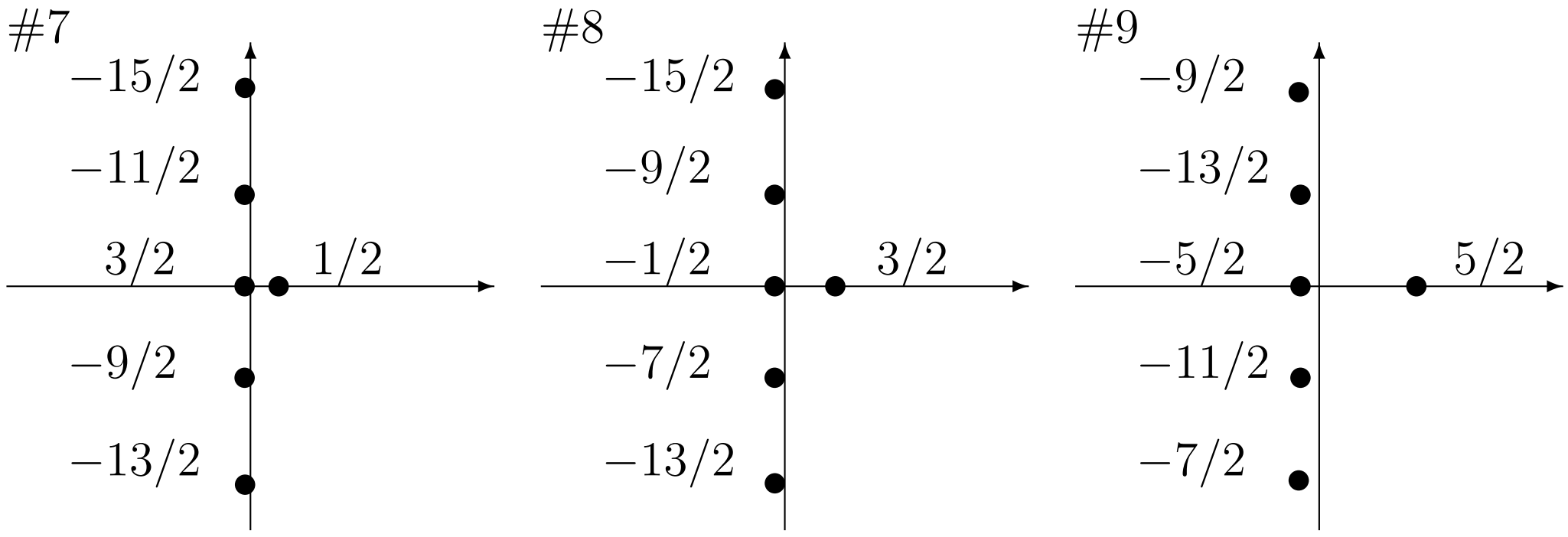

この例で分かるように個々のBetheの量子数は複雑な挙動を示す一方、(13) で定義される新しい量子数は5-string $S$ にたいして $\Re(S) > 0$ なら $J(S) = 45/2$、$\Re(S) < 0$ なら $J(S) = -45/2$ と非常に安定なふるまいをすることが分かる。

論文 [KS15] では、数値計算の結果 [GD, HNS13] に基づき、長さ $N = 12$ の場合には量子数 $J(S)$ によって全ての解が分類でき、対応する rigging を特定できることが確認されている。この際 $J(S)$ の具体的値まで特定できれば良いのだが、意外と複雑な挙動を示すため現状のデータでは推定できていない。しかし与えられた $\nu$ をもつ解全体を考え（もともと解全体が完全かどうか、などといった問題を考察しているので、解の全体が手元にあることを仮定するのは不都合ではないだろう）、その中で最大の $J(S)$ をもつ状態を探せば、後は帰納的に決定していくことができる。この際非自明な性質として、$k$-string と $l$-string との相対的な右左の位置関係が入れ替わるときに因子 $\Delta(k, l) = k + l - 3$ を $J(S)$ に加えるまたは減ずる必要がある様である。

今後の課題としては、まずBetheの量子数とは何なのかを理解すること、特に可能であればより内在的な定義を見つけられれば理解が進むのではないかと期待している。もちろん以上の観察自体も、非常に非自明な例について確認されているとはいえ、証明された結果ではないので今後様々な修正が必要となる可能性がある。

最後にストリングに関する他の話題として、$\ell = 2$ の場合に現れる例外的実数解（[EKS]、1992 年）と呼ばれる現象がある事に触れておこう。通常この現象は、存在するはずの 2-ストリングのうちいくつかが消え、その分実数解の個数が増える場合が存在する、といったように説明されることが多いようである。もう少し踏み込んで例外的実数解がどのように例外的なのかは好奇心をそそるところであるが、[Sa15] において具体例に対して実数解全体の集合を分析してみたところ、通常型の実数解は艤装配位の構造にきれいに従って出現し、例外型実数解は明瞭に艤装配位の構造を破っていることを確認した。量子数も問題なく定まる。従って例外型実数解の問題については必要以上に神経質になる事では無いのではないかという感触を持っている。

## 7　量子群の場合の文献紹介

以上通常の Lie 代数 $sl_2$ に関わる内容を紹介してきた。対称性が量子群となる場合にどのような現象が発生するのか興味をお持ちの読者もいらっしゃるかもしれないので、便宜のためにここで文献を紹介しておく。

量子群 $U_q(sl_2)$ の対称性を持つ模型が Pasquier–Saleur（[PS]、1990 年）によって考察され

ている。この模型の Hamiltonian は

$$H_N = \sum_{k=1}^{N-1} \left[ \sigma_k^1 \sigma_{k+1}^1 + \sigma_k^2 \sigma_{k+1}^2 + \frac{1}{2}(q+q^{-1}) \sigma_k^3 \sigma_{k+1}^3 \right] - \frac{1}{2}(q-q^{-1}) \left( \sigma_1^3 - \sigma_N^3 \right)$$

で与えられる。$q \to 1$ とすれば今まで考察してきた Hamiltonian に戻る。この場合にも代数的 Bethe 仮説法が適用でき、特に Bethe 仮説方程式はこの場合も本質的には多項式方程式系になる。

参考までに、次の作用素

$$e_k = -\frac{1}{2} \left( \sigma_k^1 \sigma_{k+1}^1 + \sigma_k^2 \sigma_{k+1}^2 \right) - \frac{1}{4}(q+q^{-1}) \left( \sigma_k^3 \sigma_{k+1}^3 - 1 \right) + \frac{1}{4}(q-q^{-1}) \left( \sigma_k^3 - \sigma_{k+1}^3 \right)$$

を定義すると、上記 Hamiltonian は

$$H_N = -2 \sum_{k=1}^{N-1} e_k$$

と書き直すことができる。ここで作用素 $e_k$ は Temperley–Lieb 代数を生成する ($|j-k| > 1$)。

$$e_k^2 = (q+q^{-1}) e_k, \quad e_k e_{k\pm 1} e_k = e_k, \quad e_k e_j = e_j e_k.$$

最近この系の Bethe 仮説方程式についても数値計算が実行されており [GHNS]、表現論における結果との比較がなされている。またこの場合も Bethe 仮説法が適用できることが主張されている [GN]。特に変形パラメータ $q$ が 1 の冪根となる最も興味深い場合においても Bethe 仮説法により正しく表現が構成されることが主張されている。これらの内容についてご興味をお持ちの場合は、詳細についてはそれぞれの論文をご参照ください。

# 8　結晶基底との関連

## 8.1　はじめに

Bethe 仮説法からくみ取られる数学には色々な種類がありうるであろうが、以下ではその典型的な一例としてアフィン量子群や関連する組み合わせ論への応用について解説する。なお以降の内容は前節までの内容とは独立して理解できるものである。

この様な方向の研究は 1986 年の Kerov, Kirillov, Reshetikhin [KKR, KR] による研究を嚆矢とする。基本的な思想は以下のようにまとめられる。

| Bethe 仮説法 | 組み合わせ論的類似 |
|---|---|
| Bethe ベクトル | 結晶基底のテンソル積 |
| Bethe 仮説方程式の解<br>(艤装配位との非自明な対応) | 艤装配位 |
| $B_N$ 作用素の積 | 艤装配位写像 $\Phi$<br>(結晶基底のテンソル積と艤装配位の全単射) |

今の所両者の間の関係は分かっていないのだが、以下で述べるように両者の間の対応関係は非常に深くまた密接なものである。筆者の印象では、組み合わせ論的類似物と呼んでいるも

のは Bethe 仮説に対する単なる簡便なモデルという段階をはるかに超えてそれ自身として実質的な数学を含んでいると思われる。もしそうであれば、実は Bethe 仮説法とその組み合わせ論的類似物の両者を含むようなより大きな数学的枠組みが存在し、我々はその二つの相異なる切り口を覗いていたのだった、という可能性もある。両面からの研究が今後長い時間をかけて深化していく過程で色々と予想外の景色を見せていってくれることを期待したい。

組み合わせ論サイドの研究はその後 Kashiwara の結晶基底の理論 [K91] を巻き込んで深く研究が進められ、また現在も活発な研究の対象となっている。個人的な研究の動機としては、Young 盤に関わる Lascoux–Schützenberger 理論などに見られるようないわゆる代数的組み合わせ論の分野に新しい研究の軸を与えたいというものであった。筆者が艤装配位の研究を開始したころは、この課題の研究に対する風当たりが非常に強くずいぶん苦労させられた。経験により艤装配位の重要性を確信していたため、私が適任ではないにせよせめて他の方にバトンを渡せるまでは旗を掲げていたいという一心で孤軍奮闘してきたような状態だったと認識している。その後 10 年を経て少しは事態は改善しただろうか。

以下 §8.2 ではもっとも単純な場合に基本的なアイデアを解説し、その後の内容の概説とした。引き続く節ではもう少し技術的内容にも立ち入るが、各節は概ね独立したトピックスを扱うので、興味のあるところを拾い読みして頂くことができると思う。艤装配位の理論は今の所何の問題もなく様々な種類の代数や表現のクラスに拡張され続けており、また今までの組み合わせ論の方法では到達できなかった領域にも奥深く立ち入るようになってきているので、代数的組み合わせ論の一つの大きな枠組みとなる可能性は十分あると考えている。何より艤装配位自身に数学的意味があることも分かっているので、今後の進展を期待しているところである。

## 8.2 概要：$A_1^{(1)}$ 型ベクトル表現の場合

以下アフィン Lie 代数の種類や関連する Cartan データ等は Kac [K90] の記号法に従って表記することにする。理論の大まかな雰囲気を知るために前節までで扱った $sl_2$ の 2 次元表現のテンソル積と同じ設定で様子を見てみよう。この場合は様々な異なる方法で解析することができるため文献もいろいろと存在するが[13]、その点についての詳細は筆者による別の解説 [Sa12a, Sa12b] 等を参照していただく事にして、ここでは以降の節で述べる結果の概観を与えることに専念することにする。

**艤装配位写像**　この設定において、登場人物の片方である艤装配位は §4 で導入したものと一致する（12 ページ）。そこで関連する記号や用語はそちらで定義したものをそのまま流用することにしよう。一方、Bethe ベクトルの類似は結晶基底のテンソル積となる。$B^{r,s}$ を縦 $r$、横 $s$ の長方形型 Young 準標準盤[14]で各ます目には $A_n^{(1)}$ の場合文字 $1, 2, \ldots, n+1$ のいずれかを書き込んだ物の全体としよう。現在の場合、2 次元表現は二つのベクトル $\boxed{1}$, $\boxed{2} \in B^{1,1}$ で表され、それらの間に Chevalley 作用素の類似物（柏原作用素と呼ばれる）が以下のように作用する：

$$\tilde{f} : \boxed{1} \longmapsto \boxed{2}, \qquad \tilde{f} : \boxed{2} \longmapsto 0, \qquad \tilde{e} : \boxed{2} \longmapsto \boxed{1}, \qquad \tilde{e} : \boxed{1} \longmapsto 0.$$

[13] 一つだけ有名な例を挙げておくと、Tokihiro グループによって箱玉系の初期値問題の解析 [MIT] に用いられた “10-elimination” という手法は以下に述べる艤装配位による方法と同値であることが知られている [KS09]。

[14] semi-standard tableau: Young 図のます目に縦方向には真に増加、横方向には非減少に増加するよう数字を書き込んだもの。

これらの元のテンソル積を考えることもできて、その上での柏原作用素は簡明な作用を持つ（詳細は [KN] を参照）。こうして得られるテンソル積（物としては数字 1 と 2 の列）を path と呼ぶことにしよう。その時艤装配位写像[15] $\Phi$

$$\Phi : \text{path} \longmapsto \text{rigged configurations}$$

が次のようにして定義される。

path が最高ウェイトベクトルである時、結果としてできる艤装配位は §4 で述べた条件、つまり $P_{\nu_i}(\mu,\nu) \geq 0$ かつ $0 \leq J_i \leq P_{\nu_i}(\mu,\nu)$ が成り立っているが、一般の path の場合は $P_{\nu_i}(\mu,\nu) \leq 0$ または $J_i \leq 0$ となることも可能である。例えば、与えられた path が

$$b = 1 \otimes 2 \otimes 2 \otimes 1 \otimes 2 \otimes 2 \otimes 1 \otimes 1 \otimes 1 \otimes 2 \otimes 2 \otimes 1 \otimes 1 \otimes 2 \in (B^{1,1})^{\otimes 14}$$

である場合、対応する艤装配位は

$\Phi(b) =$ 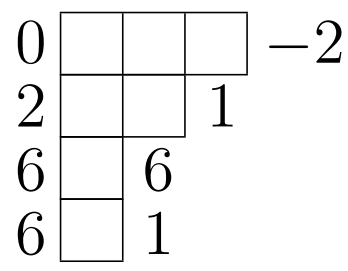

となる。しかしいつでも $J_i \leq P_{\nu_i}(\mu,\nu)$ は成り立っており、rigging が許される最大値を持つ場合、つまり $J_i = P_{\nu_i}(\mu,\nu)$ である時、ストリング $(\nu_i, J_i)$ は**特異**であるという。

さて現在の設定での艤装配位写像 $\Phi$ は次のようなアルゴリズムで与えられる。出発点となる艤装配位は空集合であるとする。以下艤装配位には必ず $(0,0)$ という特異なストリングが含まれていると考えると記述が簡明になるので、そのように仮定しよう（一般性は失われない）。path $b$ の文字 2 の位置を左端から数えて $k_1, k_2, \ldots$ 番目とするとき、順に以下の手続きを繰り返して再帰的に艤装配位を成長させていく。

1. 位置 $k_{j-1}$ まで行った結果、途中段階の艤装配位 $(\eta, I)$ が得られたとする。

2. 次の位置 $k_j$ に対して以下の操作を行う。まず $(\eta, I)$ を長さ $k_j - 1$ の状態に対する艤装配位と考え、vacancy number $P_{\eta_i}(1^{k_j-1}, \eta)$ を全ての行に対し計算し、$\eta$ の特異なストリングを全て特定する。

3. 最も長い特異なストリングから任意に一つ選び、その行の長さを 1 増やす[16]。こうして得られた新しいヤング図を $\eta'$ を書く。

4. 新しい rigging $I'$ を以下のように定める。$\eta \to \eta'$ において変化しなかった行に対応する成分は変化させない。一方 $\eta \to \eta'$ において変化した行を $\eta'_i$ とする時、対応する rigging $I'_i$ は $(\eta'_i, I'_i)$ が特異なストリングとなるように $I'_i = P_{\eta'_i}(1^{k_j}, \eta')$ と定める。

5. 以上の手続きを全ての $k_j$ について行い最終的に得られる結果を $(\nu, J)$ とするとき $\Phi(b) = (\nu, J)$ と定める。

[15] この写像を Kerov–Kirillov–Reshetikhin 写像等々と研究者の名前をつけて呼ぶこともある。しかし理論の主要な創設者である Kirillov 氏自身はこの写像に人名をつけて呼ぶことに強く反対しておられ、様々な場所でそのように述べておられる。まことに尤もなお考えだと思うので、2008 年に初めて共同研究をさせて頂いた機会に相談の上、艤装配位全単射（rigged configuration bijection）という中立的な名称を提唱している。大方のご賛同を頂ければ幸いである。

[16] 考えている艤装配位に $(0,0)$ 以外の特異なストリングが存在しない場合には、ストリング $(0,0)$ を変化させて長さ 1 のストリングが生成されると考える

計算例は次項であたえる。なお以上の手続きを逆にすると逆写像 $\Phi^{-1}$ が得られる。
　艤装配位写像の Mathematica によるプログラム例が

https://sites.google.com/site/affinecrystal/rigged-configurations

にあるのでご利用いただきたい（コード自体はより一般的に $D_n^{(1)}$ 型の場合も含んでいる）。

**箱玉系**　量子群に関する最も重要な対象の一つは $R$ 行列と呼ばれ、表現のテンソル積の左右を入れ替える非自明な同型写像である。その結晶基底における類似物である同型写像 $R: B \otimes B' \mapsto B' \otimes B$ を組み合わせ $R$ 行列と呼ぶ。Drinfeld の普遍 $R$ 行列の理論なども存在するが、結晶基底の場合でも組み合わせ $R$ 行列を具体的に求めるのは非常に難しい。§8.3 や §8.7 で述べる様にこのような問題に艤装配位は強力な手段を提供する。

　ただし現在の設定では以下のように簡単な表示を持つ。$B^{1,k}$ の元は長さ $k$ の文字 1 と 2 による準標準盤であり、$(a,b)$ で文字 1 が $a$ 個、文字 2 が $b$ 個含まれる元を表すとすれば、組み合わせ $R$ 行列

$$R: (a,b) \otimes (c,d) \longmapsto (c',d') \otimes (a',b')$$

は具体的に

$$\begin{aligned}
a' &= a + \min(b,c) - \min(a,d)\\
b' &= b - \min(b,c) + \min(a,d)\\
c' &= c - \min(b,c) + \min(a,d)\\
d' &= d + \min(b,c) - \min(a,d)
\end{aligned}$$

と表示される。

　箱玉系は $u_l = (l,0)$ なる特別な元と、与えられた path $b = b_1 \otimes \cdots \otimes b_L$ に対して時間発展 $T_l$ を次のように定義することにより与えられる力学系である。まず左側から繰り返し $R$ 行列を作用させて

$$\begin{aligned}
\underline{u_l \otimes b_1} \otimes b_2 \otimes \cdots \otimes b_L &\xrightarrow{R} b_1' \otimes \underline{u_l^{(1)} \otimes b_2} \otimes \cdots \otimes b_L \xrightarrow{R} b_1' \otimes b_2' \otimes u_l^{(2)} \otimes \cdots \otimes b_L\\
&\xrightarrow{R} \cdots \xrightarrow{R} b_1' \otimes b_2' \otimes \cdots \otimes b_L' \otimes u_l^{(L)}
\end{aligned} \tag{14}$$

としたときに

$$T_l(b) := b_1' \otimes b_2' \otimes \cdots \otimes b_L'$$

と時間発展を定める。すなわち

$$R: u_l \otimes b \longmapsto T_l(b) \otimes u_l^{(L)} \tag{15}$$

となる。結晶基底の場合の Yang–Baxter 関係式により $T_l$ どうしは可換な時間発展の族を成すので、箱玉系は量子可積分系である。箱玉系の研究は数学上の進展を結果としてもたらしたので[17]、その意味で意義のある研究であると考えている。

[17] 艤装配位写像の研究は当初ある組み合わせ論的恒等式の証明を主な目的として研究が進められていたが、その後他の方法でも当該恒等式を証明することが可能であることが分かった。しかし箱玉系の研究という観点からは艤装配位写像のようなより深い理論を避けて通ることはできず、結果として艤装配位の研究を継続することとなった。実は研究の動機となったばかりではなく、艤装配位の数学的構造には箱玉系が本質的に関わっていることが分かっている（33 ページ参照）。

非常に簡単な例を見てみよう。path $b \in (B^{1,1})^{\otimes L}$ は文字 1 と 2 の列とみなせるのでテンソル積の記号を略記すると以下のような時間発展を考えることができる。

$$
\begin{aligned}
t=0: &\quad 22211112111111\\
t=1: &\quad 11122211211111\\
t=2: &\quad 11111122122111\\
t=3: &\quad 11111111211222
\end{aligned}
$$

ここで時刻 $t$ における状態とは $(T_3)^t(b)$ を表す。文字 2 の塊はソリトン系における孤立波とみなすことができ、箱玉系はソリトン系となっている。具体的には、文字 2 の長さ $k$ の塊は他の孤立派の影響のない状況では速度 $k$ で自由伝搬し、他の孤立派と衝突すると非自明な相互作用ののち、元の孤立派が再現する（ただし位置には位相差と呼ばれる変化が生じる）。さて各時刻に対する艤装配位を計算してみよう。前項におけるアルゴリズムの定義を参考に計算すると、$t=0$ の場合文字 2 の位置 $k_j$ は $1,2,3,8$ であり、写像 $\Phi$ の計算は以下の通り。

$$
\emptyset \xrightarrow{1} \square\ {-1} \xrightarrow{2} \square\square\ {-2} \xrightarrow{3} \square\square\square\ {-3} \xrightarrow{8} \begin{array}{l} \square\square\square\ {-3} \\ \square\ 4 \end{array}
$$

一方 $t=2$ の場合は以下のようになる。

$$
\emptyset \xrightarrow{7} \square\ 5 \xrightarrow{8} \square\square\ 4 \xrightarrow{10} \begin{array}{l} \square\square\ 4 \\ \square\ 6 \end{array} \xrightarrow{11} \begin{array}{l} \square\square\square\ 3 \\ \square\ 6 \end{array}
$$

各 $t$ に対する計算をまとめると以下のような結果となる。

$$
\Phi(T^t_\infty(b)) = \begin{array}{l} \square\square\square\ {-3+3t} \\ \square\ 4+t \end{array}
$$

Young 図の形は時間によって変化せず、各 $\nu_i$ は孤立波の長さに対応していることが分かる。また riggings は対応する $\nu_i$ の長さに比例して $\nu_i t$ と変化することが分かる[18]。

もう一つ大きめの例を見てみよう。先ほども見た $(B^1)^{\otimes 27}$ での箱玉系の例:

$$
\begin{aligned}
t=0: &\quad 122122111221121111111111111\\
t=1: &\quad 111211222112212111111111111\\
t=2: &\quad 111121111221121222111111111\\
t=3: &\quad 111112111112212111222111111\\
t=4: &\quad 111111211111121221111222111\\
t=5: &\quad 111111121111112112211111222
\end{aligned}
$$

各時刻に対応する艤装配位:

$$
\Phi(T^t_\infty(b)) = \begin{array}{l} \square\square\square\ {-2+3t} \\ \square\square\ 1+2t \\ \square\ 6+t \\ \square\ 1+t \end{array}
$$

この場合も長さ $\nu_i$ の孤立波を表すストリングの rigging が $\nu_i t$ のように変化することが見て取れる。実は Bethe 仮説方程式の解に現れるストリングも、対応する Bethe ベクトルで下向きスピンの塊（マグノンと呼ばれる）に対応すると考えられているので、これも Bethe 仮説とその組み合わせ論的類似の間に存在する非常に深い関係の一端を表しているといえよう。

§8.3 で述べるように、一般に箱玉系の時間発展に対して艤装配位の Young 図は保存量（作用変数）を与え、rigging は線形に時間発展する（角変数）。すなわち艤装配位自身が明瞭な数学的意味を持つことを保証している。同じく §8.3 で述べるように、この結果は艤装配位を結晶基底の言葉で理解するための鍵を与える。その様な結果は、例えば §8.6 で述べるような艤装配位写像の Loop Schur 関数による表示の導出などに応用されている。

[18] 一般の時間発展 $T_l$ に対しては rigging は $\min(\nu_i, l)t$ と変化する。よって $\max_i \nu_i \leq l$ なら $T_l = T_\infty$ となる。

**周期系**　箱玉系には周期境界条件を課したバージョンも存在する。周期系の結果は、今の所 $A_1^{(1)}$ 型ベクトル表現の場合でのみ完成しているのでここで触れておくことにする。本項の内容は以降は登場しないので、ここでは文献についてもある程度触れておくことにする。

箱玉系の周期版はそれ以前から知られていたが [YT]、ここでは [KTT] による結晶基底を用いた定式化を紹介する。path $b \in (B^{1,1})^{\otimes L}$ で、文字 2 の個数が $L/2$ 以下のものを考える。$b$ に周期境界条件を課すと、元来は円周状のものが得られるが、適当に切り口を入れて直線状の通常の path として取り扱う。切り口としては、便宜上 $b$ が最高ウェイト条件を満たすものを考えよう。このとき、定義 (15) にならって書くと次のような性質がある。

$$R : u_l \otimes b \longmapsto T_l(b) \otimes v, \qquad R : v \otimes b \longmapsto \overline{T}_l(b) \otimes v.$$

すなわち左側の同型で得られた $v$ を用いると、右側の様に周期境界条件と適合する時間発展が定義できる。そこで $\overline{T}_l$ を周期箱玉系の定義とするのである。この場合も量子可積分系となっていることが分かる。周期箱玉系のデモンストレーションが

http://demonstrations.wolfram.com/PeriodicBoxBallSystem

にあるのでご参照頂きたい。

直線状の箱玉系の時と同様にこの場合も艤装配位上で Young 図が作用変数、rigging が角変数となることは同様である（rigging は適当な同値関係で割る必要がある）。興味深い問題として、艤装配位上で線形に時間発展させた後に結晶基底のテンソル積（文字 1 と 2 の並び）を与える解析的表示式を求める問題を考えよう [KS06]。通常の多変数テータ関数の**超離散極限**[19]**（トロピカル極限）**からスタートしよう。

$$\begin{aligned}\Theta(\mathbf{z}) &= \lim_{\epsilon \to +0} \epsilon \log \left( \sum_{\mathbf{n} \in \mathbb{Z}^g} \exp\Bigl( -\frac{{}^t\mathbf{n}A\mathbf{n}/2 + {}^t\mathbf{n}\mathbf{z}}{\epsilon} \Bigr) \right) \\ &= -\min_{\mathbf{n} \in \mathbb{Z}^g} \{ {}^t\mathbf{n}A\mathbf{n}/2 + {}^t\mathbf{n}\mathbf{z} \}.\end{aligned} \tag{16}$$

ここで簡単のため艤装配位の Young 図が

$$\nu = \{ i_1 < i_2 \cdots < i_g \} \tag{17}$$

という形をしていたと仮定しよう。そのとき各 $i_k$ を添え字とする行列 $A$ を

$$A_{i,j} = \delta_{i,j} p_i + 2 \min(i,j), \tag{18}$$

と定める（ここで $p_i = P_i(1^L, \nu)$ と書いた）。(17) に対する riggings が $J_1, \ldots, J_g$ である時ベクトル $\mathbf{J} = (J_1, \ldots, J_g)$ と定める。またベクトル $\mathbf{h}_l = (\min(i,l))_{i \in \nu}$ と $\mathbf{p} = (p_i)_{i \in \nu}$ を定めるとき、path $b$ の $k$ 番目の要素 $b_k$ に含まれる文字 2 の個数 $x(k)$ は

$$\begin{aligned} x(k) = \Theta\bigl(\mathbf{J} - \frac{\mathbf{p}}{2} - k\mathbf{h}_1\bigr) - \Theta\bigl(\mathbf{J} - \frac{\mathbf{p}}{2} - (k-1)\mathbf{h}_1\bigr) \\ - \Theta\bigl(\mathbf{J} - \frac{\mathbf{p}}{2} - k\mathbf{h}_1 + \mathbf{h}_\infty\bigr) + \Theta\bigl(\mathbf{J} - \frac{\mathbf{p}}{2} - (k-1)\mathbf{h}_1 + \mathbf{h}_\infty\bigr) \end{aligned} \tag{19}$$

[19] 正整数係数の有理関数の場合なら形式的に $\times \mapsto +$、$/ \mapsto -$、$+ \mapsto \min$ と置き換えればよい。超離散化された有理関数は入力も整数、出力も整数のデジタルな世界となる。43 ページも参照。

と表される[20]。この表示では、周期境界条件に由来する rigging の同値関係による商は自動的に満たされる。よって単純に $\mathbf{J} \mapsto \mathbf{J} + \mathbf{h}_l$ とするだけで周期箱玉系の時間発展 $T_l$ が実現される。この結果は Dubrovin, Matveev, Novikov[21]による有名な結果 [DMN] （KdV 方程式の周期解のタウ関数が多変数テータ関数となり解はその対数の 2 階微分となること、1976 年）の完全な離散化が達成されたことを意味している。

**超離散極限について**　連続的な関数の世界になじみが深い方々にとっては超離散極限（またはトロピカル極限）を取った後の世界は連続的な世界のおもちゃのように感じられてしまう事がしばしばあるようである。確かに例えば多項式であれば超離散化すると係数の情報がすべて吹き飛んでしまうし、一見実に情報不足な世界に見えてしまうのもうなずける。その様な方々のために、超離散化されることによって初めて生じる新しい性質も確かに存在するということを紹介しよう。ここでは例として超離散化によって新たに粒子性を獲得するような状況を見てみよう。以下の例が実際に数学的に重要であるのか、という問題はさておき、少なくとも離散化された世界が連続的な世界に従属する存在ではなく、それ自身独立した固有の世界であることをお示ししたいと思う。

ここで紹介するのは **Box-Basket-Ball 系** [LPS1]（略して BBB 系）と呼ばれるソリトン系である（解説 [Sa12a] も参照)。その後 Yura [Y13] によりソリトン理論の観点から研究が進められ、他の既知の系との相違点も検討された。箱玉系における状態 $(a, b)$ を拡張して状態 $(a, b, c) \in (\mathbb{Z}_{\geq 0})^3$ を考える。組み合わせ $R$ 行列の代わりに **whurl 関係式**

$$R : (a, b, c) \otimes (d, e, f) \mapsto (d', e', f') \otimes (a', b', c')$$

を考え、$u_l = (l, 0, 0)$ から出発して箱玉系の場合 (15) と同様に時間発展 $T_l$ を定義すると BBB 系が定義され、量子可積分系となる。ここで whurl 関係式とは、具体的には

$$\begin{aligned}
a' &= a - \min(a + b, a + c, b + f) + \min(e + c, d + c, d + b)\\
b' &= b - \min(a + b, a + c, b + f) + \min(a + e, d + f, e + f)\\
c' &= c - \min(e + c, d + c, d + b) + \min(a + e, d + f, e + f)\\
d' &= d + \min(a + b, a + c, b + f) - \min(e + c, d + c, d + b)\\
e' &= e + \min(a + b, a + c, b + f) - \min(a + e, d + f, e + f)\\
f' &= f + \min(e + c, d + c, d + b) - \min(a + e, d + f, e + f)
\end{aligned}$$

で与えられる変換である。この様な変換を与えるような有理関数を構成することは容易であり、かつ多数構成することができるが、BBB 系ではそのような連続的な世界には見られない以下のような粒子的性質が存在する。

---

[20] 証明は周期境界条件を課さない一般の $b$ に対して、艤装配位 $(\nu, J) = (\nu_i, J_i)_{i=1}^g$（ここで $\nu_i$, $i = 1, 2, \ldots$ は重複も許す）から

$$\begin{aligned}
x(k) &= \tau_0(k) - \tau_0(k-1) - \tau_1(k) + \tau_1(k-1),\\
\tau_r(k) &= -\min_{n \in \{0,1\}^g} \left\{ \sum_{i=1}^{g} (J_i + r\nu_i - k) n_i + \sum_{i,j=1}^{g} \min(\nu_i, \nu_j) n_i n_j \right\} \quad (r = 0, 1),
\end{aligned}$$

という式で求まることを用いる。ここで $n = (n_1, n_2, \cdots, n_g)$。詳細は §8.5 を参照。

[21] Sergei Petrovich Novikov (b. 1938). 初期のトポロジーにおける重要な研究ののち数理物理学において多大な功績を残している。1970 年フィールズ賞受賞。

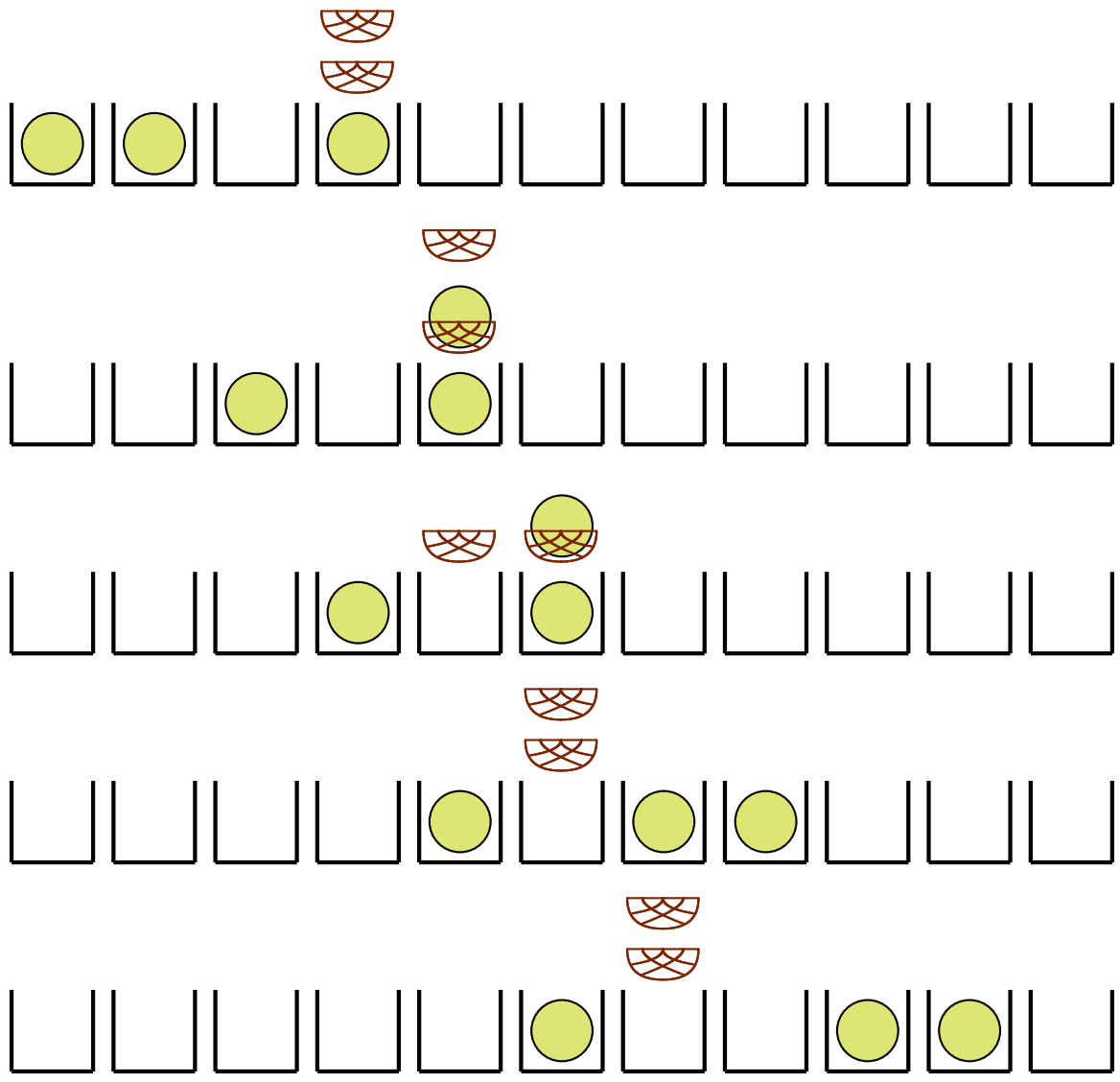

図 2: BBB 系の時間発展の例.

状態 $(a,b,c)$ の各成分について、$a$ を**空きスペース**の数、$b$ は**バスケット**の数、$c$ は**ボール**の数と解釈する。最初に容量 1 の箱を直線状に並べ、箱の上にはバスケットをいくらでも重ねることができ、一方箱とバスケットにはそれぞれ 1 つずつボールを入れられると考える（そうすると条件 $a+c=b+1$ を課すことになる）。そのとき whurl 関係式から定まる時間発展 $T_\infty$ は次のようにして粒子的に記述することができる。

> **時間発展の手順：**バスケットのうち空のものを全て一つ右の箱の上に移し、それ以外のものはそのままにしておく。次に玉を左側から順に見ていって、それぞれの玉より右側にある空箱または空きバスケットのうち最も左側のものに順に移していく。ただし既に動かした玉にはそれ以上触れない。

この手順はバスケットが一つもない場合は前述の箱玉系の時間発展 $T_\infty$ と一致する事を注意しておく（箱玉系では文字 1 を空き箱、文字 2 をボールが一つ入った箱と解釈する）。簡単な例を図 2 に挙げておこう。物理の用語でまとめておくと、ボールは一つの場所に一つしか入ることのできない Fermi 粒子、バスケットは一つの場所にいくらでも入る Bose 粒子（この二種類が自然界に存在する量子力学的粒子の大分類を与える）であり両者の間にはボールがバスケットに入るという形で非自明な相互作用を持つ。この様にして定義される BBB 系はソリトン系となっている（どのような孤立波が存在するかは解説 [Sa12a] 参照）。以上の様な性質は whurl 関係式から見て全く非自明な結果であることを指摘しておく。

まとめると、粒子とは本質的に離散的で無くてはならないので、超離散化された世界（デジタルな世界）には連続な世界には存在しない性質が存在することについてはご理解いただけたのではないかと思う。逆に超離散的な世界から出発して連続的な類似物を構成しようとすると、同一の関係式から無数の対応する連続関数を構成できてしまう。その様な中からどのような基準でもって数学的に正当性を持った関数を選び出すのかは今の所はっきりとした指針がある様には感じられないし、そもそも連続な世界での類似物が超離散的な世界で持っていたのと匹敵するような顕著な興味深い性質を持っているのかも非自明な問題であろう。筆者としては結晶基底自体も超離散化された世界で最も自然に認識できる対象なのではないかと思っているのだが、連続な世界で類似物を構成するために多大な労力を払っている方々

が多くいらっしゃるので、研究の成果を見守っているところである。

## 8.3 $A_n^{(1)}$ 型艤装配位

$A_n^{(1)}$ 型の結晶基底 (Kirillov–Reshetikhin crystal) のテンソル積に対応する艤装配位について述べる。その様な写像は Kerov–Kirillov–Reshetikhin（[KKR, KR]、1986 年）によって先駆的な研究が行われ、最終的に Kirillov–Schilling–Shimozono（[KSS]、1999 年）によって $A_n^{(1)}$ 型一般の場合に確立された（[K01] も参照）。

**クリスタル $B^{r,s}$** $A_n^{(1)}$ 型の場合 KR クリスタル $B^{r,s}$ とは、集合としては文字 $1, 2, \ldots, n+1$ による縦 $r$ 横 $s$ の長方形型 Young 準標準準盤であった。テンソル積 $B^{r,s} \otimes B^{r',s'}$ が自然に定義され、代数構造も簡明に入る。具体的な形は $a \neq 0$ に対する $\widetilde{e}_a, \widetilde{f}_a$ については [KN] を、$\widetilde{e}_0, \widetilde{f}_0$ については [Sh98] をそれぞれ参照。

テンソル積 $\bigotimes_{k=1}^{L} B^{r_k,s_k} = B^{r_1,s_1} \otimes \cdots \otimes B^{r_L,s_L}$ の元を path とよぶ。path の形状から Young 図 $\mu^{(a)}$（$a = 0, \ldots, n-1$）を以下のように定める。すなわち $\bigotimes_{k=1}^{L} B^{r_k,s_k}$ の中に $B^{a+1,l}$ が $m$ 個含まれる場合、Young 図 $\mu^{(a)}$ に長さ $l$ の行が $m$ 本存在する。

艤装配位写像 $\Phi$ は全単射

$$\Phi : \text{path} \longmapsto \text{rigged configuration}$$

を与えるが、そのアルゴリズムの定義は次節 §8.4 に回すことにする。

**最高ウェイトベクトルに対する艤装配位** 前項で定義した Young 図 $\mu^{(a)}$（$a = 0, \ldots, n-1$）と Young 図 $\nu^{(a)}$（$a = 1, \ldots, n$）の各行 $\nu_i^{(a)}$ に整数 $J_i^{(a)}$ を付け加えた以下のような対象を考える（以下 $\mu$ の部分は適宜略記する）。

$$(\nu, J) = \Big(\mu^{(0)}, \ldots, \mu^{(n-1)}, (\nu^{(1)}, J^{(1)}), \ldots, (\nu^{(n)}, J^{(n)})\Big). \tag{20}$$

Young 図 $\eta$ に対し $Q_l(\eta)$ を Young 図 $\eta$ の左側 $l$ 列分に含まれるます目の総数としよう。$(\nu, J)$ に対し **vacancy number** を以下の通り定義する：

$$P_l^{(a)}(\mu, \nu) = Q_l(\mu^{(a-1)}) + Q_l(\nu^{(a-1)}) - 2Q_l(\nu^{(a)}) + Q_l(\nu^{(a+1)}). \tag{21}$$

Bethe 仮説での性質に基づき、組 $(\nu_i^{(a)}, J_i^{(a)})$ をストリングと呼び、$\nu_i^{(a)}$ をその長さと呼ぼう。

その時 $(\nu, J)$ が最高ウェイトベクトルに対する艤装配位であるとは以下の条件を満たす場合である。

- 全ての $\nu_i^{(a)}$ に対し $P_{\nu_i^{(a)}}^{(a)}(\mu, \nu) \geq 0$。
- 全ての $J_i^{(a)}$ が $0 \leq J_i^{(a)} \leq P_{\nu_i^{(a)}}^{(a)}(\mu, \nu)$ を満たす。

$\nu^{(a)}$ の事を **configuration**、$J_i^{(a)}$ の事を **rigging** と呼ぶ。rigging が許される最大値を持つとき、そなわち $J_i^{(a)} = P_{\nu_i^{(a)}}^{(a)}(\mu, \nu)$ であるとき、**特異なストリング**と呼ぶ。

**一般の騰装配位**　本項目の内容はやや技術的であり、飛ばしてそのまま先に進んでいただいても問題はありません。組み合わせ $R$ 行列に関する話題は次項で、また箱玉系に関する基本的結果は 32 ページから始まります。

他の一般の騰装配位は Schilling [Sc06] にならって騰装配位上に柏原作用素を定義することによって得られる。用語としてストリング $(\nu_i^{(a)}, J_i^{(a)})$ の **corigging** とは

$$P^{(a)}_{\nu_i^{(a)}}(\mu, \nu) - J_i^{(a)}$$

の事を表す。以下 $(\nu^{(a)}, J^{(a)})$ に $(0,0)$ という特異なストリングを付け加えておくと記述が簡明になる。$x_\ell$ を $(\nu^{(a)}, J^{(a)})$ の最小の rigging とする。今の規約より $x_\ell \leq 0$ となる。

(1) もし $x_\ell = 0$ ならば $\widetilde{e}_a(\nu, J) = 0$ と定める。一方 $x_\ell < 0$ ならば $(\nu^{(a)}, J^{(a)})$ のストリングで rigging が $x_\ell$ であるようなもののうち長さが最小のものの長さを $\ell$ とする。その時 $\widetilde{e}_a(\nu, J)$ はストリング $(\ell, x_\ell)$ を $(\ell-1, x_\ell+1)$ に変化させ、他の rigging を元の corigging を保つように変化させる。

(2) $(\nu^{(a)}, J^{(a)})$ のストリングで rigging が $x_\ell$ であるようなもののうち長さが最大のものの長さを $\ell$ とする。そのとき $\widetilde{f}_a(\nu, J)$ はストリング $(\ell, x_\ell)$ を $(\ell+1, x_\ell-1)$ に置き換え、他の rigging を元の corigging を保つように変化させる。その結果 rigging が対応する vacancy number より真に大きなものが存在する場合 $\widetilde{f}_a(\nu, J) = 0$ と定める。

$\widetilde{f}_a$ の定義において、$(\nu^{(a)}, J^{(a)})$ の rigging が全て 0 より大きかった場合には、$(0,0)$ というストリングに $\widetilde{f}_a$ が作用してストリング $(1,-1)$ が生成されることになることに注意。例えば $(B^{1,1})^{\otimes 8}$ の元（つまり $\mu^{(0)} = (1^8)$ となる）

$$b = \boxed{1} \otimes \boxed{1} \otimes \boxed{2} \otimes \boxed{1} \otimes \boxed{1} \otimes \boxed{3} \otimes \boxed{2} \otimes \boxed{1}$$

を考え $\Phi(b)$ に $\widetilde{f}_1$ を順に作用させていくと以下のようになる。

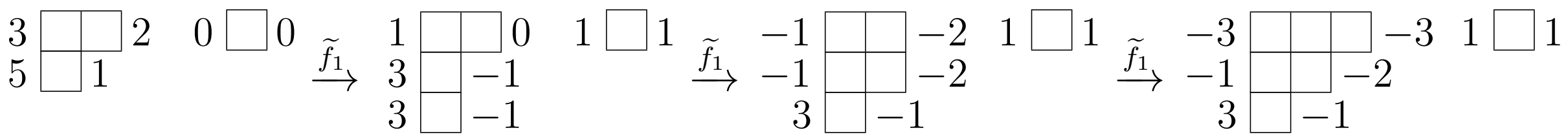

もう一度 $\widetilde{f}_1$ を作用させると 0 となる。ここで騰装配位の表示として、左から $(\nu^{(1)}, J^{(1)})$、$(\nu^{(2)}, J^{(2)})$ を書き、$\nu^{(a)}$ の右側に rigging を、左側に vacancy number を書いた。

この様にして定義される柏原作用素がクリスタルの定義を満たすことは最初 [Sc06] によって証明されたが、Stembridge の結果に依存した証明であったためそのまま拡張することはできない証明法であった。その点は [Sa14] において一般性のある方法に改められ、その後の騰装配位におけるクリスタル構造の導入において標準的な方法となった[22]。

重要な性質として騰装配位写像 $\Phi$ は柏原作用素を保つこと

$$[\widetilde{e}_a, \Phi] = [\widetilde{f}_a, \Phi] = 0 \tag{22}$$

が知られている [DS, Sa14]。写像 $\Phi$ が複雑なため極めて非自明な結果である。後で $D_n^{(1)}$ 型の場合にもふれるが、これらの性質は代数や表現の種類によらず非常に広い範囲で確認されており、騰装配位の著しい自然さを表す例となっている。

[22] 例えば結晶構造の公理に現れる $\varepsilon_a$ という関数は [Sa14] において騰装配位上に導入され、その後 $B(\infty)$ クリスタルの構造を騰装配位に導入するときに本質的であった [SalScr14]。

**組み合わせ R 行列とエネルギー関数**　この項目と引き続く計算例もアルゴリズムの詳細に興味のない方は細部は軽く流して「箱玉系の逆散乱形式」の項（32 ページ）に進んでいただいて構いません。

組み合わせ R 行列とは表現の左右を入れ替える同型写像 $R: B^{r,s} \otimes B^{r',s'} \to B^{r',s'} \otimes B^{r,s}$ の事であった（以下 $\overset{R}{\simeq}$ とも記述する）。クリスタル $B$ のアフィン化を

$$\mathrm{Aff}(B) = \{b[d] \mid b \in B,\, d \in \mathbb{Z}\} \tag{23}$$

と定義する。組み合わせ $R$ 行列により元 $b \otimes b' \in B^{r,s} \otimes B^{r',s'}$ が $\tilde{b}' \otimes \tilde{b} \in B^{r',s'} \otimes B^{r,s}$ に写されるとき、$\mathrm{Aff}(B)$ の上での $R$ 行列を

$$b[d] \otimes b'[d'] \overset{R}{\simeq} \tilde{b}'[d' - H(b \otimes b')] \otimes \tilde{b}[d + H(b \otimes b')] \tag{24}$$

で定める。ここでエネルギー関数 $H: B \otimes B' \to \mathbb{Z}$ とは、(24) で定まる $R$ 行列が $\mathrm{Aff}(B) \otimes \mathrm{Aff}(B') \otimes \mathrm{Aff}(B'')$ 上で Yang–Baxter 関係式

$$(R \otimes 1)(1 \otimes R)(R \otimes 1) = (1 \otimes R)(R \otimes 1)(1 \otimes R) \tag{25}$$

を満たすように定義されたものである [(KMN)$^2$]。

組み合わせ $R$ 行列は一般には極めて非自明な写像であるが、$A_n^{(1)}$ 型の場合は既存の代数的組み合わせ論の成果を使用して以下のように具体的に記述することができる[23]。

Young 盤の各行が上から順に $y_1, y_2, \ldots, y_r$ であったとしよう。その時各行から得られる数列を並べて、$Y$ の row word を $row(Y) = y_r y_{r-1} \ldots y_1$ と定義しよう。また $x = (x_1, x_2, \ldots)$ と $y = (y_1, y_2, \ldots)$ を二つの分割とするとき、$x$ と $y$ の連結を $(x_1 + y_1, x_2 + y_2, \ldots)$ と表す。$Y \leftarrow a$ で Young 盤に文字 $a$ を Schensted–Knuth[24][K70] の row insertion[25]することを表し、更に $Y \leftarrow ab = (Y \leftarrow a) \leftarrow b$ 等と書くしよう。

**定理 [Sh98]**

元 $b \otimes b' \in B^{r,s} \otimes B^{r',s'}$ と $\tilde{b}' \otimes \tilde{b} \in B^{r',s'} \otimes B^{r,s}$ が組み合わせ $R$ 行列によって

$$b \otimes b' \overset{R}{\simeq} \tilde{b}' \otimes \tilde{b} \tag{26}$$

となることと

$$(b' \leftarrow row(b)) = (\tilde{b} \leftarrow row(\tilde{b}')) \tag{27}$$

となることは同値。更に $H(b \otimes b')$ の値は $(b' \leftarrow row(b))$ のます目のうち分割 $(s^r)$ と $(s'^{r'})$ を連結したものの外側に存在するます目の個数で与えられる。

従って $(b' \leftarrow row(b))$ の情報から $\tilde{b}$ と $\tilde{b}'$ を見つけるアルゴリズムを記述すればよいことになる。Young 図 $Y$ とその（左上部分の）部分集合である Young 図 $Y'$ から作られる集合 $\theta = Y \setminus Y'$ のます目の総数が $r$ であり、かつ $\theta$ には同じ行に 2 つ以上のます目が存在しないとき、$\theta$ は vertical $r$-strip と呼ばれる。

[23] アフィン柏原作用素 $\tilde{e}_0$ や $\tilde{f}_0$ では jeu de taquin と呼ばれるアルゴリズムを使用したりと、$A$ 型 KR クリスタルの理論は代数的組み合わせ論の成果が総動員という感じである

[24] Donald Ervin Knuth (b. 1938). LaTeX でも日々お世話になっている数学者。

[25] $Y \leftarrow a$ とは行 $y_1$ の文字で $a$ より真に大きな文字のうち最も左のもの（文字 $x_1$ としよう）と $a$ を交換し、そうして得られる $x_1$ を用いて第二行 $y_2$ で同様の操作を行い、以下同様の操作をどこかで行 $y_k$ 中に $x_{k-1}$ より大きな文字が存在しなくなるまで繰り返す。その時 $y_k$ の右に $x_{k-1}$ を付け加えて終了する。

Young 盤 $Y = (b' \leftarrow row(b))$ の左上部分の長方形 $(s^r)$ 型の部分集合を $Y'$ とする。集合 $\theta = Y \setminus Y'$ の文字に以下のようにして番号 $1, 2, \ldots, r's'$ を振る。$\theta_1$ を $\theta$ の中でなるべく上方かつ右側に取った vertical $r'$-strip とする。その時 $\theta_1$ のます目に対し下から順に番号 $1, 2, \ldots, r'$ を振る。次に $\theta \setminus \theta_1$ に対して同様に vertical $r'$-strip $\theta_2$ をとり、同様に番号 $r'+1, \ldots, 2r'$ をふる。同様にして帰納的に $\theta$ の全てのます目に番号を振っていく。

さて $\theta$ に振られた番号の順に row insertion の逆を行う。番号 1 が振られた文字に対して row insertion の逆を行って Young 盤から放出される文字を $u_1$ とし、$Y_1$ を $(Y_1 \leftarrow u_1) = Y$ となるような Young 盤とする。次に Young 盤 $Y_1$ で番号 2 が振られていたます目にある文字から row insertion の逆を行って文字 $u_2$ と Young 盤 $Y_2$ を得る。同様にして帰納的に文字 $u_{r's'}$ と Young 盤 $Y_{r's'}$ まで得る。最終的に

$$\tilde{b}' = (\emptyset \leftarrow u_{r's'} u_{r's'-1} \cdots u_1), \qquad \tilde{b} = Y_{r's'} \tag{28}$$

が答えとなる。

**例**　以上のアルゴリズムの計算例として次のテンソル積を考える。

$$b \otimes b' = \begin{array}{|c|c|}\hline 1 & 1 \\ \hline 2 & 4 \\ \hline \end{array} \otimes \begin{array}{|c|c|}\hline 3 & 4 \\ \hline 4 & 5 \\ \hline 5 & 6 \\ \hline \end{array} \in B^{2,2} \otimes B^{3,2}.$$

$b$ の row word は $row(b) = 2411$ なので

$$\left( \begin{array}{|c|c|}\hline 3 & 4 \\ \hline 4 & 5 \\ \hline 5 & 6 \\ \hline \end{array} \leftarrow 2411 \right) = \begin{array}{|c|c|c|}\hline 1 & 1 & 4_3 \\ \hline 2 & 4 \\ \cline{1-2} 3_6 & 5_2 \\ \cline{1-2} 4_5 & 6_1 \\ \cline{1-2} 5_4 \\ \cline{1-1} \end{array}$$

を得る。ここで数字につけた添え字は row insertion の逆を行う順番を示す。最初に文字 $6_1$ から出発すると

$$\left( \begin{array}{|c|c|c|}\hline 1 & 4 & 4 \\ \hline 2 & 5 \\ \cline{1-2} 3 & 6 \\ \cline{1-2} 4 \\ \cline{1-1} 5 \\ \cline{1-1} \end{array} \leftarrow 1 \right) = \begin{array}{|c|c|c|}\hline 1 & 1 & 4 \\ \hline 2 & 4 \\ \cline{1-2} 3 & 5 \\ \cline{1-2} 4 & 6 \\ \cline{1-2} 5 \\ \cline{1-1} \end{array}, \qquad \text{従って}, \qquad Y_1 = \begin{array}{|c|c|c|}\hline 1 & 4 & 4_3 \\ \hline 2 & 5 \\ \cline{1-2} 3_6 & 6_2 \\ \cline{1-2} 4_5 \\ \cline{1-1} 5_4 \\ \cline{1-1} \end{array}, \qquad u_1 = 1.$$

次に $Y_1$ のます目 $6_2$ から同様に計算する。最後まで繰り返すと $u_6 u_5 \cdots u_1 = 321541$ および $Y_6 = \begin{array}{|c|c|}\hline 4 & 4 \\ \hline 5 & 6 \\ \hline \end{array}$ を得る。$(\emptyset \leftarrow 321541) = \begin{array}{|c|c|}\hline 1 & 1 \\ \hline 2 & 4 \\ \hline 3 & 5 \\ \hline \end{array}$ なので、結局

$$\begin{array}{|c|c|}\hline 1 & 1 \\ \hline 2 & 4 \\ \hline \end{array} \otimes \begin{array}{|c|c|}\hline 3 & 4 \\ \hline 4 & 5 \\ \hline 5 & 6 \\ \hline \end{array} \simeq \begin{array}{|c|c|}\hline 1 & 1 \\ \hline 2 & 4 \\ \hline 3 & 5 \\ \hline \end{array} \otimes \begin{array}{|c|c|}\hline 4 & 4 \\ \hline 5 & 6 \\ \hline \end{array}, \qquad H\left( \begin{array}{|c|c|}\hline 1 & 1 \\ \hline 2 & 4 \\ \hline \end{array} \otimes \begin{array}{|c|c|}\hline 3 & 4 \\ \hline 4 & 5 \\ \hline 5 & 6 \\ \hline \end{array} \right) = 3.$$

ここでエネルギー関数は $b$ と $b'$ を連結した $\begin{array}{|c|c|c|c|}\hline \phantom{0} & \phantom{0} & \phantom{0} & \phantom{0} \\ \hline \phantom{0} & \phantom{0} & \phantom{0} & \phantom{0} \\ \hline \phantom{0} & \phantom{0} \\ \cline{1-2} \end{array}$ から求めた。

**箱玉系の逆散乱形式**　$\tilde{e}_0$ を除く $A_n$ 部分に対する $B^{a,l}$ の最高ウェイト元を $u^{a,l}$ で表す。具体的に $u^{a,l}$ を Young 準標準盤として表すと、1 行目は全て文字 1、2 行目は全て文字 2、等々となる。一般のテンソル積 $b \in \bigotimes_{k=1}^{L} B^{r_k,s_k}$ と $u^{a,l}$ から出発して 23 ページの定義 (15) と同様にして箱玉系の時間発展作用素 $T^{a,l}$ が定義される [HHIKTT, FOY]。箱玉系は量子可積分系であり [FOY]、またソリトン系としてふるまうことも観察されている（箱玉系として最も基本的な $r_k = 1$ の場合は [Sa06] で一般的な状況で示されている）。

基本的な事実として、艤装配位写像を用いると、箱玉系の作用角変数の完全な系が得られる。より正確には：

**定理 [KOSTY]**

与えられた path $b$ に必要ならば右側に $u^{a,1}$ を十分付け加えて、$u^{a,l} \otimes b \overset{R}{\simeq} T^{a,l}(b) \otimes u^{a,l}$ となるようにしよう。その時 $b$ に対する艤装配位が

$$\Phi(b) = \left((\nu^{(1)}, J^{(1)}), \ldots, (\nu^{(a)}, J^{(a)}), \ldots, (\nu^{(n)}, J^{(n)})\right)$$

であれば、時間発展した path に対する艤装配位は

$$\Phi\left(T^{a,l}(b)\right) = \left((\nu^{(1)}, J^{(1)}), \ldots, \left(\nu_i^{(a)}, J_i^{(a)} + \min(\nu_i^{(a)}, l)\right)_i, \ldots, (\nu^{(n)}, J^{(n)})\right) \tag{29}$$

となる。すなわち Young 図 $\nu^{(1)}, \ldots, \nu^{(n)}$ は運動の保存量（作用変数）であり、rigging $J^{(a)}$ のみが線形に変化する（角変数）。

この事実を箱玉系の逆散乱形式と呼んでいる[26]。例は 39 ページ参照。なおこの結果より $l \geq \max_i \nu_i^{(a)}$ ならば $T^{a,l} = T^{a,l+1} = \cdots =: T^{a,\infty}$ となる。特に $T^{1,\infty} = T_\infty$ と書く。

この定理の本質は、下記の艤装配位の理論において最も深い定理の一つである。

**定理 [KSS]（$A_n^{(1)}$ 型艤装配位の $R$ 不変性）**

任意の $A_n^{(1)}$ 型 KR クリスタルのテンソル積 $b$ と $b'$ が $b \overset{R}{\simeq} b' \iff \Phi(b) = \Phi(b')$.

艤装配位写像 $\Phi$ は非常に複雑な写像であるので、この事実は極めて非自明であり、証明もとても難しい。

さて、$\Phi(b)$ と $\Phi(u^{a,l} \otimes b)$ とを比べると、28 ページの艤装配位 (20) において後者の方が $\mu^{(a-1)}$ の長さ $l$ の行が一つ多いが他の $\mu^{(a)}$ と $\nu^{(a)}$ は同一である（$\Phi(u^{a,l})$ は $\nu^{(a)}$ 部分が全て空集合）。その場合、§8.4 で与えられる写像 $\Phi^{-1}$ のアルゴリズムより、同じ $b$ を像の一部として得るためには rigging $J_i^{(a)}$ を $J_i^{(a)} + \min(\nu_i^{(a)}, l)$ に変化させればよい。あとは $u^{a,l} \otimes b \overset{R}{\simeq} T^{a,l}(b) \otimes u^{a,l}$ と艤装配位の $R$ 不変性から箱玉系の逆散乱形式が従う。

なお艤装配位の $R$ 不変性を 2 つのテンソル積の場合 $B \otimes B' \ni b \otimes b' \overset{R}{\simeq} \tilde{b}' \otimes \tilde{b}$ に適用すると、組み合わせ $R$ 行列が

$$\Phi_{B' \otimes B}^{-1} \circ \Phi_{B \otimes B'}(b \otimes b') = \tilde{b}' \otimes \tilde{b} \tag{30}$$

のようにして艤装配位の $R$ 不変性の特殊な場合として得られることになる。すなわち艤装配位写像は組み合わせ $R$ 行列の親玉のような存在である。

[26] 論文 [KOSTY] では著者間の見解の相違により明確な形でこの結果を述べることができなかった。そこでその直後に発表した論文 [KSY] の Theorem 3.5 において $a = 1$ の場合にここで述べた形の定理を与えた。本質的に大きな違いはないが、$a > 1$ の一般的な形は、例えば [Sa07] の Remark 3.5 で述べられている。

**箱玉系の数学的応用**　こうして艤装配位に関する数学的定理が箱玉系という数理物理の問題に有効に利用できることが分かった。逆に箱玉系との関連が明らかになったことにより下記のような数学的に重要な性質が発見された。筆者の印象としては、箱玉系が物理として応用される可能性はあまり高くないと思われるのだが、一方既に数学的な応用があるわけなので、その意味で箱玉系の研究は十分に意義があると考えている。

$b \in \bigotimes_{k=1}^{L} B^{r_k,s_k}$ とする。23 ページで箱玉系を定義した式 (14) の記号を変更した

$$\underline{u^{a,l} \otimes b_1} \otimes b_2 \otimes \cdots \otimes b_L \overset{R}{\simeq} b_1' \otimes \underline{v_1^{a,l} \otimes b_2} \otimes \cdots \otimes b_L \overset{R}{\simeq} b_1' \otimes b_2' \otimes \underline{v_2^{a,l} \otimes b_3} \cdots \otimes b_L \overset{R}{\simeq} \cdots \tag{31}$$

から出発しよう。以下では $v_0^{a,l} = u^{a,l}$ と書くことにする。[FOY] にならって和

$$E^{a,l}(b) = \sum_{k=1}^{L} H\left(v_{k-1}^{a,l} \otimes b_k\right) \tag{32}$$

を考える。その時、

定理 [Sa07]

艤装配位 $(\nu, J) = \Phi(b)$ に対して

$$Q_l(\nu^{(a)}) = E^{a,l}(b). \tag{33}$$

写像 $\Phi$ は複雑なアルゴリズムで定義されるので、この結果は大変非自明である。実際艤装配位の $R$ 不変性という深い結果をフルに使用して証明される。この等式は §8.6 で Loop Schur 関数との関連を考察するうえで鍵となる。

ここでは等式 (33) の別の応用として、艤装配位写像 $\Phi$ のアルゴリズムが実は代数的に記述できるものであることを指摘しておこう。Young 盤 $b_k$ の各列が左から $c_1, c_2, \ldots, c_{s_k}$ だったとしよう。その時途中段階のテンソル積として

$$C_j := b_1 \otimes b_2 \otimes \cdots b_{k-1} \otimes c_{s_k} \otimes c_{s_k-1} \otimes \cdots \otimes c_{j+1} \otimes c_j \tag{34}$$

を考えよう。$c_i$ の順番を入れ替えたのは現在の Kashiwara 流のテンソル積と Young 盤の相性が悪いためであるが、筆者は他の著者からの強い要請があった場合を除き、Kashiwara 流を使うことにしている（結晶基底のテンソル積の場合も異なる流儀が存在することによって時として混乱が生じてしまっているのは大変残念である）。この $C_j$ に対しても式 (32) に倣って和 $E^{a,l}(C_j)$ を定義できる。この時量

$$\left\{E^{a,l}(C_j) - E^{a,l-1}(C_j)\right\} - \left\{E^{a,l}(C_{j+1}) - E^{a,l-1}(C_{j+1})\right\} \tag{35}$$

は再帰的な定義を持つ $\Phi$ のアルゴリズムにおいて列 $c_j$ に対応する部分を計算するときに途中段階の $\nu^{(a)}$ の第 $l$ 列目に付け加えられるます目の総数を表している。なお式 (35) には重複が多く、整理すれば実際には局所的に決まる量である。実は [Sa07] で議論されているようにこの情報だけで rigging も含め必要な全ての情報を読み取ることができる。従って写像 $\Phi$ のアルゴリズムは、箱玉系の定義を用いると、結晶基底の言葉で記述できる代数的な性質のアルゴリズムであることが分かった。

写像 $\Phi$ のアルゴリズムを良く理解するとこれは驚くべき結果であると理解して頂けると思う。実際艤装配位の理論を数学化するのに重要な貢献をした Mark Shimozono 氏と議論した際も「あんな不思議なアルゴリズムが代数的な起源を持っていたとは衝撃的だ」と言って

おられたと記憶している。実際筆者もそれまでは艤装配位写像のアルゴリズムは精巧だが不思議なもの、といった印象を持っていた。しかしそれが実は無限次元の対称性に由来するものだったと分かり、艤装配位は単なる便利なモデルというよりそれ自身が数学的意味を持ち、従ってそれ自身を研究することに意味があると納得することができたと記憶している。

無限次元の対称性があると Young 盤も長方形のものしか現れないし一見特殊に感じる方もいらっしゃるかもしれないが、これは大雑把にいえば対称性が大きいために縛りがきつい事の反映である。数学の大海の中から特別に性質の良いものを見つけ出すのは、何の方針もなく探すのではあまり簡単なことではないであろう。そのような時無限次元の対称性のような強力な武器を使うと自然と性質の良いものが見つかる様に感じられる。

## 8.4 $A_n^{(1)}$ 型艤装配位写像のアルゴリズム

ここでは写像

$$\Phi^{-1} : \text{rigged configuration} \longmapsto \text{path}$$

のアルゴリズムを $A_n^{(1)}$ 型一般の場合に紹介する。すなわち path として最も一般的な $b \in \bigotimes_{k=1}^{L} B^{r_k,s_k}$ の場合を考える。28 ページの式 (20) のように艤装配位を

$$(\nu, J) = \left(\mu^{(0)}, \ldots, \mu^{(n-1)}, (\nu^{(1)}, J^{(1)}), \ldots, (\nu^{(n)}, J^{(n)})\right) \tag{36}$$

としよう。

**$\Phi^{-1}$ のアルゴリズム** $\Phi^{-1}$ の計算は $\mu^{(a)}$ $(0 \le a \le n-1)$ のどの行からスタートするのかを選ぶところから始まる。$\mu^{(a)}$ やその行の選び方を変えると異なる paths が得られるが、それらは艤装配位の $R$ 不変性（32 ページ）により全て $R$ 同型となる。$\mu_i^{(a)}$ を選んだとしよう。すると $B^{a+1,\mu_i^{(a)}}$ の元が以下のようにして得られる。ここでは艤装配位を Young 図を用いて図式的に表示して考えることにする。その時ストリングは Young 図の各行とその rigging の組と同一視される。

1. $\ell^{(a)}$ を $\mu_i^{(a)}$ の長さとする。$\nu^{(a+1)}, \nu^{(a+2)}, \ldots$ の特異な行（ストリング）を以下のようにして順に定める。$\ell^{(j)}$ まで定まっているとする。その時 $\nu^{(j+1)}$ に長さ $\ell^{(j)}$ 以上の特異な行が存在する場合、そのうち最も短いものを選び、その長さを $\ell^{(j+1)}$ として同じ手続きを $\nu^{(j+2)}$ に対して繰り返す。一方 $\nu^{(j_a-1)}$ まではそのような特異な行を選べたが $\nu^{(j_a)}$ には該当する特異な行が存在しなかった場合、文字 $j_a$ を出力として手続きを終了する。

2. 前項で選んだ各行の右端のます目を一つずつ削り、$\mu^{(a-1)}$ には長さ 1 の行を付け加えて新しい艤装配位を作る。新しい rigging は以下の通り。前項で選択されなかった行の rigging は変更しない。一方前項で選択され、ます目を一つ削られた行の rigging は、新しい艤装配位に対する vacancy number について特異な行となるように定める。

3. $\mu^{(a-1)}$ に付け加わった長さ 1 の行からスタートし、新たに $\ell^{(a-1)} = 1$ と定めて Step 1 と Step 2 を繰り返して出力 $j_{a-1}$ を得る。次に $\mu^{(a-2)}$ に付け加わった長さ 1 の行からスタートし、. . .、と繰り返し、最終的に $\mu^{(0)}$ に付け加わった長さ 1 の行からスタートし

て出力 $j_0$ を得る。その時 $j_a, \ldots, j_1, j_0$ を $(a+1) \times \nu_i^{(a)}$ 型 Young 盤の空きスペースの左端に下から順に書き込む。

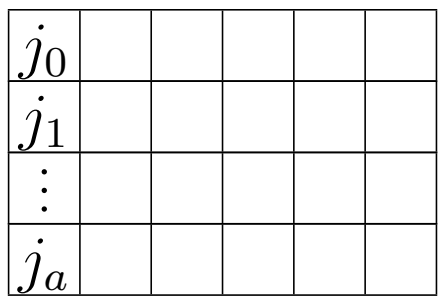

4. Step 1 から Step 3 を $\mu_i^{(a)}$（の残り）の右端からスタートして繰り返し、最終的に $\mu_i^{(a)}$ が全て取り去られるまで行う。すると上記長方形のます目がすべて埋まり、$B^{a+1,\mu_i^{(a)}}$ の元が一つ定まる。

残りの $\mu^{(0)}, \ldots, \mu^{(n-1)}$ の行についても上記 Step 1 から Step 4 を繰り返していき、全ての行が取り尽くされるまで続ける。得られた長方形を、テンソル積を挟んで右から左に並べたものを $b = \Phi^{-1}(\nu, J)$ とする。以上の手続きは矛盾なく定義されていることが知られている [KSS]。

**例**　以下の艤装配位を考える。

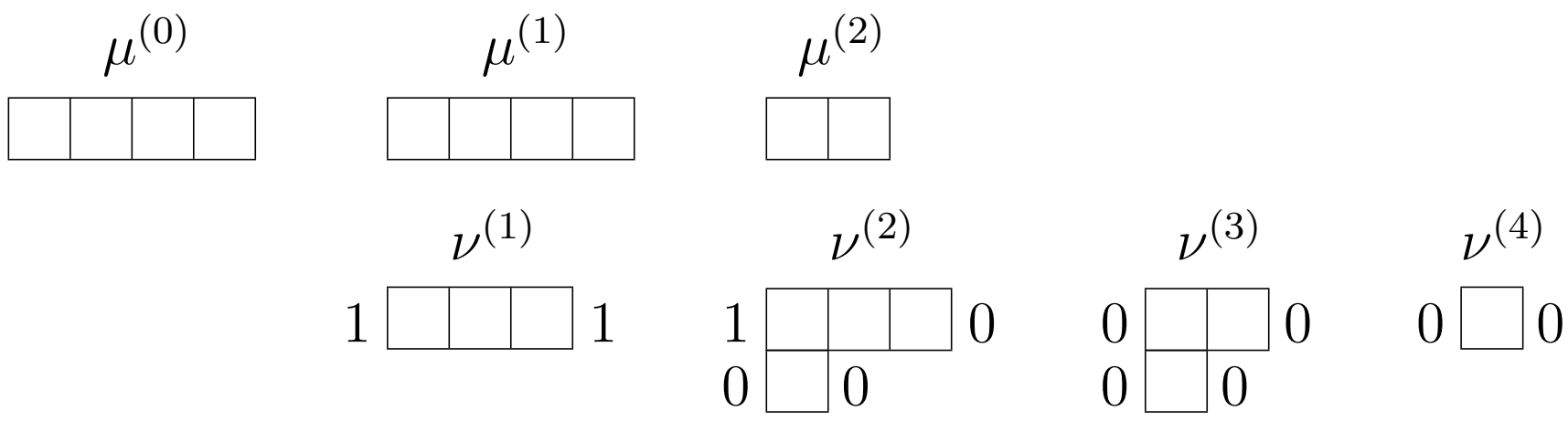

ここで riggings および vacancy number を対応する行のそれぞれ右と左に書き込んだ。最初に $\mu^{(1)}$ から削ることにすると、以下のように計算が進む（削られるます目に $\times$ を書き込んだ）。

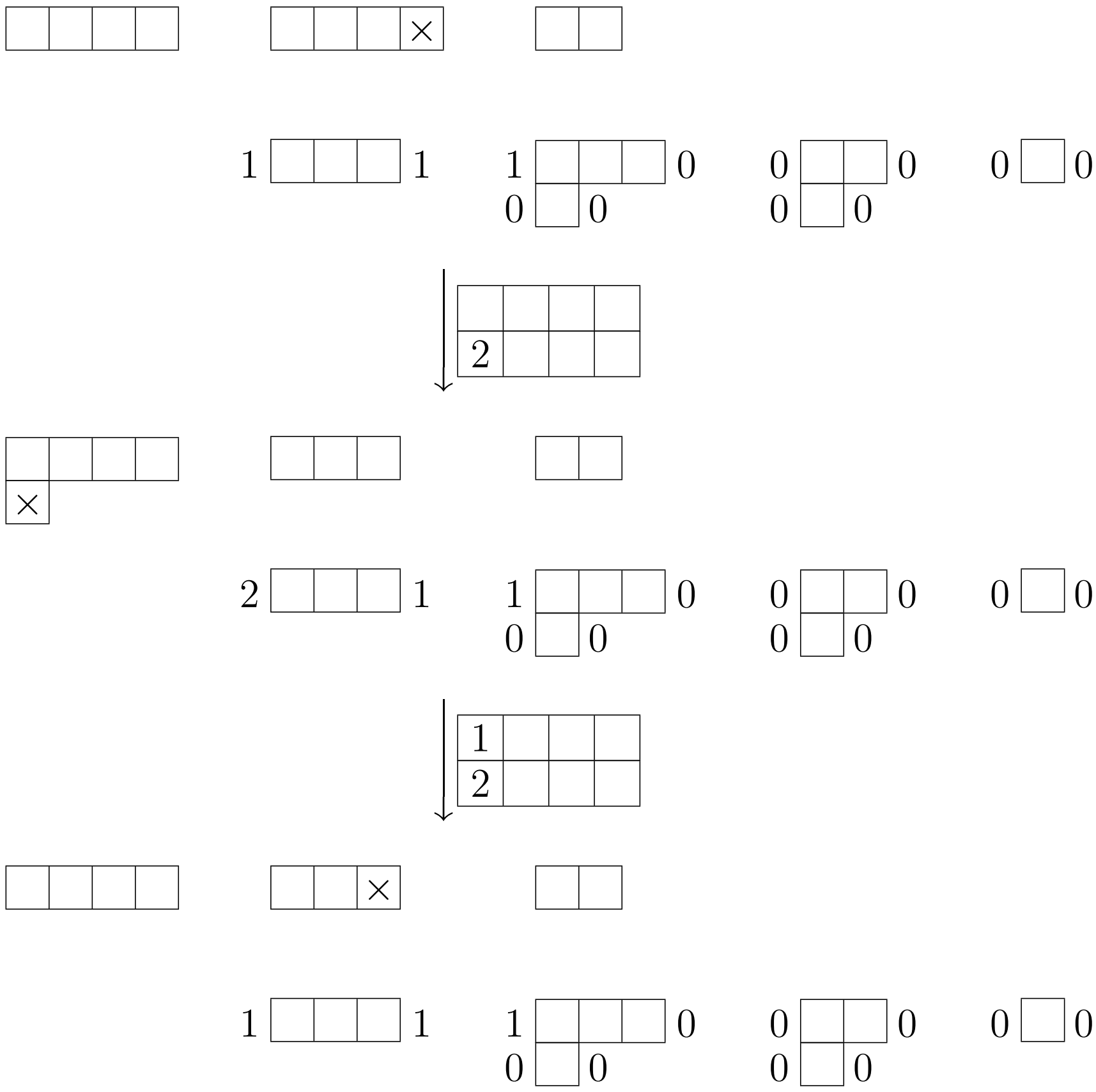

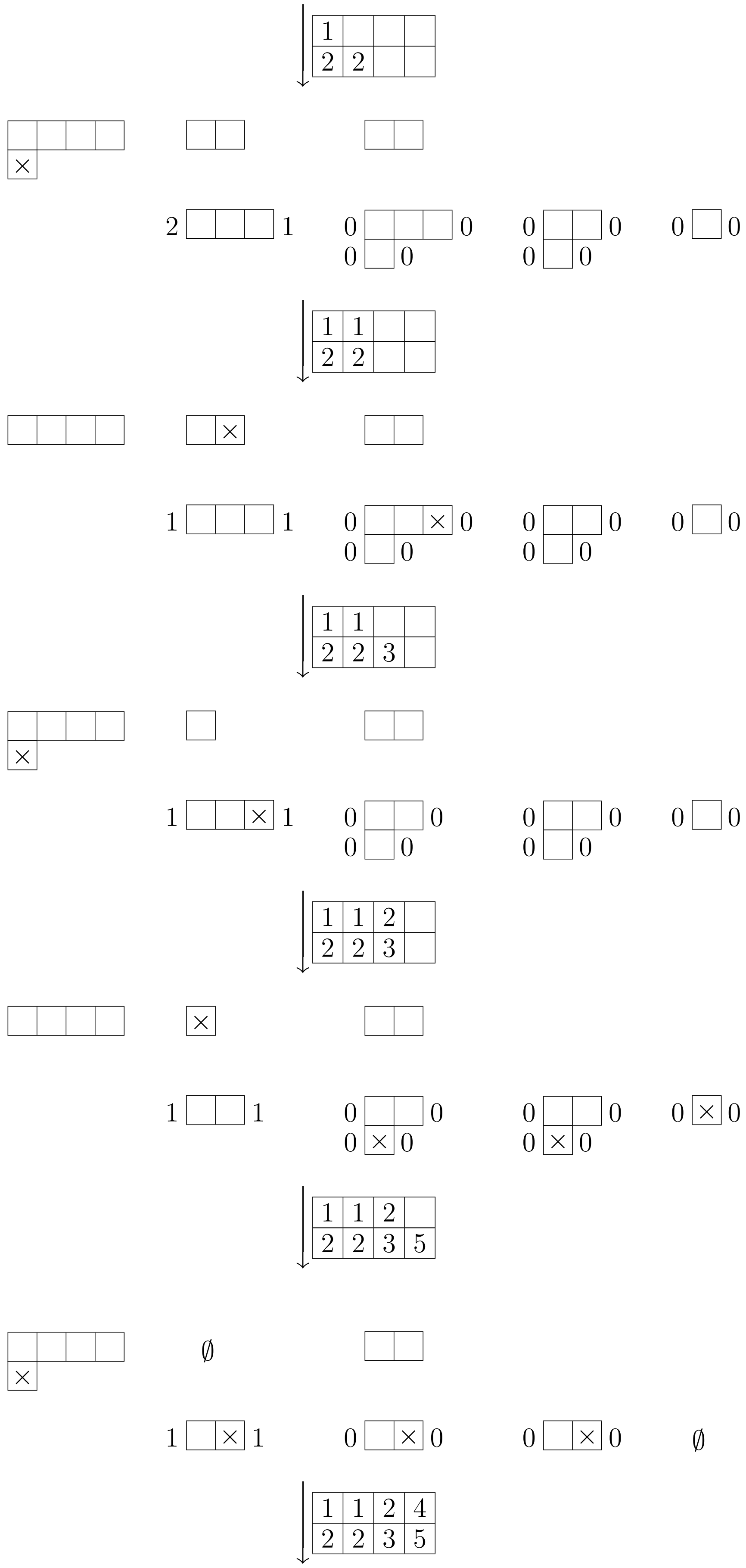

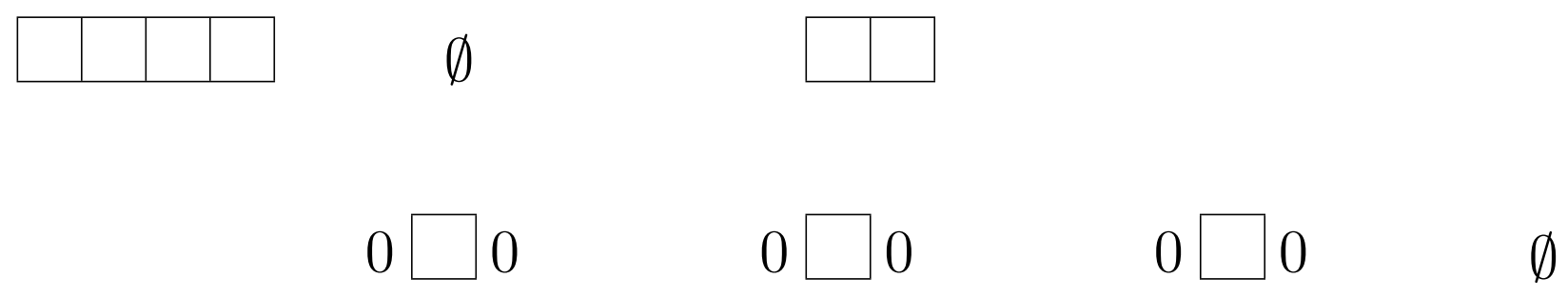

残りの計算も行うと以下の像が得られる。

$$\boxed{1\,|\,1\,|\,1\,|\,1} \otimes \begin{array}{|c|c|}\hline 1 & 2 \\ \hline 2 & 3 \\ \hline 3 & 4 \\ \hline \end{array} \otimes \begin{array}{|c|c|c|c|}\hline 1 & 1 & 2 & 4 \\ \hline 2 & 2 & 3 & 5 \\ \hline \end{array}.$$

**$\Phi$ のアルゴリズム**　基本的には上記 $\Phi^{-1}$ のアルゴリズムを逆向きに行えばよい。$A_1^{(1)}$ かつ $(B^{1,1})^{\otimes L}$ 型の場合のアルゴリズムは 22 ページに与えた。

一般の場合に最も基本となる操作は以下の通り。path の Young 盤の文字 $j$ を考えるとき、$\nu^{(j-1)}$ から以下の手続きを順に繰り返す（$\ell^{(j)} = \infty$ と仮定する)。$\nu^{(a)}$ の長さ $\ell^{(a)}$ の特異な行が選ばれたとき、$\nu^{(a-1)}$ の長さ $\ell^{(a)}$ 以下の特異な行のうち最大のものを選び、ます目を一つ加える。その様なものがなければ、$(0,0)$ という特異なストリングが存在すると考えて、長さ 1 の行を $\nu^{(a-1)}$ に付け加える。新しい rigging は、変更を受けた行については新しい艤装配位において特異な行となるように定め、他のものは元のままとする。

## 8.5　超離散タウ関数

$A_n^{(1)}$ 型 $\bigotimes_{k=1}^{L} B^{1,s_k}$ の場合、すなわち path の Young 盤が全て一行型の場合、$\Phi^{-1}$ の明示的公式が存在する [KSY]。32 ページで述べたように箱玉系の運動は艤装配位上で線形化されるので、$\Phi^{-1}$ の明示的公式と併せて箱玉系の初期値問題の完全な解が構成されたことになる。

**定式化**　具体的な式を記述するため本節内で使用する記号をいくつか用意しよう。与えられた分割 $\lambda = (\lambda_1, \ldots, \lambda_N)$ に対し

$$|\lambda| = \lambda_1 + \cdots + \lambda_N, \qquad \lambda_{[k]} = (\lambda_1, \ldots, \lambda_k) \quad (1 \leq k \leq N)$$

と書くことにしよう。更にもう一つの分割 $\mu = (\mu_1, \ldots, \mu_M)$ に対し

$$\min(\lambda, \mu) = \sum_{i=1}^{N} \sum_{j=1}^{M} \min(\lambda_i, \mu_j)$$

と書くことにする。

$(\nu, J)$ を与えられた艤装配位としよう。本節では艤装配位としては最高ウェイト元に限らず一般的なものを考えてよい。$A_n^{(1)}$ 型で最も一般的な艤装配位は 28 ページの式 (20) で与えたものであるが、現在は $\bigotimes_{k=1}^{L} B^{1,s_k}$ 型の path のみを考察しているので $\mu^{(1)} = \cdots = \mu^{(n-1)} = \emptyset$ とする。そこで $\mu^{(0)}$ を改めて $\nu^{(0)}$ と書くことにすれば、現在の設定での艤装配位は

$$(\nu, J) = \left(\nu^{(0)}, (\nu^{(1)}, J^{(1)}), \ldots, (\nu^{(n)}, J^{(n)})\right) \tag{37}$$

と表される。$\mu^{(a)}$ の添え字 $a$ をずらしていたのはこの様な状況を念頭に置いての事だったのだが、実は現在の状況では更に踏み込んで艤装配位に以下のような再帰的構造が存在する事が分かる。具体的には

$$(\nu, J)^{(a)} := \left(\nu^{(a)}, (\nu^{(a+1)}, J^{(a+1)}), \ldots, (\nu^{(n)}, J^{(n)})\right) \tag{38}$$

を艤装配位と考え $\Phi^{-1}$ を求める。その時得られた像と $J^{(a)}$ の情報とを合わせると結晶基底を用いた代数的な方法で

$$\Phi^{-1}\left(\nu^{(a-1)}, (\nu^{(a)}, J^{(a)}), \ldots, (\nu^{(n)}, J^{(n)})\right)$$

を決定することができる [Sa06]。特にこの時艤装配位の上にアフィン組み合わせ $R$ 行列 (24) の構造が入る。

艤装配位 $(\nu, J)$ の上にチャージ関数と呼ばれる量

$$c(\nu, J) = \frac{1}{2}\sum_{a,b} C_{a,b} \min\left(\nu^{(a)}, \nu^{(b)}\right) - \min\left(\nu^{(0)}, \nu^{(1)}\right) + \sum_{a,i} J_i^{(a)} \tag{39}$$

を定義する。ここで $(C_{ab})_{1\leq a,b\leq n}$ は $A_n$ 型の Cartan 行列である。数学的な背景として、いわゆる Kostka–Foulkes 多項式は艤装配位上のチャージ関数の母関数となっている [KSS]。なお、以下では艤装配位 $(\nu, J)$ を考えるとき $\nu^{(0)} = (\nu_1^{(0)}, \ldots, \nu_N^{(0)})$ の並び方は単調非増加に限らず任意の並び方を選ぶことにして、一つ固定しておこう。これは $\Phi^{-1}$ の像として $B^{1,\nu_1^{(0)}} \otimes \cdots \otimes B^{1,\nu_N^{(0)}}$ を考察することに対応する。

さて艤装配位を拡張して

$$(\nu, J)_{[k]} = \left((\nu^{(0)})_{[k]}, (\nu^{(1)}, J^{(1)}), \ldots, (\nu^{(n)}, J^{(n)})\right) \tag{40}$$

という量を考えよう。$\nu^{(0)}$ の長さが $N$ であれば $(\nu, J)_{[N]} = (\nu, J)$ である。更に

$$(\mu, I) \subseteq (\nu, J)_{[k]}$$

と書いた時、$(\nu, J)_{[k]}$ をストリング $(\nu_i^{(a)}, J_i^{(a)})$ の集合と考えた場合の部分集合であって、$(\nu^{(0)})_{[k]}$ は固定したものを考える。そうして得られる集合 $(\mu, I) \subseteq (\nu, J)_{[k]}$ の上にも形式的にチャージ関数 $c(\mu, I)$ を定義することにしよう。ここで**超離散タウ関数**を

$$\begin{aligned} \tau_{k,d} &= \max_{(\mu,I)\subseteq(\nu,J)_{[k]}} \left\{-c(\mu, I) - \left|\mu^{(d)}\right|\right\} \qquad (1 \leq d \leq n+1), \\ \tau_{k,0} &= \tau_{k,n+1} - \left|(\nu^{(0)})_{[k]}\right| \end{aligned} \tag{41}$$

と定義する。ただし $\nu^{(n+1)} = \emptyset$ と仮定した。この時、主要な結果は

**定理 [KSY]**

像 $b = \Phi^{-1}(\nu, J)$ の $k$ 番目の因子 $b_k$ に含まれる文字 $d$ の個数 $x_{k,d}$ は

$$x_{k,d} = \tau_{k,d} - \tau_{k-1,d} - \tau_{k,d-1} + \tau_{k-1,d-1}.$$

右辺は艤装配位 $(\nu, J)$ によって表される関数であるから、写像 $\Phi^{-1}$ を解析的に表示したことになっている。大きな式なのですぐには構造が分かりにくいかもしれないが、もっとも単純な $A_1^{(1)}$ 型 $(B^{1,1})^{\otimes L}$ の場合の式が 26 ページの脚注に記載されている。

**組み合わせ論的解釈**　超離散タウ関数の重要な性質として、箱玉系を用いた組み合わせ論的解釈をもつ。path $b = b_1 \otimes \cdots \otimes b_L \in \bigotimes_{k=1}^{L} B^{1,\nu_k^{(0)}}$ の Young 盤による表示において、文字 1 を空きスペース、文字 $2, 3, \ldots, n+1$ をそれぞれの文字でラベル付された玉と解釈する。$T_\infty^t(b) = b_1^{(t)} \otimes \cdots \otimes b_L^{(t)}$ と書く。path $b$ を 1 行目、$T_\infty(b)$ を 2 行目、$T_\infty^2(b)$ を 3 行目、…、と並べたものを考えるときに、$1 \le k \le L$ に対し関数 $\rho_{k,d}$ を以下のように定める。

$$\begin{aligned}\rho_{k,d} = &(b_1 \otimes b_2 \otimes \cdots \otimes b_k \text{における玉 } 2, 3, \ldots, d \text{ の総数}) \\ &+ \sum_{i \ge 1} (b_1^{(i)} \otimes b_2^{(i)} \otimes \cdots \otimes b_k^{(i)} \text{の全ての玉の個数}).\end{aligned} \tag{42}$$

箱玉系の時間発展は右向きに進行するので、第二項の和は有限和となる。

　さて超離散タウ関数は

$$\tau_{k,d} = \rho_{k,d} \tag{43}$$

という性質を持つ。この等式が示されれば、$\rho_{k,d}$ の意味を考えれば先ほどの定理はただちに従う。

　ここで $L = 20$ の $T_\infty$ 時間発展の例を見てみよう。以下テンソル積の記号は省略してスペースで表し、また Young 盤は内部の文字だけを並べて表すことにしよう。

| | | | | | | | | | | | | | | | | | | | | |
|---|---|---|---|---|---|---|---|---|---|---|---|---|---|---|---|---|---|---|---|---|
| $t=0:$ | 1 | 2 | 23 | 124 | 3 | 1111 | 3 | 1 | 1 | 1 | 1 | 1 | 1 | 1 | 1 | 1 | 1 | 1 | 1 | 1 |
| $t=1:$ | 1 | 1 | 12 | 123 | 2 | 1134 | 1 | 3 | 1 | 1 | 1 | 1 | 1 | 1 | 1 | 1 | 1 | 1 | 1 | 1 |
| $t=2:$ | 1 | 1 | 11 | 112 | 1 | 1223 | 4 | 1 | 3 | 3 | 1 | 1 | 1 | 1 | 1 | 1 | 1 | 1 | 1 | 1 |
| $t=3:$ | 1 | 1 | 11 | 111 | 2 | 1111 | 3 | 4 | 2 | 2 | 3 | 3 | 1 | 1 | 1 | 1 | 1 | 1 | 1 | 1 |
| $t=4:$ | 1 | 1 | 11 | 111 | 1 | 1112 | 1 | 3 | 1 | 1 | 2 | 2 | 4 | 3 | 3 | 1 | 1 | 1 | 1 | 1 |
| $t=5:$ | 1 | 1 | 11 | 111 | 1 | 1111 | 2 | 1 | 3 | 1 | 1 | 1 | 2 | 2 | 1 | 4 | 3 | 3 | 1 | 1 |

$\rho_{k,d}$ ないし $\tau_{k,d}$ の値は以下の通り。

| $k$ | 1 | 2 | 3 | 4 | 5 | 6 | 7 | 8 | 9 | 10 | 11 | 12 | 13 | 14 | 15 | 16 | 17 | 18 | 19 | 20 |
|---|---|---|---|---|---|---|---|---|---|---|---|---|---|---|---|---|---|---|---|---|
| $\rho_{k,1}$ | 0 | 0 | 1 | 4 | 6 | 12 | 15 | 19 | 23 | 27 | 31 | 35 | 39 | 43 | 47 | 51 | 55 | 59 | 63 | 67 |
| $\rho_{k,2}$ | 0 | 1 | 3 | 7 | 9 | 15 | 18 | 22 | 26 | 30 | 34 | 38 | 42 | 46 | 50 | 54 | 58 | 62 | 66 | 70 |
| $\rho_{k,3}$ | 0 | 1 | 4 | 8 | 11 | 17 | 21 | 25 | 29 | 33 | 37 | 41 | 45 | 49 | 53 | 57 | 61 | 65 | 69 | 73 |
| $\rho_{k,4}$ | 0 | 1 | 4 | 9 | 12 | 18 | 22 | 26 | 30 | 34 | 38 | 42 | 46 | 50 | 54 | 58 | 62 | 66 | 70 | 74 |

なお各 $t$ に対する艤装配位を計算すれば以下のようになる（vacancy number は略）。

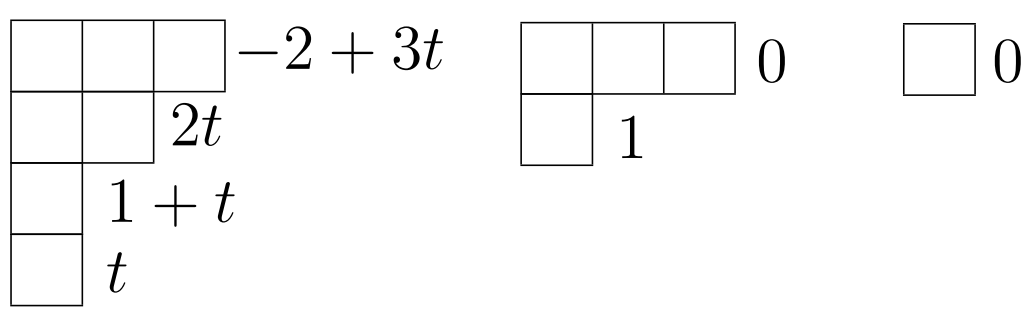

これは 32 ページで述べた箱玉系の逆散乱形式の例となっている。

**証明の概略**　上述のように主要な目標は等式 (43) の証明である。まず箱玉系の玉と運搬車による時間発展の解釈 [HHIKTT] を用いると、箱玉系の運動を表す $\rho_{k,d}$ が

$$\overline{\rho}_{k,d-1} + \rho_{k-1,d} = \max(\overline{\rho}_{k,d} + \rho_{k-1,d-1},\ \overline{\rho}_{k-1,d-1} + \rho_{k,d} - \nu_k^{(0)}) \tag{44}$$

という超離散 Hirota 双線形形式を満たすこと [KSY, Proposition 4.2] に注意する。ここで $\overline{\rho}_{k,d}$ とは $T_\infty(b)$ に対する $\rho_{k,d}$ の事を表す。

続いて $\tau_{k,d}$ も同一の方程式を満たすこと

$$\overline{\tau}_{k,d-1} + \tau_{k-1,d} = \max(\overline{\tau}_{k,d} + \tau_{k-1,d-1},\ \overline{\tau}_{k-1,d-1} + \tau_{k,d} - \nu_k^{(0)}) \tag{45}$$

を証明する [KSY, §5]。ここで $\overline{\tau}_{k,d}$ は $\rho_{k,d}$ の場合と同様の記号法であり、箱玉系の逆散乱形式により rigging を $J_i^{(1)} \mapsto J_i^{(1)} + \nu_i^{(1)}$ と変化させた $\tau_{k,d}$ の事を表す。この式を直接証明するのは難しいので、少々技巧的な方法を用いる。すなわち超離散極限で定義 (41) を与えるような行列式を見つけてくることができ（43 ページに補足あり）、後は本質的には行列式の非自明な恒等式（KP 階層の理論 [JM] で現れる双線形形式）に証明を帰着させる。なお今の所この様にしてソリトン系と関連付けられたという事実自体が数学的に意味のある結果なのか、あるいはただの計算上の方便なのかははっきりとしていない。

以上式 (44) と式 (45) が示されたことにより、$\rho_{k,d}$ と $\tau_{k,d}$ は同じ箱玉系の力学系に従うことが分かり、従って系がもっとも単純化する十分時間発展した後の状態 $T_\infty^S(b)$ $(S \gg 1)$ について $\tau_{k,d} = \rho_{k,d}$ を確かめればよいことが分かった（この様な状態を**漸近状態**と呼ぼう）。写像 $\Phi^{-1}$ の定義は組み合わせ論的なものであり、そのままでは何も進まないので、ここから先は本格的に組み合わせ論の深い議論が必要となる場面である。

まず path $b$ が最高ウェイト元である場合には艤装配位のはっきりとした特徴づけが存在するのでその場合の漸近状態について $\tau_{k,d} = \rho_{k,d}$ を確かめる。証明は代数 $A_n^{(1)}$ のランク $n$ による帰納法による。$A_1^{(1)}$ の場合は [Sa06] における艤装配位写像とエネルギー関数の関係を用いて直接示せる。一般の場合を考察するために、まず補助的にアフィンクリスタルに対する Yang–Baxter 関係式を用いて $\rho_{k,d}$ のエネルギー関数を用いた表示（[KSY, §4] 参照）を準備する。さて、帰納法を走らせるために式 (38) 周辺で指摘した艤装配位の帰納的構造を用いる。この時 $(\nu, J)^{(a)}$ に対して定義したタウ関数を $\tau_{k,d}^{(a)}$ と書くことにすれば、タウ関数に対しても帰納的構造

$$\tau_{k,d}^{(a-1)} = \max_{(\mu,I)\subseteq(\nu^{(a)},J^{(a)})} \left\{ \min(\nu_{[k]}^{(a-1)}, \mu) - \min(\mu,\mu) - \sum_i I_i^{(a)} + \tau_{\ell(\mu),d}^{(a)}(\mu, I) \right\}, \tag{46}$$

ここで $\ell(\mu)$ は分割 $\mu$ の長さ、が成り立つことに注意する。この式を漸近状態について詳しく解析すると、[Sa06] で艤装配位上に導入されたエネルギー関数の構造を絶妙に利用することができ、帰納法の仮定を用いて証明が完成する（かなめは [KSY, Lemma 6.6]）。元々論文 [Sa06] はここでの使用を念頭に置いて執筆されたもので、艤装配位上にエネルギー関数の構造を導入するにはかなりハードな組み合わせ論的証明が要求される。

最終的に最高ウェイト元に限らない一般のテンソル積に対して $\tau_{k,d} = \rho_{k,d}$ を証明するには、明らかに技巧的なやり方で（[KSY, §7] 参照）最高ウェイト元の場合から一般の結果を切り出して証明できる。箱玉系の運動の特性を巧妙に使用した議論となる。

**証明に関するコメント**　箱玉系の初期値問題の解決は当時懸案の難問となっていたが、逆散乱形式（32 ページ）にみられるように、一見目で見ることのできるソリトン $(\nu^{(1)}, J^{(1)})$ だけでなく、一段奥まったところにある他のソリトン $(\nu^{(a)}, J^{(a)})$ まで考慮しなければならない点が難関だったのだろうと思われる。実際上記証明中で [JM] におけるソリトン理論と比較した際もそれらすべてのソリトンが独立した形で寄与してくる。艤装配位写像のような強力な理論的枠組みがなければとても解けなかったのではないかと思っている。

## 8.6 Cylindric Loop Schur 関数

$A^{(1)}_{n-1}$ 型 $\bigotimes_{k=1}^{L} B^{1,s_k}$ の場合、$\Phi$ によって得られる分割 $\nu^{(a)}$ は Cylindric Loop Schur 関数によって記述できることが知られている [LPS2]。この節の内容について Mathematica によるプログラム例が

https://sites.google.com/site/affinecrystal/rigged-configurations

にて入手可能ですのでご利用ください。

以下本節では超離散極限における変数を $x_i$、対応する連続関数の世界における変数を $\mathrm{x}_i$ などと区別して表示することにする。また理解しやすくするために文字 $n$ と $L$ はそれぞれ $A^{(1)}_{n-1}$（ないし $\widehat{\mathfrak{sl}}_n$）と $\bigotimes_{k=1}^{L} B^{1,s_k}$ に出てくるものとして固定しておこう。

**Loop 対称関数**　Loop 対称関数は組み合わせ $R$ 行列の連続関数版に対する不変式の研究に伴って Yamada [Y00] によって導入された Loop 基本対称式と呼ばれるものを嚆矢とする一連の関数のクラスである。Loop 対称関数は変数 $\mathrm{x}_i^{(j)}$ に対する関数であり、上付きの添え字 $\bullet^{(j)}$ を $j \in \mathbb{Z}/n\mathbb{Z}$ の元であるとみなすところが特徴である。

その後の Lam–Pylyavskyy による研究の解説は [L] に与えられている。特に彼らの未公表の結果 [L, Theorem 4.4] によれば全ての $R$ 不変量は Loop 基本対称式により生成される環 LSym と一致する。艤装配位も $R$ 不変量であるから原理的には Loop 対称関数による表示を持つはずである。しかし「原理的に表示可能」と「実際に表示可能」は別物であるし、むしろその差が大きければ大きいほど数学としては興味深い問題と言えるだろう。実際以下で見るように艤装配位のうち $\nu^{(a)}$ の部分についてはまずまず納得のいく結果が得られているが、rigging $J^{(a)}$ の方は難しく、時間発展を十分させた後の状態 $T^S_\infty(b)$ $(S \gg 1)$ —40 ページで漸近状態と呼んだもの—について部分的な結果が得られている [Scr16] のみである。そしてそれらは既に（33 ページの式 (33) や 40 ページの議論で）示したように、艤装配位の理論の中では単純化することが知られている領域であるし、かなり正体のはっきりしている部分でもある。残る rigging の部分から興味深い数学が見つかれば良いと思う（§6 における Bethe 仮説方程式の解の解析でもその部分に興味深い構造が隠されていたのであった）。

艤装配位との関連では [LPS2] によって導入された Cylindric Loop Schur 関数が重要である。定義は通常の Schur 関数に対するタブローを用いた定義と並行したやり方でなされる。以下順を追って導入して行こう。$r \in \mathbb{Z}/n\mathbb{Z}$ を一つ選んだ時、タブローの各文字 $T(i,j)$ と変数との対応を以下の規則によって定める。

$$
\begin{array}{|c|c|c|c|c}
\hline
T(1,1) & T(1,2) & T(1,3) & T(1,4) & \cdots \\ \hline
T(2,1) & T(2,2) & T(2,3) & T(2,4) & \cdots \\ \hline
T(3,1) & T(3,2) & T(3,3) & T(3,4) & \cdots \\ \hline
T(4,1) & T(4,2) & T(4,3) & T(4,4) & \cdots \\ \hline
\vdots & \vdots & \vdots & \vdots & \ddots
\end{array}
\longmapsto
\begin{array}{|c|c|c|c|c}
\hline
\mathrm{x}^{(r)}_{T(1,1)} & \mathrm{x}^{(-1+r)}_{T(1,2)} & \mathrm{x}^{(-2+r)}_{T(1,3)} & \mathrm{x}^{(-3+r)}_{T(1,4)} & \cdots \\ \hline
\mathrm{x}^{(1+r)}_{T(2,1)} & \mathrm{x}^{(r)}_{T(2,2)} & \mathrm{x}^{(-1+r)}_{T(2,3)} & \mathrm{x}^{(-2+r)}_{T(2,4)} & \cdots \\ \hline
\mathrm{x}^{(2+r)}_{T(3,1)} & \mathrm{x}^{(1+r)}_{T(3,2)} & \mathrm{x}^{(r)}_{T(3,3)} & \mathrm{x}^{(-1+r)}_{T(3,4)} & \cdots \\ \hline
\mathrm{x}^{(3+r)}_{T(4,1)} & \mathrm{x}^{(2+r)}_{T(4,2)} & \mathrm{x}^{(1+r)}_{T(4,3)} & \mathrm{x}^{(r)}_{T(4,4)} & \cdots \\ \hline
\vdots & \vdots & \vdots & \vdots & \ddots
\end{array}
\qquad (47)
$$

右辺に現れる変数の積を $\mathrm{x}^{\mathrm{wt}^{(r)}(T)}$ と書くことにすれば、**Loop Schur 関数**とは

$$
s^{(r)}_{\lambda/\mu}(\mathrm{x}) = \sum_{T} \mathrm{x}^{\mathrm{wt}^{(r)}(T)} \qquad (48)
$$

で定義される。ここで和は文字 $1, 2, \ldots, L$ による $\lambda/\mu$ 型 Young 準標準盤の全体をわたる。

**例**　$n=3$ とする。$\lambda=(2,1)$ および $L=3$ の場合以下のタブローが存在する。

$$\begin{array}{cccc} \begin{array}{|c|c|}\hline 1 & 1 \\ \hline 2 & \\ \cline{1-1} \end{array} & \begin{array}{|c|c|}\hline 1 & 2 \\ \hline 2 & \\ \cline{1-1} \end{array} & \begin{array}{|c|c|}\hline 1 & 2 \\ \hline 3 & \\ \cline{1-1} \end{array} & \begin{array}{|c|c|}\hline 1 & 3 \\ \hline 2 & \\ \cline{1-1} \end{array} \\ \begin{array}{|c|c|}\hline 1 & 1 \\ \hline 3 & \\ \cline{1-1} \end{array} & \begin{array}{|c|c|}\hline 2 & 2 \\ \hline 3 & \\ \cline{1-1} \end{array} & \begin{array}{|c|c|}\hline 1 & 3 \\ \hline 3 & \\ \cline{1-1} \end{array} & \begin{array}{|c|c|}\hline 2 & 3 \\ \hline 3 & \\ \cline{1-1} \end{array} \end{array} \tag{49}$$

$r=1$ とした上で対応 (47) を用いると、各 Young 盤から以下の多項式が得られる。

$$\begin{aligned} s^{(1)}_{2,1}(\mathrm{x}_1,\mathrm{x}_2,\mathrm{x}_3) = {} & \mathrm{x}^{(1)}_1\mathrm{x}^{(3)}_1\mathrm{x}^{(2)}_2 + \mathrm{x}^{(1)}_1\mathrm{x}^{(3)}_2\mathrm{x}^{(2)}_2 + \mathrm{x}^{(1)}_1\mathrm{x}^{(3)}_2\mathrm{x}^{(2)}_3 + \mathrm{x}^{(1)}_1\mathrm{x}^{(3)}_3\mathrm{x}^{(2)}_2 + \\ & \mathrm{x}^{(1)}_1\mathrm{x}^{(3)}_1\mathrm{x}^{(2)}_3 + \mathrm{x}^{(1)}_2\mathrm{x}^{(3)}_2\mathrm{x}^{(2)}_3 + \mathrm{x}^{(1)}_1\mathrm{x}^{(3)}_3\mathrm{x}^{(2)}_3 + \mathrm{x}^{(1)}_2\mathrm{x}^{(3)}_3\mathrm{x}^{(2)}_3, \end{aligned}$$

$\mathrm{x}^{(1)}_i=\mathrm{x}^{(2)}_i=\mathrm{x}^{(3)}_i=\mathrm{x}_i$ と置けば通常の Schur 多項式 $s_{2,1}(\mathrm{x}_1,\mathrm{x}_2,\mathrm{x}_3)$ が得られる。

**Cylindric Loop Schur 関数**　Loop Schur 関数を以下のように拡張する。二つの分割 $\lambda$ と $\mu$ から歪分割 $\lambda/\mu$ を作る。その時 $\mathscr{D}_a(\lambda/\mu)$ とは分割 $\lambda/\mu$ を 2 次元平面に配置したと考え、$(n-a,a)$ ずつシフトさせていって無限につなげた円柱状分割とする。同様にしてタブローに対しても円柱版を考えるが、そうしてできたタブローの各行、各列が準標準準盤の条件を満たす時、円柱状タブローと呼ぶことにしよう。対応 (47) においてシフト $(n-a,a)$ は添え字 $\bullet^{(j)}$ に影響しないので、Loop Schur 関数 (48) と同様にして円柱状タブローについても多項式が定まり、**Cylindric Loop Schur 関数**とよぶ。もし $\lambda_1<n-a$ であれば円柱状分割は通常の分割の非連結和となることから $s^{(r)}_{\mathscr{D}_a(\lambda/\mu)}=s^{(r)}_{\lambda/\mu}$ となる。

重要な性質として、Cylindric Loop Schur 関数は全て $R$ 不変量となっている [LPS2, Theorem 4.4]。

**例**　$n=3$ および $a=1$ の時シフトとしては $(n-a,a)=(2,1)$ を考える。分割 $\lambda=(2,1)$ を考えると対応する円柱状分割 $\mathscr{D}_1(\lambda)$ は以下の通り：

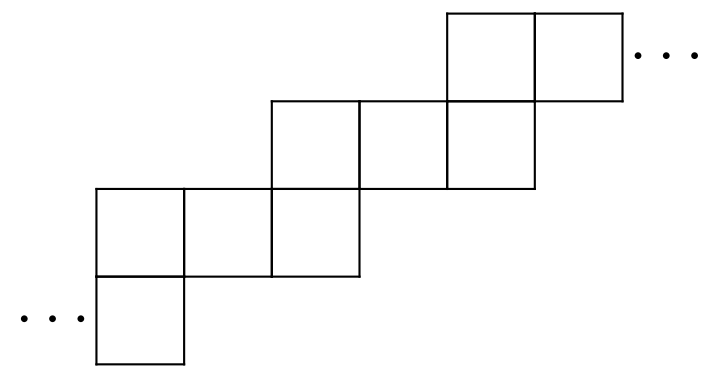

先ほどの例と同じ状況で考えよう。$L=3$ とすると、Cylindric Loop Schur 関数は以下の通り。

$$\begin{aligned} s^{(1)}_{\mathscr{D}_1(2,1)}(\mathrm{x}_1,\mathrm{x}_2,\mathrm{x}_3) = {} & \mathrm{x}^{(1)}_1\mathrm{x}^{(3)}_1\mathrm{x}^{(2)}_2 + \mathrm{x}^{(1)}_1\mathrm{x}^{(3)}_2\mathrm{x}^{(2)}_2 + \mathrm{x}^{(1)}_1\mathrm{x}^{(3)}_2\mathrm{x}^{(2)}_3 + \\ & \mathrm{x}^{(1)}_1\mathrm{x}^{(3)}_1\mathrm{x}^{(2)}_3 + \mathrm{x}^{(1)}_2\mathrm{x}^{(3)}_2\mathrm{x}^{(2)}_3 + \mathrm{x}^{(1)}_1\mathrm{x}^{(3)}_3\mathrm{x}^{(2)}_3 + \mathrm{x}^{(1)}_2\mathrm{x}^{(3)}_3\mathrm{x}^{(2)}_3. \end{aligned}$$

円柱状ではない場合 (49) にあげた 8 個の Young 準標準準盤が存在する。そのうち最初の 3 つをシフト $(2,1)$ により円柱状に伸ばすと、以下のように準標準準盤の条件を満たす。

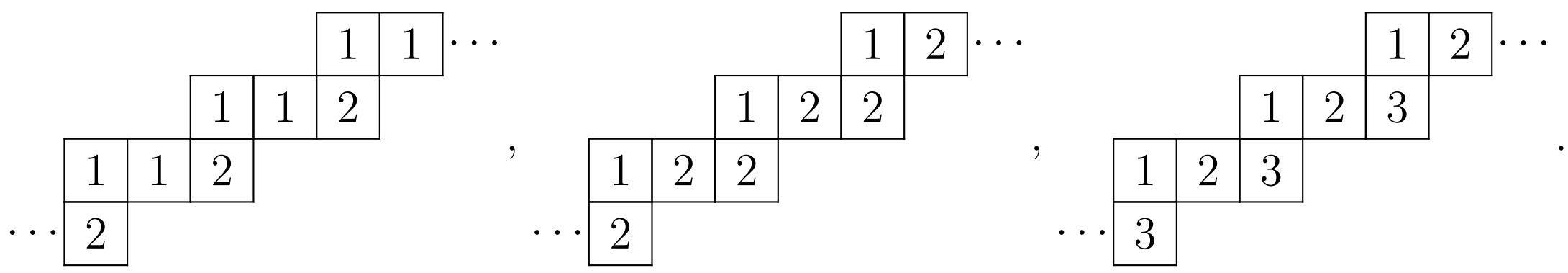

しかしながら 4 つ目のタブローは以下のように円柱状タブローとはならない。

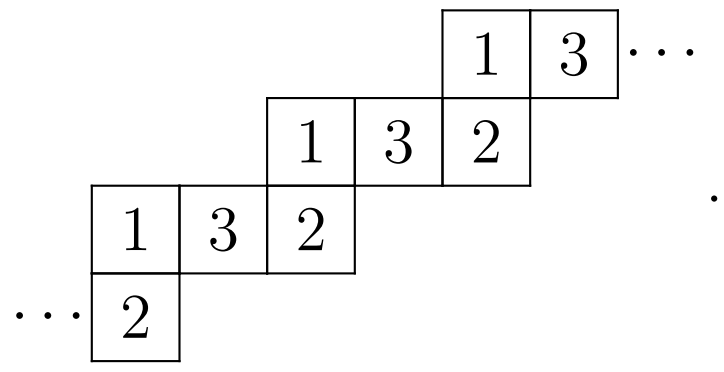

従ってこの項は $s^{(1)}_{\mathscr{D}_1(2,1)}$ には寄与しない。残る 4 つのタブローは $s^{(1)}_{\mathscr{D}_1(2,1)}$ に寄与することも同様にして確かめることができる。

**艤装配位との関連**　結果を定式化するためにいくつか準備をしよう。テンソル積 $b = b_1 \otimes \cdots \otimes b_L \in B^{1,s_1} \otimes \cdots \otimes B^{1,s_L}$ をタブロー表示したときに $x_j^{(i+j-1)}$ で $b_{L+1-j}$ における文字 $i$ の個数を表すものとする。例えば $n=3$ かつ $L=4$ の時

$$\boxed{1^{x_4^{(1)}}2^{x_4^{(2)}}3^{x_4^{(3)}}} \otimes \boxed{1^{x_3^{(3)}}2^{x_3^{(1)}}3^{x_3^{(2)}}} \otimes \boxed{1^{x_2^{(2)}}2^{x_2^{(3)}}3^{x_2^{(1)}}} \otimes \boxed{1^{x_1^{(1)}}2^{x_1^{(2)}}3^{x_1^{(3)}}}$$

となる。

次に $1 \le a \le (n-1)$ に対し分割 $\lambda(a,i)$ を以下のように再帰的に定める。最初に $\lambda(a,0) = (n-a)^L$、すなわち長さ $(n-a)$ の行 $L$ 本からなる長方形とする。$\lambda(a,i)$ まで定まったとき、$\lambda(a,i+1)$ は $\lambda(a,i)$ の最も左下のます目からスタートして、右または上向きに $\lambda(a,i)$ の外縁のます目を最大 $n$ 個削除して得られる分割とする。$i$ が十分大きくなると $\lambda(a,i)$ は空集合となる。例えば $n=6,\ a=2$ および $L=7$ の時は以下のようになる。

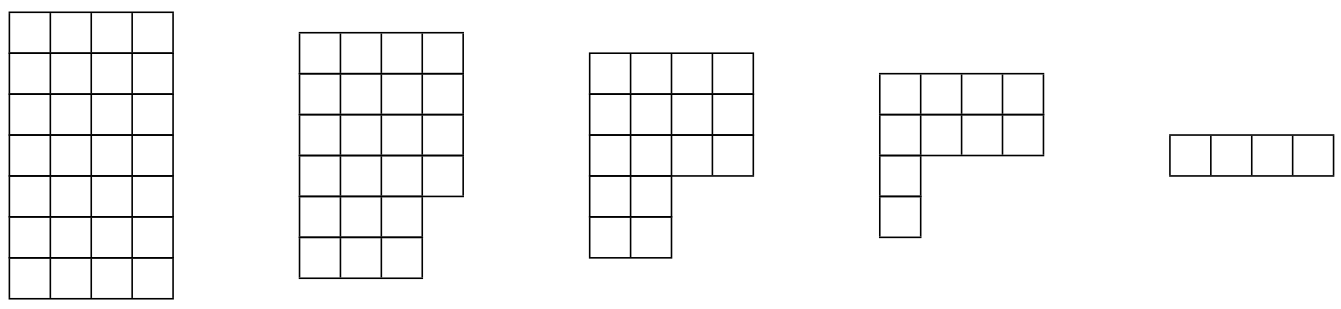

一方 $n=6, a=3$ および $L=3$ の時は各ステップで必ずしも $n$ 個のます目を取ることはできず、以下のようになる。

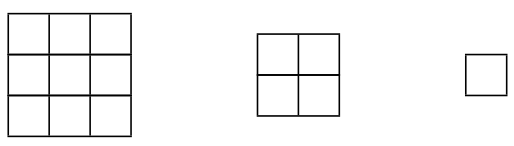

最後に有理関数 $P(\mathrm{x}_1,\mathrm{x}_2,\ldots,\mathrm{x}_r)$ の**超離散化**ないし**トロピカル化** $\mathrm{trop}(P)(x_1,x_2,\ldots,x_r)$ とは有理関数 $P(\mathrm{x}_1,\mathrm{x}_2,\ldots,\mathrm{x}_r)$ の中で形式的に $+ \mapsto \min$,、$\times \mapsto +$ および $\div \mapsto -$ と置き換えることとする。この極限の概念的問題については既に 26 ページから始まる項目『超離散極限について』で詳しく検討したが、ここでは技術的側面について若干補足しておこう。今与えた定式化では $P(\mathrm{x}_1,\mathrm{x}_2,\ldots,\mathrm{x}_r)$ は負符号が含まれない正整数係数多項式の商であることが望ましいが、対数と極限を用いた理解（25 ページ式 (16) 参照）によればこれは本質的な問題ではなく、特異性の問題を何らかの形で回避できればより一般的な係数を持っていても良いことになる。実際 40 ページでの議論はそのようなものであり、結果として行列式の超離散極限を取っている。式 (16) の場合も厳密にいえば正整数係数有理関数の範疇を超えている。よって超離散化の問題を考察するときに頭から正整数係数有理関数であることを仮定する訳にはいかない。なお、当然のことながら超離散化の過程で係数の情報は消し去られてしまっているので、超離散化された式から元の式を復元することは不可能である。

以上の準備の下で結果を述べよう。$\Phi(b)$ に現れる分割 $\nu^{(a)} = (\nu_1^{(a)}, \nu_2^{(a)}, \ldots)$ は

$$\nu_i^{(a)} = \mathrm{trop}\left(\frac{s^{(0)}_{\mathscr{D}_a(\lambda(a,i-1))}}{s^{(0)}_{\mathscr{D}_a(\lambda(a,i))}}\right)(x_l^{(k)}) \tag{50}$$

で与えられる。$a=1$ の場合これは定理であり [LPS2, Theorem 6.1]、$a>1$ の場合は予想である [LPS2, Conjecture 5.3]。証明は 33 ページの式 (33) で与えた $\nu^{(a)}$ のエネルギー関数による表示と [LP] によるエネルギー関数の Loop Schur 関数による表示とを用いる。

$a=1$ の場合の定理 (50) から得られる興味深い結果の一つは、$A^{(1)}_{n-1}$ 型の長さ $L$ の path に含まれるソリトンの個数は最大

$$\left\lceil \frac{(n-1)L}{n} \right\rceil$$

である。ここで天井関数 $\lceil N \rceil$ とは $N$ 以上の最小の整数とする。

**例** $a=2$ の場合の予想 (50) の例を考えよう。具体的には $n=4$ および $L=4$ の場合を考える。この時 $\lambda(2,0), \lambda(2,1), \ldots$ は以下の通り。

$$\begin{array}{|c|c|}\hline \phantom{1} & \phantom{1} \\ \hline \phantom{1} & \phantom{1} \\ \hline \phantom{1} & \phantom{1} \\ \hline \phantom{1} & \phantom{1} \\ \hline \end{array} \qquad \begin{array}{|c|c|}\hline \phantom{1} & \phantom{1} \\ \hline \phantom{1} \\ \cline{1-1} \phantom{1} \\ \cline{1-1} \end{array} \qquad \emptyset \qquad \emptyset \qquad \cdots$$

ここから $\mathscr{D}_2(\lambda(a,i))$ を構成するにはシフト $(n-a,a)=(2,2)$ を考えればよい。文字 $1,\ldots,L$ によるタブローを考えるので、$s^{(0)}_{\mathscr{D}_2(2,2,2,2)}$ は一項のみである。

$$s^{(0)}_{\mathscr{D}_2(2,2,2,2)} = \mathrm{x}^{(4)}_1\mathrm{x}^{(3)}_1\mathrm{x}^{(1)}_2\mathrm{x}^{(4)}_2\mathrm{x}^{(2)}_3\mathrm{x}^{(1)}_3\mathrm{x}^{(3)}_4\mathrm{x}^{(2)}_4.$$

次に $s^{(0)}_{\mathscr{D}_2(2,1,1)}$ には以下の 14 個のタブローが該当するが

$$\begin{array}{|c|c|}\hline 1 & 1 \\ \hline 2 \\ \cline{1-1} 3 \\ \cline{1-1} \end{array} \quad \begin{array}{|c|c|}\hline 1 & 1 \\ \hline 2 \\ \cline{1-1} 4 \\ \cline{1-1} \end{array} \quad \begin{array}{|c|c|}\hline 1 & 1 \\ \hline 3 \\ \cline{1-1} 4 \\ \cline{1-1} \end{array} \quad \begin{array}{|c|c|}\hline 1 & 2 \\ \hline 2 \\ \cline{1-1} 3 \\ \cline{1-1} \end{array} \quad \begin{array}{|c|c|}\hline 1 & 2 \\ \hline 2 \\ \cline{1-1} 4 \\ \cline{1-1} \end{array} \quad \begin{array}{|c|c|}\hline 1 & 2 \\ \hline 3 \\ \cline{1-1} 4 \\ \cline{1-1} \end{array} \quad \begin{array}{|c|c|}\hline 1 & 3 \\ \hline 2 \\ \cline{1-1} 3 \\ \cline{1-1} \end{array}$$

$$\begin{array}{|c|c|}\hline 1 & 3 \\ \hline 2 \\ \cline{1-1} 4 \\ \cline{1-1} \end{array} \quad \begin{array}{|c|c|}\hline 1 & 3 \\ \hline 3 \\ \cline{1-1} 4 \\ \cline{1-1} \end{array} \quad \begin{array}{|c|c|}\hline 1 & 4 \\ \hline 2 \\ \cline{1-1} 4 \\ \cline{1-1} \end{array} \quad \begin{array}{|c|c|}\hline 1 & 4 \\ \hline 3 \\ \cline{1-1} 4 \\ \cline{1-1} \end{array} \quad \begin{array}{|c|c|}\hline 2 & 2 \\ \hline 3 \\ \cline{1-1} 4 \\ \cline{1-1} \end{array} \quad \begin{array}{|c|c|}\hline 2 & 3 \\ \hline 3 \\ \cline{1-1} 4 \\ \cline{1-1} \end{array} \quad \begin{array}{|c|c|}\hline 2 & 4 \\ \hline 3 \\ \cline{1-1} 4 \\ \cline{1-1} \end{array}$$

次のタブローは円柱状タブローの条件を満たさず $s^{(0)}_{\mathscr{D}_2(2,1,1)}$ には含まれない。

$$\begin{array}{|c|c|}\hline 1 & 4 \\ \hline 2 \\ \cline{1-1} 3 \\ \cline{1-1} \end{array}$$

従って以下のような多項式を得る。

$$\begin{aligned} s^{(0)}_{\mathscr{D}_2(2,1,1)} &= \mathrm{x}^{(4)}_1\mathrm{x}^{(3)}_1\mathrm{x}^{(1)}_2\mathrm{x}^{(2)}_3 + \mathrm{x}^{(4)}_1\mathrm{x}^{(3)}_1\mathrm{x}^{(1)}_2\mathrm{x}^{(2)}_4 + \mathrm{x}^{(4)}_1\mathrm{x}^{(3)}_1\mathrm{x}^{(1)}_3\mathrm{x}^{(2)}_4 + \mathrm{x}^{(4)}_1\mathrm{x}^{(3)}_2\mathrm{x}^{(1)}_2\mathrm{x}^{(2)}_3 \\ &+ \mathrm{x}^{(4)}_1\mathrm{x}^{(3)}_2\mathrm{x}^{(1)}_2\mathrm{x}^{(2)}_4 + \mathrm{x}^{(4)}_1\mathrm{x}^{(3)}_2\mathrm{x}^{(1)}_3\mathrm{x}^{(2)}_4 + \mathrm{x}^{(4)}_1\mathrm{x}^{(3)}_3\mathrm{x}^{(1)}_2\mathrm{x}^{(2)}_3 + \mathrm{x}^{(4)}_1\mathrm{x}^{(3)}_3\mathrm{x}^{(1)}_2\mathrm{x}^{(2)}_4 \\ &+ \mathrm{x}^{(4)}_1\mathrm{x}^{(3)}_3\mathrm{x}^{(1)}_3\mathrm{x}^{(2)}_4 + \mathrm{x}^{(4)}_1\mathrm{x}^{(3)}_4\mathrm{x}^{(1)}_2\mathrm{x}^{(2)}_4 + \mathrm{x}^{(4)}_1\mathrm{x}^{(3)}_4\mathrm{x}^{(1)}_3\mathrm{x}^{(2)}_4 + \mathrm{x}^{(4)}_2\mathrm{x}^{(3)}_2\mathrm{x}^{(1)}_3\mathrm{x}^{(2)}_4 \\ &+ \mathrm{x}^{(4)}_2\mathrm{x}^{(3)}_3\mathrm{x}^{(1)}_3\mathrm{x}^{(2)}_4 + \mathrm{x}^{(4)}_2\mathrm{x}^{(3)}_4\mathrm{x}^{(1)}_3\mathrm{x}^{(2)}_4. \end{aligned}$$

path $b=b_1\otimes b_2\otimes b_3\otimes b_4$ に対して以下のような座標を用いることにしよう。

$$\begin{aligned} b &= \boxed{1^{x^{(4)}_4}2^{x^{(1)}_4}3^{x^{(2)}_4}4^{x^{(3)}_4}} \otimes \boxed{1^{x^{(3)}_3}2^{x^{(4)}_3}3^{x^{(1)}_3}4^{x^{(2)}_3}} \otimes \boxed{1^{x^{(2)}_2}2^{x^{(3)}_2}3^{x^{(4)}_2}4^{x^{(1)}_2}} \otimes \boxed{1^{x^{(1)}_1}2^{x^{(2)}_1}3^{x^{(3)}_1}4^{x^{(4)}_1}} \\ &= \boxed{1^a2^b3^c4^d} \otimes \boxed{1^e2^f3^g4^h} \otimes \boxed{1^i2^j3^k4^l} \otimes \boxed{1^m2^n3^o4^p}. \end{aligned}$$

その時

$$
\begin{aligned}
\operatorname{trop} s^{(0)}_{\mathscr{D}_2(2,2,2,2)} &= p+o+l+k+h+g+d+c,\\
\operatorname{trop} s^{(0)}_{\mathscr{D}_2(2,1,1)} &= \min(p+o+l+h,\ p+o+l+c,\ p+o+g+c,\ p+j+l+h,\\
&\qquad p+j+l+c,\ p+j+g+c,\ p+e+l+h,\ p+e+l+c,\\
&\qquad p+e+g+c,\ p+d+l+c,\ p+d+g+c,\ k+j+g+c,\\
&\qquad k+e+g+c,\ k+d+g+c).
\end{aligned}
$$

すると式 (50) は $\nu^{(2)}_1 = \operatorname{trop} s^{(0)}_{\mathscr{D}_2(2,2,2,2)} - \operatorname{trop} s^{(0)}_{\mathscr{D}_2(2,1,1)}$ および $\nu^{(2)}_2 = \operatorname{trop} s^{(0)}_{\mathscr{D}_2(2,1,1)}$ となる。
数値的な例を考えると、以下の paths

$$\{a,b,c,d,e,f,g,h,i,j,k,l,m,n,o,p\} = \{3,2,c,3,3,3,1,0,0,3,0,2,1,0,3,3\}$$

に対して $\Phi(b)$ は以下のようになる。

| $c$ | $(\nu^{(1)}_1, J^{(1)}_1)$ | $(\nu^{(1)}_2, J^{(1)}_2)$ | $(\nu^{(1)}_3, J^{(1)}_3)$ | $(\nu^{(2)}_1, J^{(2)}_1)$ | $(\nu^{(2)}_2, J^{(2)}_2)$ | $(\nu^{(3)}_1, J^{(3)}_1)$ |
|---|---|---|---|---|---|---|
| 0 | $(8,-2)$ | $(8,-2)$ | $(4,-1)$ | $(8,4)$ | $(4,0)$ | $(8,-7)$ |
| 1 | $(9,-2)$ | $(8,-1)$ | $(4,-1)$ | $(8,2)$ | $(5,-1)$ | $(8,-6)$ |
| 2 | $(10,-2)$ | $(8,0)$ | $(4,-1)$ | $(8,0)$ | $(6,-2)$ | $(8,-5)$ |
| 3 | $(11,-2)$ | $(8,1)$ | $(4,-1)$ | $(8,-2)$ | $(7,-3)$ | $(8,-4)$ |
| 4 | $(12,-2)$ | $(8,2)$ | $(4,-1)$ | $(8,-4)$ | $(8,-4)$ | $(8,-3)$ |
| 5 | $(13,-2)$ | $(8,2)$ | $(4,-1)$ | $(9,-5)$ | $(8,-4)$ | $(8,-3)$ |
| 6 | $(14,-2)$ | $(8,2)$ | $(4,-1)$ | $(10,-6)$ | $(8,-4)$ | $(8,-3)$ |
| 7 | $(15,-2)$ | $(8,2)$ | $(4,-1)$ | $(11,-7)$ | $(8,-4)$ | $(8,-3)$ |

この結果は式 (50) による結果に一致する。一見して明らかなように riggings の振る舞いは $\nu^{(a)}$ の振る舞いよりはるかに複雑である。

## 8.7 $D^{(1)}_n$ 型艤装配位

$A^{(1)}_n$ 型艤装配位の理論を振り返ってみると、最も核心部を成す艤装配位写像とは組み合わせ $R$ 行列の親玉のような存在である、という点に注目すべきであろう。ところが $A^{(1)}_n$ 型の場合必要となる組み合わせ論的構造は Lascoux–Schützenberger 理論等に代表される代数的組み合わせ論の世界では既によく知られた構造であり、ある意味では結晶基底導入以前に既に存在していたと言っても良いほどのものであった。この点に関する一部の研究者の不満は根深いもののようであり、「あなたたちの研究は既に知られている事実の再定式化や単純な拡張に過ぎず、本質的には何も新しい発見をしていない」という意見が噴出するのもやむを得ないものと思っている。

幸いにして艤装配位の理論は従来の代数的組み合わせ論の手法では全く手が付けられなかったような問題でも十分に取り扱える場合がある。以下ではそのような事例について詳しく紹介し、今紹介したようなご意見に対する一つの答えとしようと思う。もちろんこれまで組み合わせ論に携わる多くの研究者にとって難問であった問題であれば、艤装配位の枠組みで研究してもその過程で大きな困難を伴うのは当然のことである。壁が高ければ高いほど新しい結果を得ている証拠であると思って先へ進むことにする。

**クリスタル $B^{r,s}$**　以下では $D_n^{(1)}$ 型の場合を考える。Dynkin 図は以下の通り。

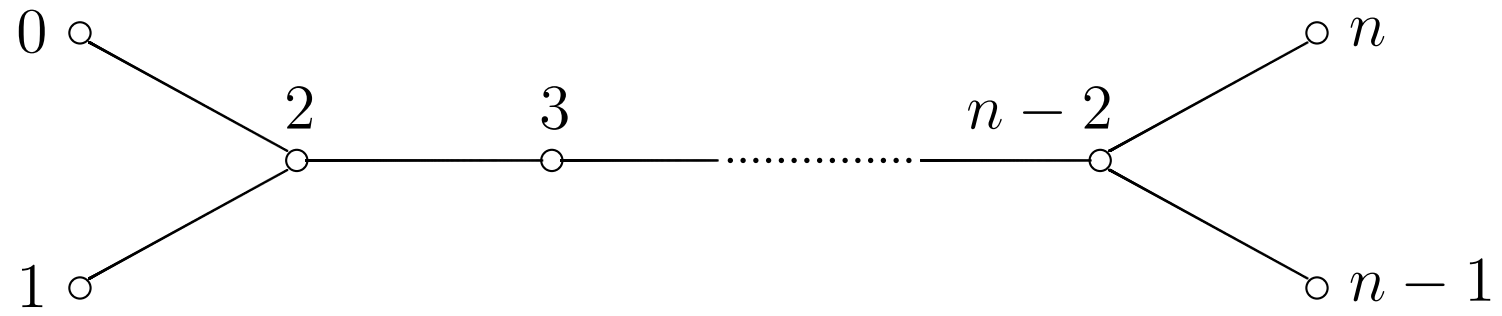

$A_n^{(1)}$ 型の場合であればクリスタル $B^{r,s}$ とは長方形型の Young 準標準盤の全体、と簡単に言ってしまえるのでイメージもわきやすいのであるが、$D_n^{(1)}$ 型の場合はすでにこの段階からかなり厄介な問題が生じる。

手始めは簡単な $B^{1,1}$ から始めよう。この場合 $B^{1,1}$ は一ますのタブロー$\boxed{1}$, $\boxed{2}$, . . . , $\boxed{n}$, $\boxed{\overline{n}}$, . . . , $\boxed{\overline{2}}$, $\boxed{\overline{1}}$からなる。柏原作用素 $\tilde{f}_a$ を矢印にラベル $a$ を添えて表示することにすると、$B^{1,1}$ 上の $D_n$ の作用（クリスタルグラフという）は以下の通り。

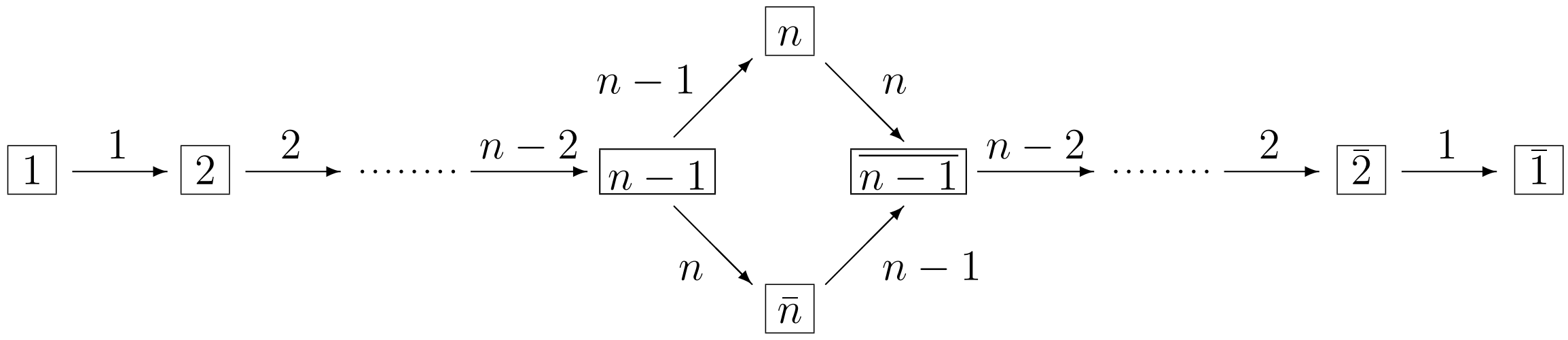

念のため付け加えておくと、$D_n^{(1)}$ 型の場合これらの作用に加えて

$$\tilde{f}_0 : \boxed{\overline{2}} \longmapsto \boxed{1}, \qquad \tilde{f}_0 : \boxed{\overline{1}} \longmapsto \boxed{2}$$

という作用が存在するので $B^{1,1}$ は最高ウェイト表現ではなくなる。

一般の $B^{r,s}$ では、部分代数 $D_n$ に制限した時、ウェイトとして縦 $r$ 横 $s$ の長方形から縦ドミノ ⌸ を任意個数取り除いたものを持つ表現が一つずつ発生する。Kashiwara–Nakashima タブロー [KN] とはそれら古典部分代数の表現のウェイトと同じ形の Young 盤に $B^{1,1}$ の時と同じ文字を埋めて $B^{r,s}$ の元を実現するという方法である。しかし無限次元の代数 $D_n^{(1)}$ を扱う場合にはこのやり方では深刻な困難に直面する。

そもそも $D_n^{(1)}$ の Kirillov–Reshetikhin クリスタルといった場合 $B^{r,s}$ という形で表されることから想像がつくように、基本的には長方形型の対象を取り扱う方が自然である。よって Kashiwara–Nakashima タブローを縦 $r$ 横 $s$ の長方形に置いた時に現れる「すき間」は一見空っぽに見えるが本当は情報が含まれているべき場所である。実際 $A$ 型の場合の真似をして組み合わせ $R$ 行列を $D$ 型版 row insertion によって計算してみようとすればすぐに分かるように、重要な情報が次々消去されて「すき間」に飲み込まれていってしまい計算不能になる。

よって $D_n^{(1)}$ のクリスタル $B^{r,s}$ の表示は長方形型の対象であるべきである。その様な問題意識を持っていた研究者も存在していたが何の方針もなく探してもうまく行かなかったようである。艤装配位写像を用いることによりそのような長方形型タブローを構成することができる [OSS11]。この様なタブローは KR タブローと呼ばれ、詳細な定義は §8.8 で与えることとするが、以下で必要となる重要な点のみ述べる。

KR タブローは縦 $r$ 横 $s$ の長方形タブローであり、古典部分代数 $D_n$ に関する最高ウェイト元についてのみ明示的に定義される。すなわち他のタブローは適切な $\tilde{f}_a$ の列を作用させて定義される。その場合一般のタブローの形はかなり複雑であり、Kashiwara–Nakashima タブローにみられるような特徴づけを与えることは今の所成功していない。なお $D$ 型以外の非例外型アフィン代数について行われた観察 [ScScr14] によれば、$D$ 型以外でもほぼ同じ定義で良いようであり、普遍性がある。同様な方向で定義すれば、艤装配位写像を構成すると懸案だった例外型の $E$ 型代数に対してもタブロー表示を定義できると思われる。

**$D_n^{(1)}$ 型艤装配位**　$A_n^{(1)}$ 型の場合と同様にして $D_n^{(1)}$ 型でも古典部分代数 $D_n$ 型 Dynkin 図の各頂点に Young 図 $\nu^{(a)}$ を配置することにより艤装配位を考える。以下 28 ページにおける定義と並行して定義する。rigging を $J_i^{(a)}$ と書くとき、次のような対象を考える。

$$(\nu, J) = \Big(\mu^{(1)}, \dots, \mu^{(n)}, (\nu^{(1)}, J^{(1)}), \dots, (\nu^{(n)}, J^{(n)})\Big). \tag{51}$$

ただし $\mu^{(a)}$ における添え字 $a$ のつけ方は $A_n^{(1)}$ 型の場合とはずらしてあることに注意（表現論の比重が重くなるとそのような記法の方が自然に見える）。$\nu^{(a)}$ を $D_n$ 型 Dynkin 図の各頂点と対応付ける必然性が生じるのは vacancy number の定義においてである。$A_n^{(1)}$ 型の場合の定義 (21) と同様にして、Dynkin 図上で単線でつながっている頂点には係数 1 を、自分自身には係数 $-2$ を与えることにより以下のように定義する（Cartan データを想起されたい）。

$$\begin{aligned}
P_l^{(a)}(\mu,\nu) &= Q_l(\mu^{(a)}) + Q_l(\nu^{(a-1)}) - 2Q_l(\nu^{(a)}) + Q_l(\nu^{(a+1)}) \qquad (1 \le a \le n-3),\\
P_l^{(n-2)}(\mu,\nu) &= Q_l(\mu^{(n-2)}) + Q_l(\nu^{(n-3)}) - 2Q_l(\nu^{(n-2)}) + Q_l(\nu^{(n-1)}) + Q_l(\nu^{(n)}),\\
P_l^{(a)}(\mu,\nu) &= Q_l(\mu^{(a)}) + Q_l(\nu^{(n-2)}) - 2Q_l(\nu^{(a)}) \qquad (a = n-1, n).
\end{aligned} \tag{52}$$

以降の定義は $A_n^{(1)}$ 型の場合と全く同様である。念のため書くと、$(\nu, J)$ が最高ウェイトベクトルに対する艤装配位であるとは以下の条件を満たす場合である。

- 全ての $\nu_i^{(a)}$ に対し $P_{\nu_i^{(a)}}^{(a)}(\mu,\nu) \ge 0$。
- 全ての rigging が $0 \le J_i^{(a)} \le P_{\nu_i^{(a)}}^{(a)}(\mu,\nu)$ を満たす。

一般の艤装配位は Kashiwara 作用素を艤装配位上に定義することによって得られる。ここで艤装配位上の $\widetilde{e}_a, \widetilde{f}_a$ は 29 ページにおける $A_n^{(1)}$ 型の場合と全く同様である [Sc06]。艤装配位の著しい自然さの一例である。

**$D_n^{(1)}$ 型艤装配位写像の構成**　詳しい定義は §8.8 に譲るが、$D_n^{(1)}$ 型でも完全に一般の $\bigotimes_{k=1}^{L} B^{r_k,s_k}$ 型 path に対して艤装配位写像

$$\Phi : \text{path} \longmapsto \text{rigged configuration}$$

が確立されている（[OSSS]、2016 年）。$A_n^{(1)}$ 型の場合の達成（[KSS]、1999 年）から 17 年もかかってしまったが、これは艤装配位の理論が新たな段階に到達したことを示す。最終的な解決は以下の三つの論文によってなされた。まず 2011 年に発表した論文 [OSS11] において艤装配位写像とそれを実現する KR タブローの定義および基本的な性質を予想した。同時に以下で述べる Dynkn 図の頂点 $0 \leftrightarrow 1$ の入れ替えに関し艤装配位が持っている重要な性質を与えた。次いで執筆された論文 [Sa14] において、艤装配位写像の持つ最も基本的な性質の一つ、柏原作用素との可換性

$$[\widetilde{e}_a, \Phi] = [\widetilde{f}_a, \Phi] = 0 \tag{53}$$

が確立された。証明は極めて困難であるが、$D_n^{(1)}$ 型艤装配位の最終的構成に当たってどうしても解決しなければならない最難関の壁として立ちはだかっていた問題であった。最終的に論文 [OSSS] において艤装配位写像の構成及びもう一つの基本定理である $R$ 不変性が確立された。$D_n^{(1)}$ 型の艤装配位写像の場合、アルゴリズムが複雑すぎて $R$ 不変性の証明はそのま

までは歯が立たない。そうした時上述の柏原作用素との可換性は強力な武器を提供する。これは path が最高ウェイト元の場合であっても同様の状況であり、本質的なステップであると考えている。

なお [OSSS] において艤装配位写像を構成する際に di Francesco–Kedem の結果 [dFK] を引用している個所があるが、筆者の印象ではおそらくそれは本質的なことではなく、除去可能なのではないかという気がしている。その場合 $D_n^{(1)}$ 型の写像を完成させた三つの論文の議論は原理的には全ての非例外型代数すべてについて適用可能な方法であろうと思われる。

最終的な解決だけでなく、開拓者的な仕事にも触れておく必要があるだろう。$D_n^{(1)}$ 型艤装配位写像の最初の形は 2003 年の論文 [OSS03] において $(B^{1,1})^{\otimes L}$ 型の path について与えられた。その後論文 [Sc04] において $\bigotimes_{k=1}^L B^{r_k,1}$ 型の path について、また論文 [ScSh] において $\bigotimes_{k=1}^L B^{1,s_k}$ 型の path について拡張された。特に論文 [Sc04, Appendix C] において予告された結果は最終的な構成でも重要な役割を果たした結果なのだが、残念ながら現在も未出版である。しかし $D_n^{(1)}$ 型といえども以上の結果に含まれるような単純な場合はそれまでの代数的組み合わせ論の手法の拡張である程度手がつけられたので（[HKOT] 参照)、艤装配位の理論が代数的組み合わせ論にとって真に新しい領域に到達するにはどうしても一般の $\bigotimes_{k=1}^L B^{r_k,s_k}$ 型 path に対して構成しなければならないと考えた次第である。

**様々な双対性**　艤装配位が単なるモデルではなく、無限次元対称性の核心部とかかわる数学的に意味のある量である、という事実の一つの証左として、対称性の持っている双対性などの非常に深い性質がごく自然な形で実現するという事を指摘しておこう。組み合わせ $R$ 行列が艤装配位上で自明化する $R$ 不変性は一つの典型的な例である。他の例として、Lusztig involution と呼ばれる深い双対性は艤装配位上以下のように現れる。与えられた艤装配位に対して、rigging $J_i^{(a)}$ と corigging $P_{\nu_i^{(a)}}^{(a)}(\mu,\nu) - J_i^{(a)}$ の役割を全て入れ替えて定義される艤装配位写像 $\widetilde{\Phi}$ を実行して得られる像は、もともとの path を Lusztig involution で写した像が所属する古典部分代数の表現の最高ウェイトベクトルと一致する [OSSS]。この様な深い性質を持った量がそういくつも存在するはずがないのは明らかであろう。なお、同様の双対性は 16 ページの式 (10) にも登場し、Bethe 仮説方程式の解の解析においても基本的であったことを注意しておく。

さらに次のような興味深い例がある [OSS11]。話の順序として Schilling によるアフィン柏原作用素 $\widetilde{e}_0, \widetilde{f}_0$ の構成 [Sc08] から必要な部分を見ておこう。アフィン代数としての $D_n^{(1)}$ の構造を KR クリスタル $B^{r,s}$ に導入するには避けて通れない部分であるが、単純にはいかない。そこで Dynkin 図の 0 番目と 1 番目の頂点を入れ替えるような involution $\sigma$ を構成し、良く知られた（そして単純な）$\widetilde{e}_1$ と $\widetilde{f}_1$ の作用に帰着させることとなる。よって核心部は involution $\sigma$ の構成となるが、これも一筋縄には行かない。$A$ 型の楽園から外に出てしまうと何かと難しくなってしまう。

以下簡単のため $r \leq n-2$ とする。Dynkin 図の頂点 $I = \{0,1,\ldots,n\}$ の部分集合 $J \subset I$ に対し、$\widetilde{e}_a b = 0$ が全ての $a \in J$ に対し成り立つ時 $b$ を $J$-最高ウェイト元と呼ぶことにしよう。上で述べた involution $\sigma$ を一般的に構成するのは至難の業だが、頂点 0 と 1 を除いた $J = \{2,3,\ldots,n\}$-最高ウェイト元に対してなら構成することができる。構成の中心は $J$-最高ウェイト元 $b$ から $\pm$-ダイヤグラムと呼ばれるものへの非自明な全単射である。$\pm$-ダイヤグラムとは Young 盤の組 $\lambda \supset \mu$ に対し $\lambda/\mu$ の部分に $+,-$ の符号を入れたもので、各列で $+,-$ の符号はそれぞれ最大一つ（両方無しも可)、かつ両符号が共に存在するときは上側に

+ を入れるものとしよう。例えば、

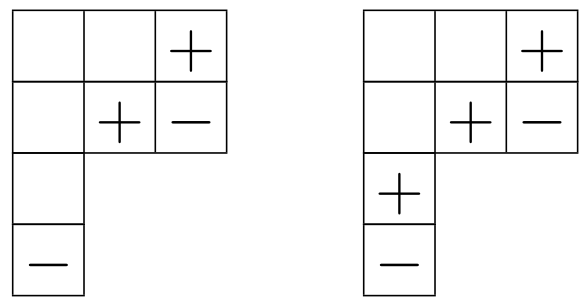

といったものである。±-ダイヤグラムの形 $\lambda$ は古典部分代数 $D_n$ の最高ウェイトの形に合わせるので（つまり Kashiwara–Nakashima タブローと同じ形）、$B^{r,s}$ ならば縦 $r$ 横 $s$ の長方形から縦ドミノ $\genfrac{}{}{0pt}{}{\square}{\square}$ を任意個数取り除いたものとなる。求める involution $\sigma$ は ±-ダイヤグラムの上で以下のような形をとる。すなわち内側の形 $\mu$ は固定したうえで、+ および −、± および空集合をそれぞれ入れ替える。図式的には、

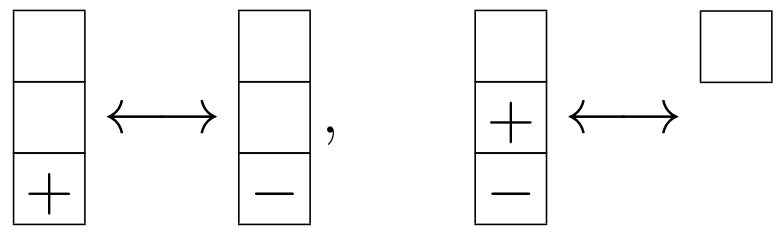

という形である。

まとめると、$J$-最高ウェイト元 $b$ から ±-ダイヤグラムへの非自明な全単射を用いると involution $\sigma$ が非常に自然な実現をすることが分かった。単に「表現を構成する」という事だけが目的であればこれで十分であるかもしれない。しかし表現の構成というのは、いかに重要な課題であったとしても、「対称性そのものを理解すること」という最終的な目標に対するそれ以外の接近法の重要性を減ずるものではないだろう。殊に無限次元の対称性は有限次元の対称性とは圧倒的に異なる巨大さを持つので、従来の理論的枠組みに加えて何か本質的に新しい考え方が必要とされるであろうという事を常に肝に銘じておかなければならない。例えば我々の問題としている $D_n^{(1)}$ 代数が表している対称性の理解という観点からは、±-ダイヤグラムおよびそこへの全単射というものをブラックボックスとせずに、それらが数学的に何を表しているのかという点まで踏み込んで理解してみることも有用であろう。論文 [OSS11] で筆者が念頭としていたのはその様な問題意識であり、実際艤装配位の理論という大きな枠組みを通してみると ±-ダイヤグラムの正しい居場所が明らかになる。

より具体的にいえば、以下のような対応により、±-ダイヤグラムとは本質的には艤装配位の事であり、$J$-最高ウェイト元からの全単射とは本質的には艤装配位写像の事であった。クリスタル $B^{r,s}$ の ±-ダイヤグラム $P$ のある列 $c$ について考える。$c$ の高さが $x$、また $x+y=r$ とすると、$y$ はタブロー $P$ を $r\times s$ の長方形に置いた時のすき間の高さとなる。この「すき間」は縦ドミノ $\genfrac{}{}{0pt}{}{\square}{\square}$ を取り除いて得られるので、$y$ は偶数となる。さてまず各列 $c$ に対応する艤装配位を与える。以下分割 $\nu^{(a)}$ の形を一行目に、対応する rigging を二行目に書く。

(A) $c$ が符号を含まないとき。

$$\nu=(\overbrace{(\ \ 1),(1),\cdots,(1)}^{x},\overbrace{(1),(1,1),\cdots,(1,\ldots,1)}^{y},(1^y),\cdots,(1^y),(1^{\frac{y}{2}}),(1^{\frac{y}{2}}))$$
$$J=(\overbrace{(-1),(0),\cdots,(0)}^{x},\overbrace{(1),(0,0),\cdots,(0,\ldots,0)}^{y},(0^y),\cdots,(0^y),(0^{\frac{y}{2}}),(0^{\frac{y}{2}}))$$

(B) $c$ が符号 + を含むとき。

$$\nu=(\overbrace{\emptyset,\emptyset,\cdots,\emptyset}^{x},\overbrace{(1),(1,1),\cdots,(1,\ldots,1)}^{y},(1^y),\cdots,(1^y),(1^{\frac{y}{2}}),(1^{\frac{y}{2}}))$$
$$J=(\overbrace{\emptyset,\emptyset,\cdots,\emptyset}^{x},\overbrace{(0),(0,0),\cdots,(0,\ldots,0)}^{y},(0^y),\cdots,(0^y),(0^{\frac{y}{2}}),(0^{\frac{y}{2}}))$$

(C) $c$ が符号 $-$ を含むとき。

$$\nu = (\overbrace{(\ 2),(2),\cdots,(2)}^{x-1},\overbrace{(1,1),(1,1,1),\cdots,(1,\ldots,1)}^{y+1},(1^{y+2}),\cdots,(1^{y+2}),(1^{\frac{y+2}{2}}),(1^{\frac{y+2}{2}}))$$

$$J = (\overbrace{(-2),(0),\cdots,(0)}^{x-1},\overbrace{(0,0),(0,0,0),\cdots,(0,\ldots,0)}^{y+1},(0^{y+2}),\cdots,(0^{y+2}),(0^{\frac{y+2}{2}}),(0^{\frac{y+2}{2}}))$$

(D) $c$ が符号 $\pm$ を含むとき。

$$\nu = (\overbrace{(\ 1),(1),\cdots,(1)}^{x-1},\overbrace{(1,1),(1,1,1),\cdots,(1,\ldots,1)}^{y+1},(1^{y+2}),\cdots,(1^{y+2}),(1^{\frac{y+2}{2}}),(1^{\frac{y+2}{2}}))$$

$$J = (\overbrace{(-1),(0),\cdots,(0)}^{x-1},\overbrace{(0,0),(0,0,0),\cdots,(0,\ldots,0)}^{y+1},(0^{y+2}),\cdots,(0^{y+2}),(0^{\frac{y+2}{2}}),(0^{\frac{y+2}{2}}))$$

若干補足をすると、$\pm$ の場合の $\nu^{(x-1)}$ を除きすべて特異なストリングである。また $c$ が空の時（つまり $x=0$ の時）は $+$ の場合の特別な場合として扱う。また $-$ の場合で $x=1$ ならば $(\nu^{(1)},J^{(1)})=((1,1),(-1,-1))$ と取る。

一般の $\pm$-ダイヤグラム $P$ に対応する艤装配位は、各列 $c$ に対応する艤装配位を「加えて」得られる。ここで「加える」とは各分割 $\nu^{(a)}$ では Young 図を横に全て結合し、また rigging は対応するものを全て足しあげることによって得られる。こうして得られる艤装配位が確かに艤装配位写像 $\Phi$ の下で $\pm$-ダイヤグラム $P$ に対応する艤装配位であることは、論文 [OSS11] においてクリスタル同型であることが確かめられたことと論文 [Sa14] において艤装配位写像 $\Phi$ と柏原作用素 $\widetilde{e}_a,\widetilde{f}_a$ $(a=1,\ldots,n)$ とが可換であることが確立されたことを合わせて証明された（[Sa14, §4.1] 参照）。$D_n^{(1)}$ 型の艤装配位写像は複雑であるから、こういった結果を証明するのにも柏原作用素との可換性のような飛び道具がないと極めて困難である。

2010 年から 2011 年にかけて $D_n^{(1)}$ 型艤装配位写像の形が大体予想出来てきたころ、折角今まで誰も足を踏み入れたことのない領域に到達したのだからちょっとは面白い現象があるだろう、と思いあれこれの論文を出来立ての写像を使って調べているうちに以上のような結果に行きつき、大変面白く感じたのを思い出す。当時写像を計算していると 8.8 節で述べる Kirillov–Reshetikhin タブローという実に奇妙なもの（特に $\pm$-ダイヤグラムのように最高ウェイト元でない場合）が出てきてなかなか自信が持てなかったのだが、この現象や、組み合わせ $R$ 行列の計算をしているうちに少しずつ正しさを確信して行ったように思う。最終的には論文 [OSS11] を公表した後 Scrimshaw 氏が別のプログラムを組み私のもろもろの予想を独立に確認して頂き、ほぼ間違いがないことが分かってほっとしたことを思い出す。

最後に一つ例を見ておこう。$D_n^{(1)}$ 型 $(n\geq 10)$ のクリスタル $B^{8,5}$ の以下の元を考える：

$$P = \begin{array}{|c|c|c|c|}\hline \phantom{+} & \phantom{+} & \phantom{+} & + \\ \hline & & & - \\ \hline & & \\ \cline{1-3} & & \\ \cline{1-3} & + \\ \cline{1-2} - & - \\ \cline{1-2}\end{array}$$

最後の列は空である。一行目に $\nu_P$ を、二行目に $J_P$ を書くことにすると、対応する艤装配位は以下の通り：

$$((\ 6),(6,2),(6,2,2),(6,2,2,2),(6,2,2,2,2),(5,5,2,2,2,2),(5,5,5,2,2,2,2),\cdots)$$
$$((-5),(0,0),(0,0,0),(0,0,0,0),(1,0,0,0,0),(0,0,0,0,0,0),(0,0,0,0,0,0,0),\cdots)$$

ここでは例として $\nu_P^{(5)}$ を見てみよう。これは $P$ の各列に由来する下記の艤装配位の「和」である：

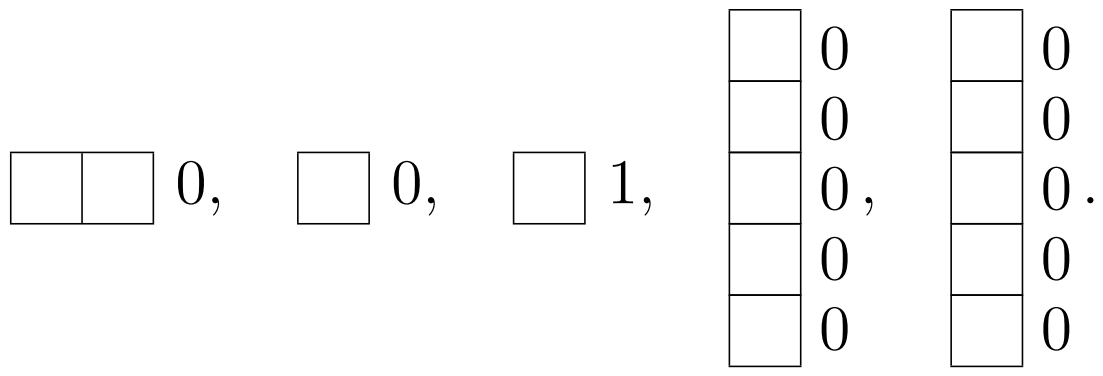

ここでは対応する riggings のみ表示した。これらを加えて $\nu_P^{(5)}$

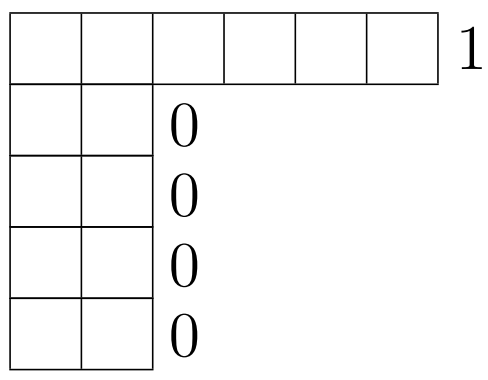

が得られる。

**箱玉系の逆散乱形式**　この話は $D_n^{(1)}$ 型艤装配位写像の研究の主要な動機の一つであり、また主要な結果の一つなのであるが、論文 [OSSS] では残念ながらある共著者の強い反対により触れることができなかったので[27]、ここで述べることにする。

結果自体は32ページで述べた $A_n^{(1)}$ 型の場合と全く同様であるが、念のため結果を書いておこう。時間発展 $T^{a,l}$ は $A_n^{(1)}$ 型の場合と全く同じ元 $u^{a,l}$ を用いて定義される。$D_n^{(1)}$ 型の一般的な path $b \in \bigotimes_{k=1}^{L} B^{r_k,s_k}$ に必要ならば右側に $u^{a,1}$ を十分付け加えて、$u^{a,l} \otimes b \overset{R}{\simeq} T^{a,l}(b) \otimes u^{a,l}$ となるようにしよう。その時 $b$ に対する艤装配位が

$$\Phi(b) = \Big((\nu^{(1)}, J^{(1)}), \dots, (\nu^{(a)}, J^{(a)}), \dots, (\nu^{(n)}, J^{(n)})\Big)$$

であれば、時間発展した path に対する艤装配位は

$$\Phi\left(T^{a,l}(b)\right) = \Big((\nu^{(1)}, J^{(1)}), \dots, \big(\nu_i^{(a)}, J_i^{(a)} + \min(\nu_i^{(a)}, l)\big)_i, \dots, (\nu^{(n)}, J^{(n)})\Big) \tag{54}$$

となる。すなわち Young 図 $\nu^{(1)}, \dots, \nu^{(n)}$ は運動の保存量（作用変数）であり、rigging $J^{(a)}$ のみが線形に変化する（角変数）。スピン表現かどうかにかかわりなく $A_n^{(1)}$ 型の場合と同じ結果になるわけであり、艤装配位の自然さを表している。

証明も $A_n^{(1)}$ 型の場合と全く同様である。32ページと全く同様に以下の定理が成り立つ。

**定理 [OSSS]（$D_n^{(1)}$ 型艤装配位の $R$ 不変性）**

任意の $D_n^{(1)}$ 型 KR クリスタルのテンソル積 $b$ と $b'$ が $b \overset{R}{\simeq} b' \iff \Phi(b) = \Phi(b')$.

証明はまずスピン表現と他の元のテンソル積 $B^{n,1} \otimes B^{r,s}$ の場合に不変性を証明する。といっても $D$ 型艤装配位写像は複雑であるから直接証明するのはおそらく絶望的であり、艤装配位写像と柏原作用素の可換性 [Sa14] を用いて path と艤装配位双方に適当な柏原作用素の列を作用させていって証明することになる。この辺りは既に論文 [OSS11] を執筆している最中から共著者間の議論で認識しており、最終的な困難は柏原作用素との可換性にあり、と観念して論文 [Sa14] を準備した。さて、この部分が解決すると、あとは箱玉系の時間発展 $T^{n,1}$ を十分繰り返すことにより $D_n^{(1)}$ 型の場合の $R$ 不変性が $A_n^{(1)}$ 型の場合の $R$ 不変性に帰着され

[27] 結局 $A$ 型の場合も $D$ 型の場合もこの重要な結果はおかしな出版の形態となってしまい誠に遺憾である。

てしまい、[KSS] の定理より結果が従う。$A_n^{(1)}$ 型の場合も [KSS] 論文の長大な議論をほぼ完全に使って証明しているので、論文 [Sa14] の同じく長大な証明もあわせると全体の証明は膨大なものである。真に深い結果であると言えよう。

この後の式 (54) の証明は $A_n^{(1)}$ 型の場合と全く同様である。

$R$ 不変性の定理について補足しておく。$A_n^{(1)}$ 型の場合と同様この結果は $D_n^{(1)}$ 型の組み合わせ $R$ 行列を特殊な場合として含んでいる。$D_n^{(1)}$ 型の場合従来の代数的組み合わせ論の手法では全く歯が立たない難問であることが知られていたので、艤装配位の理論が代数的組み合わせ論の世界に大きな手法をもたらしたことがようやく証明されたと言えよう。

**Misra–Wilson の箱玉系**　せっかく非常に一般的な逆散乱形式が確立されたのだから非自明な例を考察してみよう。幸い Misra–Wilson [MW] が $D_n^{(1)}$ 型代数の $(B^{2,1})^{\otimes L}$ 型テンソル積に対する箱玉系の $T^{2,s}$ 時間発展という極めて非自明な系を考察しているので、彼らの結果を艤装配位と比較してみることにしよう。彼らの結果をかいつまんで紹介すると、$B^{2,s} \otimes B^{2,1}$ 型の組み合わせ $R$ 行列を書き下すことにより可能なソリトンが分類された [MW, Proposition 3]。更に二つのソリトンの散乱を考察し、それらが $D_n$ 型 Dynkin 図の頂点 2 を取り除いて得られる代数に対する組み合わせ $R$ 行列によって記述されることが示された。

以下では彼らの論文に出ている例を艤装配位を用いて検討してみよう。逆散乱形式 (54) から分かるように、与えられた path $b \in (B^{2,1})^{\otimes L}$ に対し $s \geq \max_i \nu_i^{(2)}$ とすれば $T^{2,s}(b) = T^{2,s+1}(b) = \cdots$ が成り立つ。以下この様な時間発展を $T_\infty$ と書くことにすると、論文 [MW, §IV-B] で解析されている以下の三つの例では確かに各ソリトンがストリング $(\nu_i^{(2)}, J_i^{(2)})$ に対応していることが確認できる。

**例 1**　まず $D_4^{(1)}$ 型の例を考察しよう。以下の図では第一行目に与えられた path $b \in (B^{2,1})^{\otimes 23}$ を表示し、以下第 $t$ 行目には $T_\infty^{t-1}(b)$ を表示する。見やすくするために全てのテンソル積の記号 $\otimes$ やタブローの枠を省略し、また今の場合真空ないし背景と解釈される元 $u^{2,1} = \begin{array}{|c|}\hline 1 \\ \hline 2 \\ \hline \end{array} \in B^{2,1}$ を $\begin{array}{|c|}\hline \cdot \\ \hline \cdot \\ \hline \end{array}$ と表すことにする。

| | | | | | | | | | | | | | | | | | | | | | | |
|---|---|---|---|---|---|---|---|---|---|---|---|---|---|---|---|---|---|---|---|---|---|---|
| 1 | 1 | 1 | · | · | 2 | 2 | · | · | · | · | · | · | · | · | · | · | · | · | · | · | · | · |
| $\bar{3}$ | $\bar{4}$ | $\bar{4}$ | · | · | 4 | 3 | · | · | · | · | · | · | · | · | · | · | · | · | · | · | · | · |
| · | · | · | 1 | 1 | 1 | · | 2 | 2 | · | · | · | · | · | · | · | · | · | · | · | · | · | · |
| · | · | · | $\bar{3}$ | $\bar{4}$ | $\bar{4}$ | · | 4 | 3 | · | · | · | · | · | · | · | · | · | · | · | · | · | · |
| · | · | · | · | · | · | 1 | 1 | 1 | 2 | 2 | · | · | · | · | · | · | · | · | · | · | · | · |
| · | · | · | · | · | · | $\bar{3}$ | $\bar{4}$ | $\bar{4}$ | 4 | 3 | · | · | · | · | · | · | · | · | · | · | · | · |
| · | · | · | · | · | · | · | · | · | 1 | 1 | $\bar{4}$ | 2 | · | · | · | · | · | · | · | · | · | · |
| · | · | · | · | · | · | · | · | · | $\bar{3}$ | $\bar{4}$ | 4 | 3 | · | · | · | · | · | · | · | · | · | · |
| · | · | · | · | · | · | · | · | · | · | · | · | 1 | $\bar{4}$ | 2 | 1 | · | · | · | · | · | · | · |
| · | · | · | · | · | · | · | · | · | · | · | · | $\bar{3}$ | $\bar{3}$ | 3 | 3 | · | · | · | · | · | · | · |
| · | · | · | · | · | · | · | · | · | · | · | · | · | · | 1 | 1 | 2 | 2 | 1 | · | · | · | · |
| · | · | · | · | · | · | · | · | · | · | · | · | · | · | $\bar{3}$ | $\bar{4}$ | $\bar{3}$ | 3 | 3 | · | · | · | · |
| · | · | · | · | · | · | · | · | · | · | · | · | · | · | · | · | 1 | 1 | · | 2 | 2 | 1 | · |
| · | · | · | · | · | · | · | · | · | · | · | · | · | · | · | · | $\bar{3}$ | $\bar{4}$ | · | $\bar{3}$ | 3 | 3 | · |

状態 $T_\infty^t(b)$ に対応する艤装配位 $\Phi(T_\infty^t(b))$ は以下の通り。

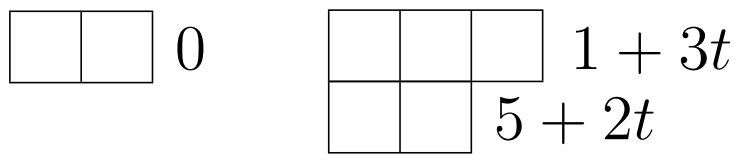

$\nu^{(2)}$ には長さ 3 と 2 の行があり、それぞれの rigging は速さ 3 と 2 で時間発展する。これは path のレベルで視覚的に了解される事実と一致する。

**例 2**　次に同論文で解析されている $D_5^{(1)}$ 型の例について、先ほどの例と同様に表示すると以下のようになる。

| | | | | | | | | | | | | | | | | | | | | | | |
|---|---|---|---|---|---|---|---|---|---|---|---|---|---|---|---|---|---|---|---|---|---|---|
| 2 | 2 | 1 | 1 | · | · | 2 | 1 | · | · | · | · | · | · | · | · | · | · | · | · | · | · | · |
| $\bar{3}$ | 5 | 4 | 3 | · | · | $\bar{4}$ | $\bar{5}$ | · | · | · | · | · | · | · | · | · | · | · | · | · | · | · |
| · | · | · | · | 2 | 2 | 1 | 1 | 2 | 1 | · | · | · | · | · | · | · | · | · | · | · | · | · |
| · | · | · | · | $\bar{3}$ | 5 | 4 | 3 | $\bar{4}$ | $\bar{5}$ | · | · | · | · | · | · | · | · | · | · | · | · | · |
| · | · | · | · | · | · | · | · | 2 | 2 | 4 | 1 | 1 | · | · | · | · | · | · | · | · | · | · |
| · | · | · | · | · | · | · | · | $\bar{3}$ | 5 | $\bar{4}$ | $\bar{5}$ | 3 | · | · | · | · | · | · | · | · | · | · |
| · | · | · | · | · | · | · | · | · | · | · | 2 | 1 | 2 | 2 | 1 | 1 | · | · | · | · | · | · |
| · | · | · | · | · | · | · | · | · | · | · | 5 | 4 | $\bar{3}$ | $\bar{4}$ | $\bar{5}$ | 3 | · | · | · | · | · | · |
| · | · | · | · | · | · | · | · | · | · | · | · | · | 2 | 1 | · | · | 2 | 2 | 1 | 1 | · | · |
| · | · | · | · | · | · | · | · | · | · | · | · | · | 5 | 4 | · | · | $\bar{3}$ | $\bar{4}$ | $\bar{5}$ | 3 | · | · |

状態 $T_\infty^t(b)$ に対応する艤装配位 $\Phi(T_\infty^t(b))$ は以下の通り。

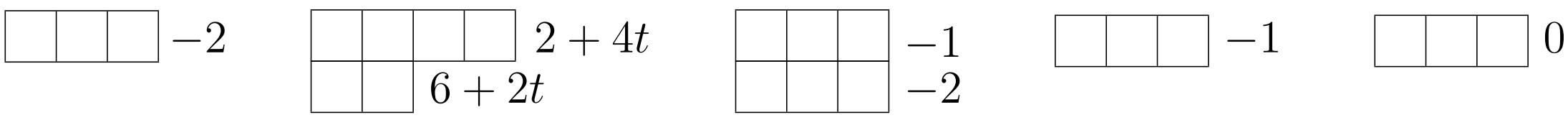

**例 3**　最後に同論文の $D_6^{(1)}$ 型の例を見てみよう。

| | | | | | | | | | | | | | | | | | | | | | | | | |
|---|---|---|---|---|---|---|---|---|---|---|---|---|---|---|---|---|---|---|---|---|---|---|---|---|
| 2 | 2 | 1 | 1 | 1 | · | · | 2 | 2 | · | · | · | · | · | · | · | · | · | · | · | · | · | · | · | · |
| $\bar{3}$ | $\bar{5}$ | 6 | 5 | 4 | · | · | $\bar{5}$ | $\bar{5}$ | · | · | · | · | · | · | · | · | · | · | · | · | · | · | · | · |
| · | · | · | · | · | 2 | 2 | 1 | 1 | 4 | 2 | · | · | · | · | · | · | · | · | · | · | · | · | · | · |
| · | · | · | · | · | $\bar{3}$ | $\bar{5}$ | 6 | 5 | $\bar{5}$ | $\bar{5}$ | · | · | · | · | · | · | · | · | · | · | · | · | · | · |
| · | · | · | · | · | · | · | · | · | · | 2 | $\bar{5}$ | 6 | 2 | 1 | · | · | · | · | · | · | · | · | · | · |
| · | · | · | · | · | · | · | · | · | · | $\bar{3}$ | $\bar{4}$ | $\bar{5}$ | 4 | 4 | · | · | · | · | · | · | · | · | · | · |
| · | · | · | · | · | · | · | · | · | · | · | · | · | 1 | 1 | 2 | 2 | 2 | 2 | 1 | · | · | · | · | · |
| · | · | · | · | · | · | · | · | · | · | · | · | · | $\bar{5}$ | 6 | $\bar{3}$ | $\bar{4}$ | $\bar{5}$ | 4 | 4 | · | · | · | · | · |
| · | · | · | · | · | · | · | · | · | · | · | · | · | · | · | 1 | 1 | · | · | · | 2 | 2 | 2 | 2 | 1 |
| · | · | · | · | · | · | · | · | · | · | · | · | · | · | · | $\bar{5}$ | 6 | · | · | · | $\bar{3}$ | $\bar{4}$ | $\bar{5}$ | 4 | 4 |

状態 $T_\infty^t(b)$ に対応する艤装配位 $\Phi(T_\infty^t(b))$ は以下の通り。

$$
\begin{aligned}
&(\nu^{(1)}, J^{(1)}) = \{(4), (-2)\}, && (\nu^{(2)}, J^{(2)}) = \{(5,2), (3+5t, 7+2t)\}, \\
&(\nu^{(3)}, J^{(3)}) = \{(5,3), (-2,-1)\}, && (\nu^{(4)}, J^{(4)}) = \{(4,3), (0,0)\}, \\
&(\nu^{(5)}, J^{(5)}) = \{(5), (-3)\}, && (\nu^{(6)}, J^{(6)}) = \{(4), (-1)\}.
\end{aligned}
$$

ここで艤装配位を以下のように記述した。

$$(\nu^{(a)}, J^{(a)}) = \{(\nu_1^{(a)}, \nu_2^{(a)}, \ldots), (J_1^{(a)}, J_2^{(a)}, \ldots)\}.$$

以上 Misra–Wilson による極めて非自明な $D$ 型箱玉系は艤装配位の理論からも明瞭に理解できることが確認された。

## 8.8 $D_n^{(1)}$ 型艤装配位写像のアルゴリズム

$D_n^{(1)}$ 型の path $b \in \bigotimes_{k=1}^{L} B^{r_k,s_k}$ と艤装配位

$$(\nu, J) = \left(\mu^{(1)}, \ldots, \mu^{(n)}, (\nu^{(1)}, J^{(1)}), \ldots, (\nu^{(n)}, J^{(n)})\right)$$

との間の全単射

$$\Phi^{-1} : (\nu, J) \longmapsto \text{path}$$

の記述を行う。ただしここでも艤装配位の $\mu^{(a)}$ のラベル $a$ は $A_n^{(1)}$ 型の場合とずらしてあることに注意。

主に定義するべき内容は $B^{r,s}$ を記述するための Kirillov–Reshetikhin タブローと写像 $\Phi^{-1}$ のアルゴリズムの二種類に大別される。なお以下簡単のため $r \le n-2$ と仮定する。スピン表現の場合については論文 [Sa14, OSSS] を参照されたい。

**Kirillov–Reshetikhin タブロー** KR タブローは最初論文 [OSS11] において導入された。その後論文 [ScScr14] において $A$ 型と $D$ 型以外の非例外型代数についても数値的に検討され、それぞれの代数に対する柏原作用素を用いれば細かな修正を除き実質的に同じ定義が使用できることが確かめられた。すなわち普遍性のある自然な概念であると言えよう。

$I_0 = \{1, 2, \ldots, n\}$ とする。KR タブローは $B^{r,s}$ の古典部分代数 $D_n$ に対する最高ウェイト元 $u_\lambda$、すなわち $\tilde{e}_a u_\lambda = 0$ $(a \in I_0)$ に対してのみ明示的に定義される。その他の元は柏原作用素 $\tilde{f}_a$ $(a \in I_0)$ を通常通り作用させることにより得られる。元 $u_\lambda$ のウェイト $\lambda$ の形を、高さ $h$ のコラムが $k_h$ 個あるというように表記しよう。基本ウェイトの記号を使えば、$\lambda = k_r\bar{\Lambda}_r + k_{r-2}\bar{\Lambda}_{r-2} + \cdots$ となる。

$u_\lambda$ に対応する KR タブローを $\mathrm{fill}(u_\lambda)$ と書き、その対応を **filling map** と呼ぼう。以下 filling map を定義する。$k_c$ を数列 $k_{r-2}, k_{r-4}, \ldots$ における最初の奇数としよう。その様な $k_c$ が存在しないときは $k_c = k_{-1}$、すなわち $c = -1$ とおく。$t$ を $u_\lambda$ の Kashiwara–Nakashima タブロー表示とする。すなわち Young 図 $\lambda$ の一行目には文字 1 を、二行目には文字 2 を、と順に埋めたものとする。タブロー $t$ を $r \times s$ の長方形の左上におき、すき間を以下のように埋めていく。手続きは $t$ の左端の列から始め、以下の手順で右向きに進む。

1. 高さ $r$ の列には何もしない。$c \ge 0$ である場合、高さ $c$ の列を一つ減らし、また高さ $c$ 未満の列を全て左に一列ずらす。

2. 高さが $c$ 以上の列の充填法。高さ $h$ の列に対するすき間であれば以下の文字列の転置で埋める。

   | | | | |
   |---|---|---|---|
   | $\overline{r}$ | $\overline{r-1}$ | $\cdots$ | $\overline{h+1}$ |
   | $h+1$ | $h+2$ | $\cdots$ | $r$ |

   .

   $c = -1$ であればこれで全てのすき間が埋まるので手続きが終了する。

3. 右端を除く他の列の充填法。文字 $x$ を再帰的に再定義しながら下記の内容で埋める。初期条件として $x = c+1$ と定める。対応する列の高さが $h$ であれば、

   | | | | | | | | |
   |---|---|---|---|---|---|---|---|
   | $y$ | $\cdots$ | $r-1$ | $r$ | $\overline{r}$ | $\cdots$ | $\overline{x+1}$ | $\overline{x}$ |

   の転置をすき間に埋める。すき間のます目の数は $r-h$ なので $y = r-(x-h-2)$ となる。$x = y$ と再定義して次の列に同じ手順を繰り返す。

4. 最終的に得られた $x$ を用いて

| | | | | | | | |
|---|---|---|---|---|---|---|---|
| 1 | 2 | $\cdots$ | $y$ | $\overline{y}$ | $\cdots$ | $\overline{x+1}$ | $\overline{x}$ |

の転置を右端の列とする。ここで $c \geq 0$ の場合 Step 1 の結果右端の列は空白になっている事に注意。また以上の定義から $y = (r + x - 1)/2$ と定まる。

**例 1**　$B^{8,7}$ の $I_0$-最高ウェイト元 $u_\lambda$ で $\lambda = \bar{\Lambda}_8 + 2\bar{\Lambda}_6 + \bar{\Lambda}_4 + 2\bar{\Lambda}_2$ なるものを考える。この場合 $(k_6, k_4, k_2, k_0) = (2, 1, 2, 1)$ なので $c = 4$ となる。$u_\lambda$ の KN タブロー表示を $8 \times 7$ の長方形の左上に置くと、

| | | | | | | |
|---|---|---|---|---|---|---|
| 1 | 1 | 1 | 1 | 1 | 1 | |
| 2 | 2 | 2 | 2 | 2 | 2 | |
| 3 | 3 | 3 | 3 | | | |
| 4 | 4 | 4 | 4 | | | |
| 5 | 5 | 5 | | | | |
| 6 | 6 | 6 | | | | |
| 7 | | | | | | |
| 8 | | | | | | |

$\xrightarrow{\text{Step 1}}$

| | | | | | | |
|---|---|---|---|---|---|---|
| 1 | 1 | 1 | 1 | 1 | | |
| 2 | 2 | 2 | 2 | 2 | | |
| 3 | 3 | 3 | | | | |
| 4 | 4 | 4 | | | | |
| 5 | 5 | 5 | | | | |
| 6 | 6 | 6 | | | | |
| 7 | | | | | | |
| 8 | | | | | | |

.

そこですき間を次のように埋める。

| | | | | | | |
|---|---|---|---|---|---|---|
| | | | | | 5 | 1 |
| | | | | | 6 | 2 |
| | | | 7 | 5 | 7 | 3 |
| | | | 8 | 6 | 8 | 4 |
| | | | $\bar{8}$ | 7 | $\bar{8}$ | 5 |
| | | | $\bar{7}$ | 8 | $\bar{7}$ | 6 |
| | $\bar{8}$ | 7 | $\bar{6}$ | $\bar{8}$ | $\bar{6}$ | $\bar{6}$ |
| | $\bar{7}$ | 8 | $\bar{5}$ | $\bar{7}$ | $\bar{5}$ | $\bar{5}$ |

, よって, $\mathrm{fill}(u_\lambda) =$

| | | | | | | |
|---|---|---|---|---|---|---|
| 1 | 1 | 1 | 1 | 1 | 5 | 1 |
| 2 | 2 | 2 | 2 | 2 | 6 | 2 |
| 3 | 3 | 3 | 7 | 5 | 7 | 3 |
| 4 | 4 | 4 | 8 | 6 | 8 | 4 |
| 5 | 5 | 5 | $\bar{8}$ | 7 | $\bar{8}$ | 5 |
| 6 | 6 | 6 | $\bar{7}$ | 8 | $\bar{7}$ | 6 |
| 7 | $\bar{8}$ | 7 | $\bar{6}$ | $\bar{8}$ | $\bar{6}$ | $\bar{6}$ |
| 8 | $\bar{7}$ | 8 | $\bar{5}$ | $\bar{7}$ | $\bar{5}$ | $\bar{5}$ |

.

**例 2**　見慣れない概念でしょうからもう少し $r = 12$ の例を。

| | | | | | | | | |
|---|---|---|---|---|---|---|---|---|
| 1 | 1 | 1 | 1 | 1 | 1 | 1 | 7 | 1 |
| 2 | 2 | 2 | 2 | 2 | 2 | 2 | 8 | 2 |
| 3 | 3 | 3 | 3 | 7 | 9 | 7 | 9 | 3 |
| 4 | 4 | 4 | 4 | 8 | 10 | 8 | 10 | 4 |
| 5 | 5 | 5 | 5 | 9 | 11 | 9 | 11 | 5 |
| 6 | 6 | 6 | 6 | 10 | 12 | 10 | 12 | 6 |
| 7 | 7 | 7 | 7 | 11 | $\overline{12}$ | 11 | $\overline{12}$ | 7 |
| 8 | 8 | 8 | 8 | 12 | $\overline{11}$ | 12 | $\overline{11}$ | 8 |
| 9 | 9 | $\overline{12}$ | 9 | $\overline{12}$ | $\overline{10}$ | $\overline{12}$ | $\overline{10}$ | 9 |
| 10 | 10 | $\overline{11}$ | 10 | $\overline{11}$ | $\overline{9}$ | $\overline{11}$ | $\overline{9}$ | $\overline{9}$ |
| $\overline{12}$ | 11 | $\overline{10}$ | 11 | $\overline{10}$ | $\overline{8}$ | $\overline{10}$ | $\overline{8}$ | $\overline{8}$ |
| $\overline{11}$ | 12 | $\overline{9}$ | 12 | $\overline{9}$ | $\overline{7}$ | $\overline{9}$ | $\overline{7}$ | $\overline{7}$ |

| | | | | | | |
|---|---|---|---|---|---|---|
| 1 | 1 | 1 | 1 | 1 | 7 | 1 |
| 2 | 2 | 2 | 2 | 2 | 8 | 2 |
| 3 | 3 | 7 | 9 | 7 | 9 | 3 |
| 4 | 4 | 8 | 10 | 8 | 10 | 4 |
| 5 | 5 | 9 | 11 | 9 | 11 | 5 |
| 6 | 6 | 10 | 12 | 10 | 12 | 6 |
| 7 | 7 | 11 | $\overline{12}$ | 11 | $\overline{12}$ | 7 |
| 8 | 8 | 12 | $\overline{11}$ | 12 | $\overline{11}$ | 8 |
| 9 | 9 | $\overline{12}$ | $\overline{10}$ | $\overline{12}$ | $\overline{10}$ | 9 |
| 10 | 10 | $\overline{11}$ | $\overline{9}$ | $\overline{11}$ | $\overline{9}$ | $\overline{9}$ |
| $\overline{12}$ | 11 | $\overline{10}$ | $\overline{8}$ | $\overline{10}$ | $\overline{8}$ | $\overline{8}$ |
| $\overline{11}$ | 12 | $\overline{9}$ | $\overline{7}$ | $\overline{9}$ | $\overline{7}$ | $\overline{7}$ |
| 1 | 1 | 1 | 1 | 1 | 5 | 1 |
| 2 | 2 | 2 | 2 | 2 | 6 | 2 |
| 3 | 3 | 9 | 7 | 9 | 7 | 3 |
| 4 | 4 | 10 | 8 | 10 | 8 | 4 |
| 5 | 5 | 11 | 9 | 11 | 9 | 5 |
| 6 | 6 | 12 | 10 | 12 | 10 | 6 |
| 7 | 7 | $\overline{12}$ | 11 | $\overline{12}$ | 11 | 7 |
| 8 | 8 | $\overline{11}$ | 12 | $\overline{11}$ | 12 | 8 |
| 9 | 9 | $\overline{10}$ | $\overline{12}$ | $\overline{10}$ | $\overline{12}$ | $\overline{8}$ |
| 10 | 10 | $\overline{9}$ | $\overline{11}$ | $\overline{9}$ | $\overline{11}$ | $\overline{7}$ |
| $\overline{12}$ | 11 | $\overline{8}$ | $\overline{10}$ | $\overline{8}$ | $\overline{10}$ | $\overline{6}$ |
| $\overline{11}$ | 12 | $\overline{7}$ | $\overline{9}$ | $\overline{7}$ | $\overline{9}$ | $\overline{5}$ |

左から順に $(k_{10}, k_8, k_6, k_4, k_2, k_0) = (2, \underline{3}, 0, 0, 3, 1),\ (2, \underline{1}, 0, 0, 3, 1),\ (2, 0, \underline{1}, 0, 3, 1)$ に対する $\mathrm{fill}(u_\lambda)$。下線を引いた文字が $k_c$ を表す。$k_h$ に対応する列の文字 $1, 2, \ldots, h$ を対応する各 Step 2, 3, 4 に応じてピンク、黄色、緑で表した。

**$\Phi^{-1}$ のアルゴリズム**　基本的には $A_n^{(1)}$ 型のアルゴリズムと類似であるが、46 ページに与えた古典部分代数 $D_n$ のクリスタルグラフにそって実行されるため以下のような相違点が生じる。

$\mu^{(a)}$ のある行 $\mu_i^{(a)}$ を選択し、$B^{a,\mu_i^{(a)}}$ の元を以下のようにして定める。$\ell^{(a-1)} = \mu_i^{(a)}$ とする。次いで $\nu^{(a)}$ の特異なストリングで長さが $\ell^{(a-1)}$ 以上で最も短いものを選び、その長さを $\ell^{(a)}$ と定める。再帰的に以上の手続きを繰り返していき、途中 $\nu^{(j)}$ で初めて該当する特異なストリングが選べなくなった場合、$\ell^{(j)} = \infty$ としてそこでストップし、出力の文字を $k = j$ とする。以上は $A_n^{(1)}$ 型の場合と同一であるが、$\nu^{(n-2)}$ から先は Dynkin 図ないし $D_n$ のクリスタルグラフが二股に分かれることから以下のような手続きとなる。

$\ell^{(n-2)} < \infty$ だったとしよう。その時 $(\nu^{(n-1)}, J^{(n-1)})$ または $(\nu^{(n)}, J^{(n)})$ の特異なストリングで $\ell^{(n-2)}$ 以上の長さで最小のものを探し、それぞれの長さから同様にして $\ell^{(n-1)}$ および $\ell^{(n)}$ を定める（該当するものがなければ $\infty$ とする）。

1. $\ell^{(n-1)} = \infty$ かつ $\ell^{(n)} = \infty$ ならば出力を $k = n-1$ として停止。
2. $\ell^{(n-1)} < \infty$ かつ $\ell^{(n)} = \infty$ ならば出力を $k = n$ として停止。
3. $\ell^{(n-1)} = \infty$ かつ $\ell^{(n)} < \infty$ ならば出力を $k = \bar{n}$ として停止。
4. $\ell^{(n-1)} < \infty$ かつ $\ell^{(n)} < \infty$ ならば $\ell_{(n-1)} = \max(\ell^{(n-1)}, \ell^{(n)})$ として継続する。

以降は手続きの左右を入れ替えて繰り返す。すなわち $\ell_{(j+1)}$ が定まってるとき $(\nu^{(j)}, J^{(j)})$ の特異なストリングで $\ell_{(j+1)}$ 以上のもののうち最小のものの長さを $\ell_{(j)}$ と定義する。初めて $\ell_{(j)} = \infty$ となったとき出力を $k = \overline{j+1}$ として停止する。一方最終的に $\ell_{(1)} < \infty$ となったときには出力を $k = \bar{1}$ として停止する。

その後のステップは $A_n^{(1)}$ 型の場合と全く同一である。すなわち選択したストリングから一ますずつ削り、新しい rigging は、削られたストリングについては新しい艤装配位で特異なストリングとなるように定め、一方削られなかったストリングの rigging は変更しない。

**例 1**　$D_5^{(1)}$ とする。以下の $B^{3,2} \otimes B^{2,2}$ 型艤装配位を考えよう。

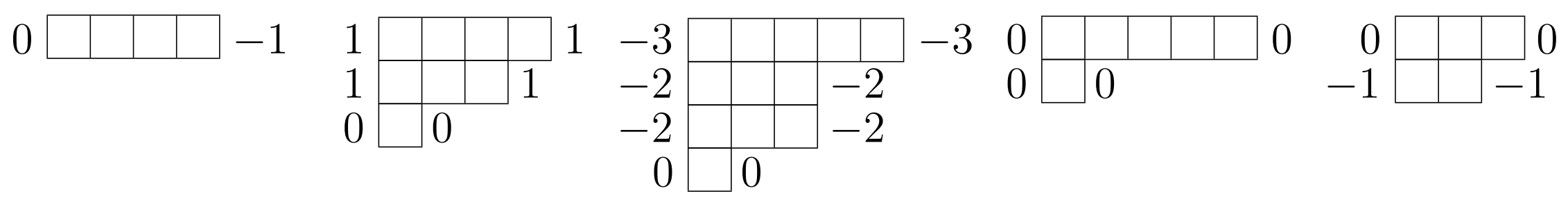

$\mu^{(a)}$ は略してあるが、今のテンソル積の形状から $\mu^{(2)} = (2)$ および $\mu^{(3)} = (2)$ である。計算を見やすくするために以下の記号を使用する。$\mu^{(a)}$ の長さ $j$ の行から削るとき（つまり $\ell^{(a-1)} = j$ のとき）、ます目を選択してから一ますずつ削除し新しい rigging を決定するまでの一連の操作を $\delta_j^{(a)}$ と表すことにしよう。与えられた艤装配位で $\mu^{(a)}$ の長さ $j$ の行を選択すると順に $\delta_j^{(a)}, \delta_1^{(a-1)}, \delta_1^{(a-2)}, \ldots, \delta_1^{(1)}$ を行うことで $B^{a,j}$ の一つの列が得られることになる。次いで $\delta_{j-1}^{(a)}, \delta_1^{(a-1)}, \delta_1^{(a-2)}, \ldots, \delta_1^{(1)}$ を行うと二列目が得られる、等となる。この手順は $A_n^{(1)}$ 型の場合と共通である。

ここでは $\mu^{(2)}$ の長さ 2 の行を選択して $B^{2,2}$ の元を得るところから始めよう。すなわち最初のステップは $\delta_2^{(2)}$ であり、$\nu^{(2)} = (4,3,1)$ の長さ 2 かそれ以上の特異なストリングを探すことから開始する。同様の操作を続けていって $B^{2,2}$ に対する全ての手順を実行すると以下のようになる。ます目に "×" のしるしを付したものが選択され除去されるます目である。また $\delta_j^{(a)}$ と矢印を挟んで途中段階の像を記述した。

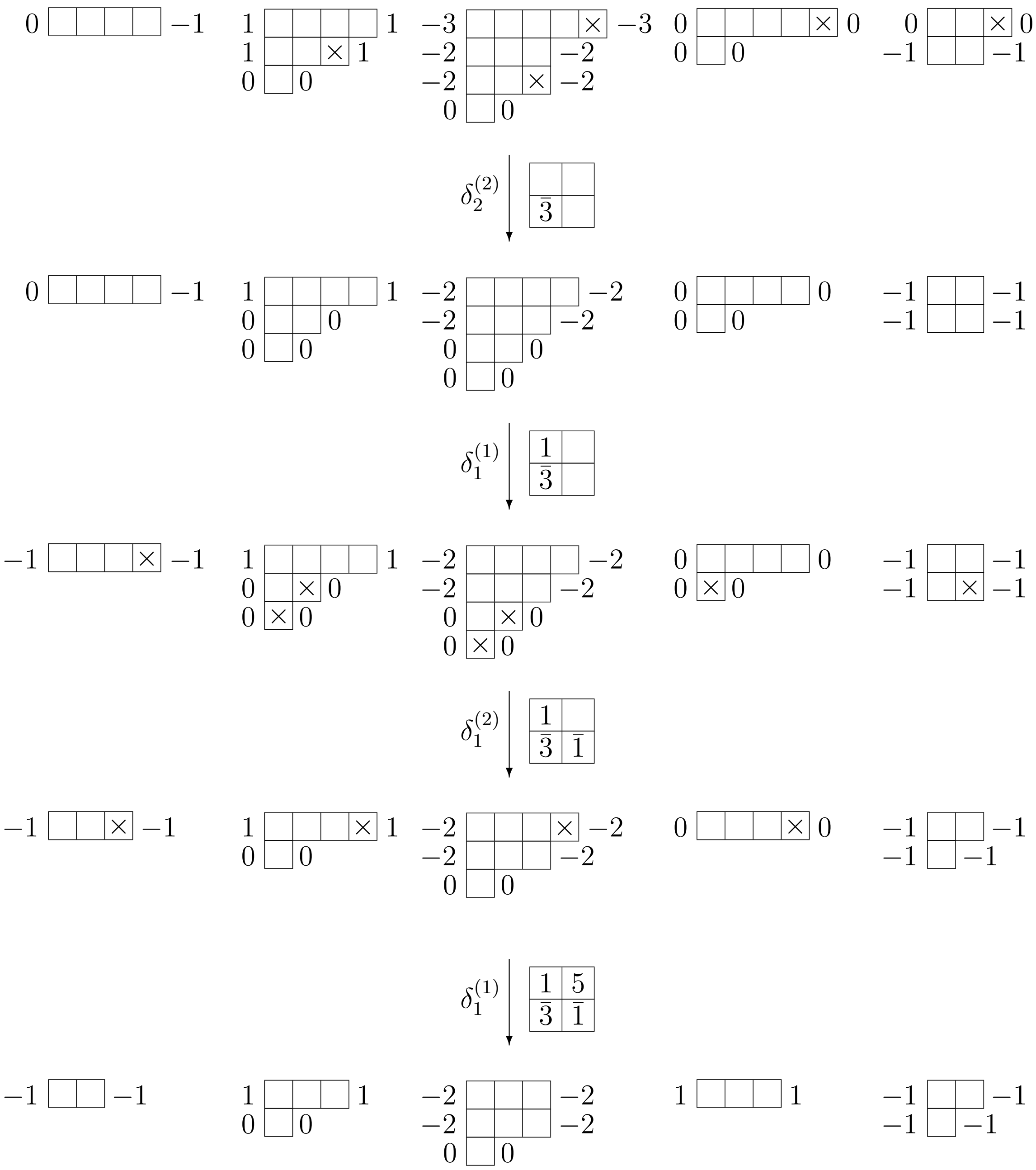

いくつか注釈をつけておく。

- $\delta_2^{(2)}$ の後新しい vacancy numbers を計算する時には path が $B^{3,2} \otimes B^{2,1} \otimes B^{1,1}$ であるとして計算する。つまり $\mu^{(1)} = (2,1)$, $\mu^{(2)} = (1)$ とする。
- 最初の $\delta_1^{(1)}$ は $\nu^{(1)}$ に特異なストリングが存在しないため削れずに出力が 1 である。
- $\delta_1^{(2)}$ は $\nu^{(2)}$ から削り始めて右端まで行って折り返した後 $\nu^{(1)}$ まで削っている。

残りの計算も実行すると、最終的な像 $\Phi^{-1}(\nu, J)$ は以下の通り。

$$\Phi^{-1}_{B^{3,2}\otimes B^{2,2}}(\nu, J) = \begin{array}{|c|c|}\hline 1 & \bar{5} \\ \hline 4 & \bar{3} \\ \hline 5 & \bar{1} \\ \hline \end{array} \otimes \begin{array}{|c|c|}\hline 1 & 5 \\ \hline \bar{3} & \bar{1} \\ \hline \end{array}. \tag{55}$$

参考までに $\mu^{(a)}$ の選び方の順を逆にすると以下の像が得られる。

$$\Phi^{-1}_{B^{2,2}\otimes B^{3,2}}(\nu,J)=\begin{array}{|c|c|}\hline 2&\bar{5}\\\hline 5&\bar{3}\\\hline\end{array}\otimes\begin{array}{|c|c|}\hline 1&4\\\hline 5&\bar{2}\\\hline \bar{3}&\bar{1}\\\hline\end{array}. \tag{56}$$

籏装配位写像の $R$ 不変性により (55) と (56) は $R$-行列により移りあう。

$$R:\begin{array}{|c|c|}\hline 1&5\\\hline \bar{3}&\bar{1}\\\hline\end{array}\otimes\begin{array}{|c|c|}\hline 1&\bar{5}\\\hline 4&\bar{3}\\\hline 5&\bar{1}\\\hline\end{array}\longmapsto\begin{array}{|c|c|}\hline 1&4\\\hline 5&\bar{2}\\\hline \bar{3}&\bar{1}\\\hline\end{array}\otimes\begin{array}{|c|c|}\hline 2&\bar{5}\\\hline 5&\bar{3}\\\hline\end{array}.$$

よってこの場合 $R$-行列は $R=\Phi^{-1}_{B^{3,2}\otimes B^{2,2}}\circ\Phi_{B^{2,2}\otimes B^{3,2}}$ と実現される。長さ 2 以上のテンソル積に対して任意の回数 $R$-行列を作用させたものでも同様に $R=\Phi^{-1}_{B'}\circ\Phi_B$ と一回で計算できる。計算の効率だけではなく、組み合わせ $R$-行列を定義から原理的に計算を行ったり、又はそれを改良したような再帰的な方法ではなく、明示的なアルゴリズムによって計算する方法はおそらく存在していなかったと思われる。アルゴリズムが明示的に与えられたからこそ箱玉系の逆散乱形式などの重要な結果が導出できるのである。

**例 2**　今考えた例は出てくるタブローが Kashiwara–Nakashima タブローばかりだったので、そうではない例をひとつあげる。$D^{(1)}_6$ とする。$B^{4,3}\otimes B^{2,3}\simeq B^{2,3}\otimes B^{4,3}$ 型のテンソル積についての以下の最高ウェイト籏装配位を考える。

$$\emptyset\qquad\begin{array}{lll}0&\square\square&0\\0&\square&0\end{array}\qquad\begin{array}{lll}1&\square\square\square&0\\0&\square&0\\0&\square&0\end{array}\qquad\begin{array}{lll}0&\square\square\square&0\\0&\square\square\square&0\\0&\square&0\\0&\square&0\end{array}\qquad\begin{array}{lll}0&\square\square\square&0\\0&\square&0\end{array}\qquad\begin{array}{lll}0&\square\square\square&0\\0&\square&0\end{array}$$

二通りの方法で籏装配位写像 $\Phi^{-1}$ を計算すると以下の像を得る。

$$R:\begin{array}{|c|c|c|}\hline 1&1&1\\\hline 2&2&2\\\hline \bar{4}&3&3\\\hline \bar{3}&4&\bar{3}\\\hline\end{array}\otimes\begin{array}{|c|c|c|}\hline 1&1&1\\\hline 2&3&\bar{2}\\\hline\end{array}\simeq\begin{array}{|c|c|c|}\hline 1&1&1\\\hline 2&2&2\\\hline\end{array}\otimes\begin{array}{|c|c|c|}\hline 1&1&1\\\hline 2&4&3\\\hline 4&\bar{4}&4\\\hline \bar{4}&\bar{2}&\bar{4}\\\hline\end{array}.$$

$\begin{array}{|c|c|c|}\hline 1&1&1\\\hline 2&4&3\\\hline 4&\bar{4}&4\\\hline \bar{4}&\bar{2}&\bar{4}\\\hline\end{array}$ に対応する最高ウェイト元は $\begin{array}{|c|c|c|}\hline 1&1&1\\\hline 2&2&2\\\hline 3&\bar{4}&3\\\hline 4&\bar{3}&4\\\hline\end{array}=\mathrm{fill}\left(\begin{array}{|c|c|c|}\hline 1&1&1\\\hline 2&2&2\\\hline 3\\\cline{1-1} 4\\\cline{1-1}\end{array}\right)$。最高ウェイト元を得るために作用させた $\widetilde{e}_a$ の列を逆順にして $\tilde{f}_a$ の列としてかければ対応する Kashiwara–Nakashima タブローが得られる。結局上記 $R$-行列の結果を KN タブロー表示すれば

$$R:\begin{array}{|c|c|c|}\hline 1&1&1\\\hline 2&2&2\\\hline\end{array}\otimes\begin{array}{|c|c|c|}\hline 1&1&1\\\hline 2&3&\bar{2}\\\hline\end{array}\simeq\begin{array}{|c|c|c|}\hline 1&1&1\\\hline 2&2&2\\\hline\end{array}\otimes\begin{array}{|c|c|c|}\hline 1&1&1\\\hline 2&3&\bar{2}\\\hline 4\\\cline{1-1} \bar{4}\\\cline{1-1}\end{array}$$

という形になる。

籏装配位写像 $\Phi$ の定義をよくご理解いただければ、KR タブローの様に風変わりな、しかし代数構造を持つ対象が発生することは誠に驚異的なこととご納得頂けると思う。しかし実際にはなかなか受け入れがたい事実であるようで、かつて関係者の方々に概念を正しく理解して頂くのに大変な労力を要したことを思い出す。

## 8.9 Littlewood–Richardson タブロー

本稿の締めくくりにデザートのような話題を提供しようと思う[28]。この話題における一つの主役は Littlewood–Richardson（LR）タブローである。普通 LR タブローを扱う時はタブローそのものではなくある条件を満たす LR タブローの総数について関心をもたれる場合が多いのであるが、以下で述べる構成では LR タブローそのものに重要な意味があるところが興味深い。

考察する問題をより具体的に述べると、任意の非例外型アフィンリー代数 $\mathfrak{g}$ の最高ウェイト表現と $A_n^{(1)}$ 型の最高ウェイト表現および LR タブローの組との間には興味深い全単射が存在する。そのような全単射は母関数のレベルでの恒等式として Shimozono–Zabrocki [SZ] により提唱され、また論文 [Sh05, LS05, LOS] において恒等式の証明が行われた。その様な恒等式の成り立つ数学的根拠として論文 [OS] において具体的に表現の間の全単射が構成されたので、以下はその内容を紹介する。Cartan データを含む詳細は原論文に譲ることとして、議論の本質的な部分をあまり細部に立ち入らずに概観することを目標とする。

なお本節では代数として一般の非例外型アフィンリー代数を考察する。代数のランクが低い場合を除き全ての代数について普遍的な構成となっている。

**LR タブロー**　まずは LR タブローの復習から。$T$ が LR タブローであるとは以下の二つの条件を満たす時をいう。

- $T$ は歪準標準盤である。つまり $T$ は左上のある Young 図の部分を空白として他の部分に数字を書き込んだタブローであり、書き込まれた文字について通常の準標準盤としての条件を課す。

- $T$ に書き込まれた文字を右から左、上から下へ、という順に読んで得られる数列を row word と呼ぶとき、row word $w = x_1 \cdots x_l$ が Yamanouchi 条件を満たす。ここで Yamanouchi 条件とは $w$ を左端から任意の部分まで読んで得られる任意の部分列 $x_1 \cdots x_k$ に対して、部分列に含まれる文字 1 の総数が文字 2 の総数以上であり、文字 2 の総数が文字 3 の総数以上であり、等々という条件が満たされることをいう[29]。

LR タブロー $T$ の型 $\eta/\lambda$ とは、$T$ の外形が $\eta$ であり、左上の文字が書き込まれない領域の形が $\lambda$ であることをいう。また LR タブロー $T$ のウェイト $\mu = (\mu_1, \ldots, \mu_l)$ とは、$T$ に含まれる文字 $j$ の個数を $\mu_j$ として得られる分割の事である。

例えば以下のタブローは $(43^31^2)/(2^21)$ 型の LR タブローであり、ウェイトは $(3,3,2,2)$ である。

| | | | |
|---|---|---|---|
| | | 1 | 1 |
| | | 2 | |
| | 1 | 3 | |
| 2 | 2 | 4 | |
| 3 | | | |
| 4 | | | |

この場合の row word は 1123142234 であり、確かに Yamanouchi 条件を満たす。

[28]大分こってりした味付けで胃もたれしてしまうかもしれないが。

[29]Yamanouchi 条件は $A_n$ 型クリスタルのテンソル積 $(B^{1,1})^{\otimes L}$ の元が最高ウェイト元であるための条件と同一である。

型 $\eta/\lambda$ でウェイト $\mu$ の LR タブローの総数は Littlewoood–Richardson 係数と呼ばれ、通常 $c^{\eta}_{\lambda\mu}$ と表される。良く知られているように $c^{\eta}_{\lambda\mu}$ は、$\mathfrak{gl}_n$ のウェイト $\lambda$ の既約表現を $V_\lambda$ と書くとき、テンソル積 $V_\lambda \otimes V_\mu$ の中に含まれる表現 $V_\eta$ の重複度に等しい。

**艤装配位の安定化**　元々考えられていた問題では代数 $\mathfrak{g}$ に対する KR クリスタルのテンソル積 $\bigotimes_k B^{r_k,s_k}$ を考える事になるが、我々の立場では艤装配位上で考察することになる。いくつか用語を準備しよう。KR クリスタルの議論では非例外型アフィンリー代数を Dynkin 図の 0 番目の頂点付近の形状に応じて以下のように分類すると便利なことが多い。

| $\diamondsuit$ | $\diamondsuit$ 型の $\mathfrak{g}$ |
|---|---|
| $\varnothing$ | $A^{(1)}_n$ |
| $\square$ | $D^{(2)}_{n+1}, A^{(2)}_{2n}$ |
| $\square\square$ | $C^{(1)}_n$ |
| $\genfrac{}{}{0pt}{}{\square}{\square}$ | $A^{(2)}_{2n-1}, B^{(1)}_n, D^{(1)}_n$ |

ここで代数の型を表す $\diamondsuit$ の形は、クリスタル $B^{r,s}$ を古典部分代数の表現に分解したときに現れるウェイトの形（すなわち Kashiwara–Nakashima タブローの取りうる形）が $r \times s$ の長方形から $\diamondsuit$ の形を任意個数取り除いて得られる形となることによる。以下簡単のため各型 $\diamondsuit = \varnothing, \square, \square\square, \genfrac{}{}{0pt}{}{\square}{\square}$ について

$$\mathfrak{g}^{\diamondsuit} = A^{(1)}_n, D^{(2)}_{n+1}, C^{(1)}_n, D^{(1)}_n \tag{57}$$

という場合を念頭に置いて議論するが、この選択は必須のものではない[30]。なお特に断りがない場合、本節での艤装配位は最高ウェイト元の場合を考える。

$\diamondsuit$ 型の代数に対する艤装配位を考えよう。詳細は論文 [OSS03] に譲るが、一般の代数に対する艤装配位も Dynkin 図の 0 番を除く各頂点に $\nu^{(a)}$ を配置して

$$(\nu, J) = \left(\mu^{(1)}, \ldots, \mu^{(n)}, (\nu^{(1)}, J^{(1)}), \ldots, (\nu^{(n)}, J^{(n)})\right) \tag{58}$$

の様に得られる。ここで $\mu^{(a)}$ の番号付は $D^{(1)}_n$ 型と同様に行 $\mu^{(a)}_i$ はクリスタル $B^{a,\mu^{(a)}_i}$ に対応するとしよう。記号 $a^{\diamondsuit}$ を以下のように定める。

$$a^{\square\square} = a^{\square} = n-1, \qquad a^{\genfrac{}{}{0pt}{}{\square}{\square}} = n-2.$$

その心は Dynkin 図の真ん中にある $A$ 型的な頂点のうち最大のものを表している。一般の代数の場合の vacancy numbers も $A$ 型的な頂点に対しては $A^{(1)}_n$ 型艤装配位と同様に定義される（例えば [OS, §2.3] を参照されたい）。

興味深い点として与えられた path に対し考察する代数のランク $n$ を変化させていくと、$n$ が十分小さい場合を除き艤装配位が安定化し、特徴的な振る舞いをするようになる [OS, §2.4]。具体的には、ある $k$ が存在して

$$\nu^{(k)} = \nu^{(k+1)} = \cdots = \nu^{(a^{\diamondsuit})} \tag{59}$$

となり、かつ $\nu^{(a^{\diamondsuit})}$ は $\diamondsuit$ を単位として敷き詰めることのできる形状となる（$A^{(1)}_n$ 型なら空集合となる）。この時 (59) に対する vacancy numbers は全て 0 となり、艤装配位に対する最

30 一方 [LOS] の議論においては必須のものだった

高ウェイト条件から riggings は全て正となることと合わせると、実質的にこの部分に対する rigging の自由度は消滅する。この様な状況下では艤装配位の $\nu^{(a^\diamondsuit)}$ より右側の部分に現れる各代数への依存性は無視することができてしまい、艤装配位は型 $\diamondsuit$ のみで分類されることになる。

艤装配位を用いてもウェイトを計算することができ、例えば [Sa14, Eq.(13)] や [OS, Eq.(2.7)] に一般的な式が与えられているが、さしあたって以下の議論を理解するために代数のランクが十分大きくて安定化 (59) が成立している場合の式を示す。その場合 $\lambda = \mathrm{wt}(\nu, J)$ を Young 図で表した時の $a$ 行目 $\lambda_a$ は

$$\lambda_a = \sum_{b \geq a} \left|\mu^{(b)}\right| + \left|\nu^{(a-1)}\right| - \left|\nu^{(a)}\right| \tag{60}$$

で与えられる。ここで分割 $\lambda$ に対し $|\lambda|$ でます目の総数を表す。以下型 $\diamondsuit$ でウェイト $\lambda$ の（最高ウェイト）艤装配位の全体を $\mathrm{RC}^{\diamondsuit}(\lambda)$ で表すことにする。path の形状を表す $\mu^{(a)}$ は固定しておき、特に明示しないことにする。

**全単射の定式化**　本題の全単射を定義しよう。いくつか記号を準備する。$\mathrm{LR}^{\eta}_{\lambda\mu}$ で型 $\eta/\lambda$ をもつウェイト $\mu$ の LR タブローの全体としよう。また分割 $\lambda$ に対し $\ell(\lambda)$ で分割の長さを表すことにしよう。

考察する全単射は

$$\Psi : \mathrm{RC}^{\diamondsuit}(\lambda) \longrightarrow \mathrm{RC}^{\varnothing}(\eta) \times \mathrm{LR}^{\eta}_{\lambda\mu} \tag{61}$$

というものである。ここで各要素の対応が

$$\Psi : (\nu, J) \longmapsto \{(\nu', J'), T\}$$

であるとき

$$\lambda = \mathrm{wt}(\nu, J), \qquad \mu = \nu^{(a^\diamondsuit)}, \qquad \eta = \mathrm{wt}(\nu', J')$$

である。$\diamondsuit = \varnothing$ すなわち $A^{(1)}_n$ 型の場合も $\lambda = \eta$ として含まれている。この全単射は抽象的なものではなく具体的なものであり、本節末で述べるように代数の持っている双対性とも関連した深い内容を持ったものである。

全単射 (61) が成立するためには、代数のランクは安定化 (59) が成立しているようなものであれば十分であるが、より精密に述べることができる。一つ艤装配位 $(\nu, J) \in \mathrm{RC}^{\diamondsuit}(\lambda)$ が与えられたとしよう。この時代数のランク $n$ が条件

$$a^{\diamondsuit} \geq \ell(\mathrm{wt}(\nu, J)) + \ell\left(\nu^{(a^\diamondsuit)}\right) \tag{62}$$

を満たす時、写像 $\Psi$ が成立する。この条件の意味は、具体的なアルゴリズムを見て頂ければご理解いただけると思うが、大体以下のようなものである。簡単のため $D^{(1)}_n$ 型の場合で考えよう。LR タブロー $T \in \mathrm{LR}^{\eta}_{\lambda\mu}$ に対し $N_{out} = \ell(\eta)$ および $N_{in} = \ell(\lambda)$ と定義しよう。この時アルゴリズムの定義から $\Psi^{-1}$ は代数のランク $n$ が $N_{out} + 2 \leq n$ を満たす時定義される。一方 $N_{out}$ と $N_{in}$ の差の最大値は $N_{out} - N_{in} = \ell(\nu^{(a^\diamondsuit)})$ なので合わせると $N_{in} + \ell(\nu^{(a^\diamondsuit)}) + 2 \leq n$ すなわち (62) が得られる。割と精密な評価であるから必要に応じて定理が成立するための代数のランクに関する条件を導出するのに利用できる（例えば [OS, Remark 4.1] を参照）。

艤装配位写像 $\Phi$ : path $\mapsto$ rigged configuration が確立されている場合は全単射 $\Psi$ はクリスタルのテンソル積 $\bigotimes_k B^{r_k,s_k}$ のレベルに拡張される。その様な写像を直接求めることは非常に難しい問題のようであり [Sh05]、おそらく一般には未解決であると思われる。もちろん艤装配位上に直接柏原作用素を導入することもできるので、表現の間の全単射という意味では (61) で完成しているといっても良いが、クリスタルのテンソル積のレベルで実現されることも重要であろう。

**アルゴリズムの定義**　では写像 $\Psi$ のアルゴリズムの記述をしよう。大雑把にいってアルゴリズムは $A_n^{(1)}$ 型の艤装配位写像のアルゴリズムをそのまま左右を入れ替えたものである。これは節末で関連を述べる代数の双対性が Dynkin 図の左右を入れ替えるものであるのと見事に整合する。また LR タブローは写像の recording tableau としてごく自然に現れる。考えている問題の複雑さから見るとあっけないほど単純になってしまい、艤装配位の底知れぬ強力さに驚かされる。

アルゴリズムは $\nu^{(a^\diamondsuit)}$ の各列を右から順に見て以下のように行う（$\nu^{(n)}$ 等には自然に拡張されるが本質的ではないので省略）。$\nu^{(a^\diamondsuit)}$ の右端の列が左から数えて $l$ 列目であり、高さが $h_l$ であるとしよう。その時操作 $\delta_l$ が以下のようにして定まる。

1. $\ell^{(a^\diamondsuit)}=l$ とし、再帰的に $\ell^{(a)}$ を定める。$\ell^{(a)}$ が定まっているとき、$\nu^{(a-1)}$ の特異なストリングで長さが $\ell^{(a)}$ 以上のもののうち最小のものの長さを $\ell^{(a-1)}$ と定める。その様なものがないとき、$\ell^{(a-1)}=\infty$ とし、出力を $k=a$ として停止する。一方 $\ell^{(1)}<\infty$ ならば出力を $k=1$ として停止する。

2. 前項で選んだ $\nu^{(a)}$ の各行の右端から一ますずつ削る。新しい rigging は、削られなかった行に対しては元と同じとし、一方削られた行に対しては新しい艤装配位において特異なストリングとなるように定める。ただし vacancy number の計算において path の形状を表す $\mu^{(a)}$ は変化させない。

3. recording tableau $T$ には以下の要領で文字を書き込む。写像 $\Psi$ の過程では $\delta_l$ を $h_l$ 回繰り返して艤装配位に作用させる。もし考えている $\delta_l$ が $\delta_l^{h_l}$ の中の $j$ 番目であったとき、$T$ の $k$ 行目の右端に文字 $j$ を書き込む。

$\nu^{(a^\diamondsuit)}$ の $i$ 列目の高さが $h_i$ である時、写像 $\Psi$ は

$$\Psi(\nu,J)=\{(\nu',J'),T\}=\delta_1^{h_1}\circ\delta_2^{h_2}\circ\cdots\circ\delta_l^{h_l}(\nu,J)$$

で与えられる。

写像 $\Psi$ は各ステップで逆にすることができ、それが逆写像 $\Psi^{-1}$ を与える。

**写像の例**　ここでは $D_n^{(1)}$ $(n\geq 8)$ 型 $(B^{1,3})^{\otimes 3}\otimes(B^{1,2})^{\otimes 2}\otimes(B^{1,1})^{\otimes 2}$ テンソル積の以下の元を考える：

$$b=\boxed{1}\boxed{1}\boxed{1}\otimes\boxed{2}\boxed{\bar{1}}\boxed{\bar{1}}\otimes\boxed{1}\boxed{2}\boxed{\bar{2}}\otimes\boxed{2}\boxed{3}\otimes\boxed{2}\boxed{\bar{2}}\otimes\boxed{\bar{2}}\otimes\boxed{2}.$$

その時対応する艤装配位上で写像 $\Psi$ は以下のように計算される。以下の図では次のような記述をする。一番最初の艤装配位は上で与えた元 $b$ に対応する。紫色の印をつけた箱は矢印の左側で指定されたそれぞれの操作 $\delta$ によって取り除かれるます目を表す。各段階で得られるレコーディングタブロー $T$ はそれぞれの矢印の右側に与えた。

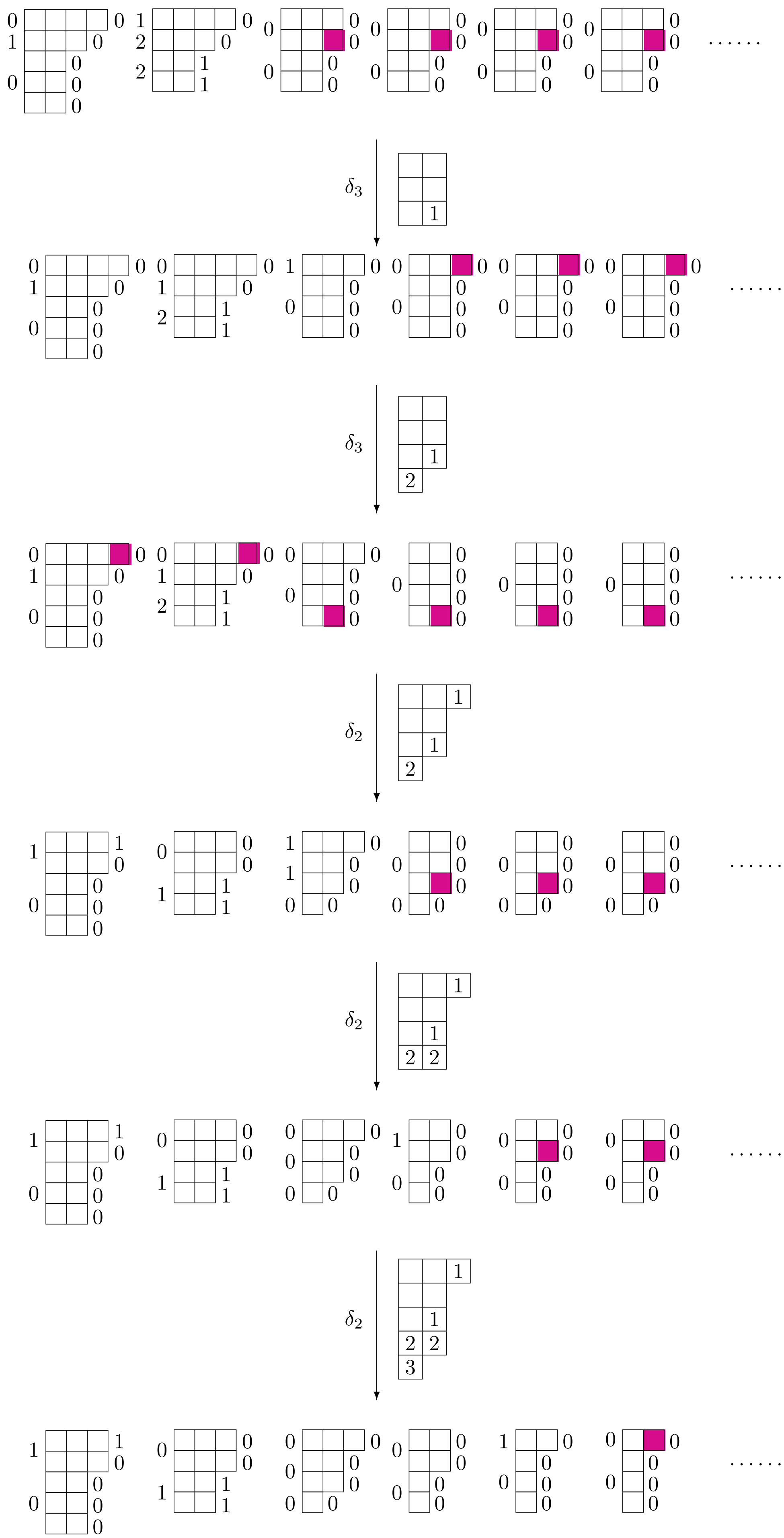

$\delta_3$
$\delta_3$
$\delta_2$
$\delta_2$
$\delta_2$

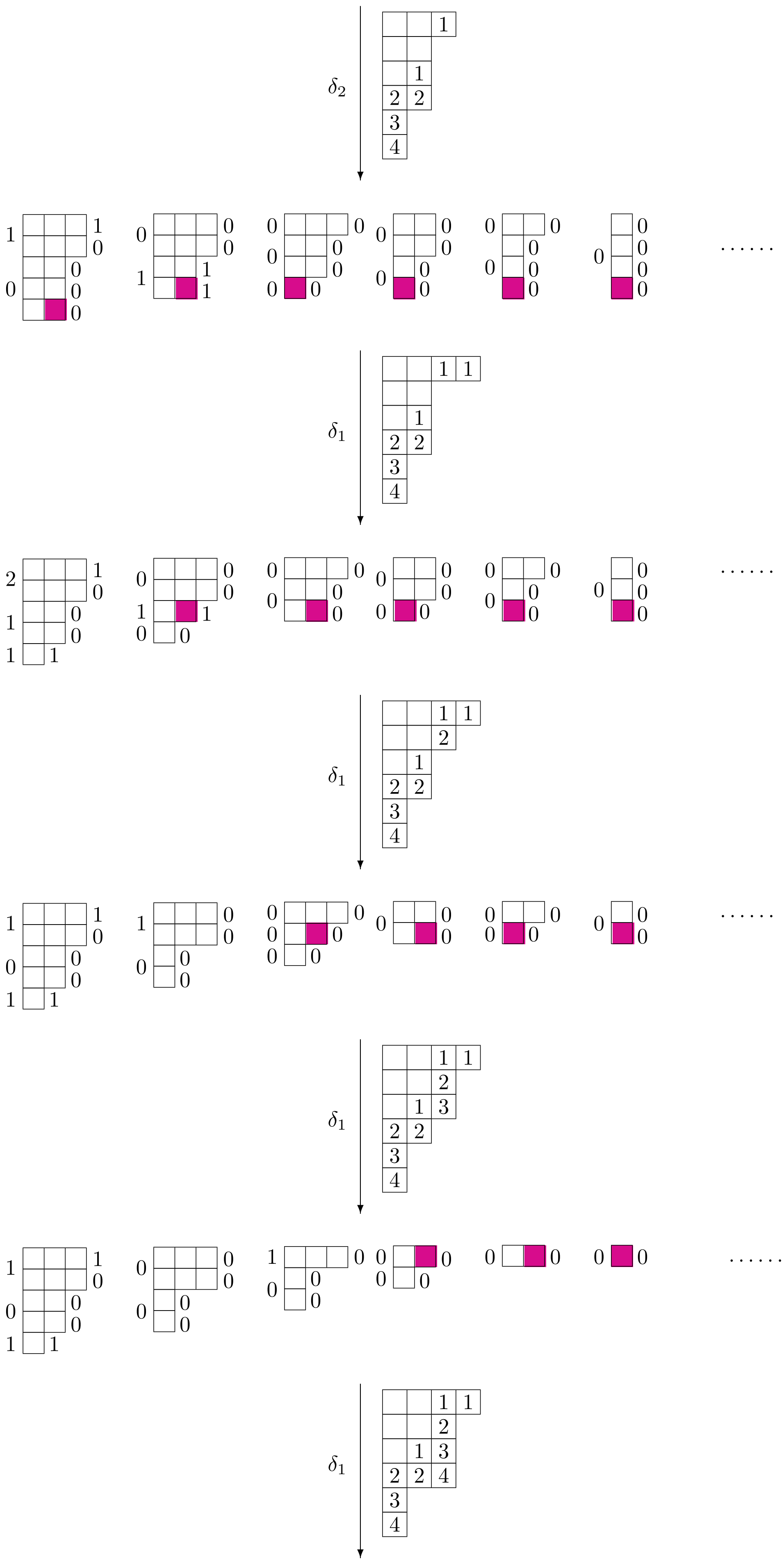
$\delta_2$
$\delta_1$
$\delta_1$
$\delta_1$

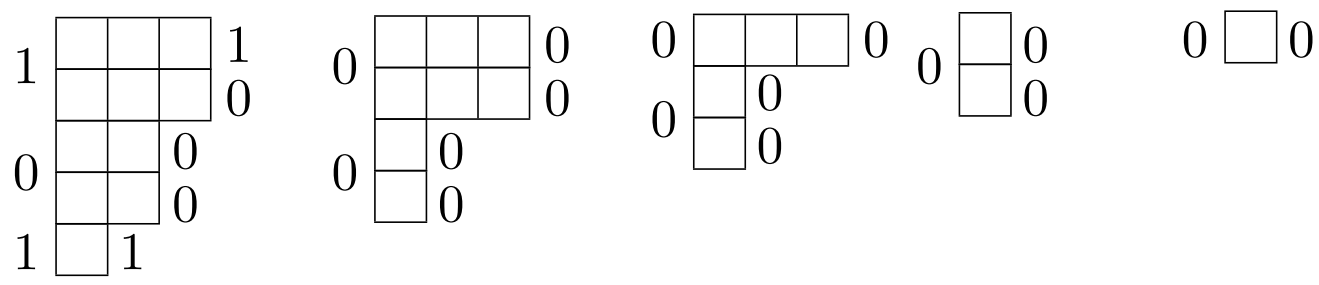

図式の最後に与えられた艤装配位および $T$ が今の場合の写像 $\Psi$ の像となる。特に $T$ は準標準盤となっており、更にその row word は 1123142234 となっていて確かに Yamanouchi 条件を満たしていることが確かめられる。なお艤装配位写像 $\Phi^{-1}$ のもとで最終的な艤装配位は以下のテンソル積の元に対応している：

$$b' = \fbox{1}\fbox{1}\fbox{1} \otimes \fbox{2}\fbox{2}\fbox{2} \otimes \fbox{1}\fbox{3}\fbox{3} \otimes \fbox{4}\fbox{4} \otimes \fbox{3}\fbox{5} \otimes \fbox{4} \otimes \fbox{6}.$$

なお逆写像 $\Psi^{-1}$ の計算は上記の計算を完全に逆にたどっていけばよい。

**双対性**　$b$ および $b'$ は先ほどの例で与えられたものと同じとし、それらを $D_8^{(1)}$ の元であるとみなす。これらを文献 [LOS] 第 5.3 節に与えられた対合 (involution) $\sigma$ と比較してみよう。$\sigma$ は Dynkin 図の左右を入れ替える双対性、すなわち

$$\sigma \circ \widetilde{e}_a = \widetilde{e}_{n-a} \circ \sigma$$

が任意の $a \in I = \{0, 1, \ldots, n\}$ について成り立つものとしよう。今考えている path $b$ について具体的に計算すると、

$$\sigma(b) = \fbox{$\bar{8}$}\fbox{$\bar{8}$}\fbox{$\bar{8}$} \otimes \fbox{8}\fbox{8}\fbox{$\bar{7}$} \otimes \fbox{6}\fbox{$\bar{8}$}\fbox{$\bar{6}$} \otimes \fbox{$\bar{7}$}\fbox{$\bar{6}$} \otimes \fbox{6}\fbox{$\bar{6}$} \otimes \fbox{7} \otimes \fbox{$\bar{7}$}.$$

すると $b'$ は $\sigma(b)$ に対応する $I_0$-最高ウェイト元と一致する（$I_0 = I \setminus \{0\}$）。同様の関係が任意の $\Psi$ の像に対して成り立っているものと予想している。

なお $\sigma$ の定義では代数として (57) に与えたものを取る必要があるが、艤装配位上の全単射 $\Psi$ の定義では任意の非例外型アフィンリー代数を選ぶことができ、その定式化は代数の種類にほとんど依存しない普遍性のあるものであった。これもまた艤装配位のもつ顕著に深い性質の一つであり、艤装配位自身が数学的に意味のある量であることの証左となっている。

## 8.10　その他の話題

本稿では触れられなかったが興味深い結果について簡単に紹介して本稿の〆とする。

**他の代数への拡張**　艤装配位写像 $\Phi$ は様々な代数に対して拡張されている。現状を知るために筆者の知る範囲での星取表を以下に載せておく。

- $A_n^{(1)}$ 型：完全に一般の場合に完成 [KSS]。
- $D_n^{(1)}$ 型：完全に一般の場合に完成 [OSSS]。
- その他の非例外型代数：[OSS03] により $(B^{1,1})^{\otimes L}$ 型のテンソル積について完成。
- $E_6^{(1)}$ 型：[OSano] により $(B^{1,1})^{\otimes L}$ 型のテンソル積について完成。
- $D_4^{(3)}$ 型：[Scr15] により $(B^{1,1})^{\otimes L}$ 型のテンソル積について完成。

様子を見るために代数が simply laced ではない場合の例として $C_n^{(1)}$ 型の場合を紹介することにしよう。[OSS03] により完成されているのは $(B^{1,1})^{\otimes L}$ 型の場合なので $\mu^{(1)}=(1^L)$ としよう。すると艤装配位は

$$(\nu, J)=\big((\nu^{(1)},J^{(1)}),(\nu^{(2)},J^{(2)}),\dots,(\nu^{(n)},J^{(n)})\big)$$

となる。ここで $\nu^{(n)}$ の各行の長さは偶数 $\nu_i^{(n)}\in 2\mathbb{Z}$ とする。vacancy numbers は

$$\begin{aligned}
&P_l^{(1)}(\nu)=L-2Q_l(\nu^{(1)})+Q_l(\nu^{(2)}),\\
&P_l^{(a)}(\nu)=Q_l(\nu^{(a-1)})-2Q_l(\nu^{(a)})+Q_l(\nu^{(a+1)}), \qquad (1<a<n),\\
&P_l^{(n)}(\nu)=Q_l(\nu^{(n-1)})-Q_l(\nu^{(n)}).
\end{aligned}$$

今までと異なるのは $a=n$ の場合である。最高ウェイト元の場合の艤装配位の定義は $A_n^{(1)}$ 型や $D_n^{(1)}$ 型の場合と全く同一である。すなわち全ての vacancy numbers が 0 以上であり、かつ rigging は常に 0 以上対応する vacancy number 以下であることとする。

艤装配位写像のかなめであるます目を選択し除去する手続き $\delta$ は以下のようになる。

1. $\ell^{(0)}=1$ とし以下のように再帰的に $\ell^{(a)}$ を定めていく。$\ell^{(a-1)}$ まで定まったとする。その時 $\nu^{(a)}$ の特異なストリングで長さが $\ell^{(a-1)}$ 以上のものが存在するとき、そのうち最小のものの長さを $\ell^{(a)}$ と定め、手順を続行する。もし $\nu^{(a)}$ に該当する特異なストリングが存在しない時は $\ell^{(a)}=\infty$ と定め、出力を $a$ として手続きを停止する。

2. もし $\ell^{(n)}<\infty$ であった場合、$\nu^{(n)}$ の長さ $\ell^{(n)}$ の行からは二ます並んで削ることとする（そうすると $\nu^{(n)}$ の行はいつでも偶数にできる）。$\overline{\ell}^{(n)}=\ell^{(n)}$ と定義し、$\ell^{(n)}=\overline{\ell}^{(n)}-1$ と再定義すると、$\nu^{(n)}$ の $\ell^{(n)}$ 列目と $\overline{\ell}^{(n)}$ 列目のます目が削られることになる。

3. もし $\overline{\ell}^{(a+1)}=\ell^{(a)}$ であった場合 Case (S) であると言い、以下のように該当する各 $\nu^{(a)}$ の行から二ますずつ削る。$\overline{\ell}^{(a)}=\ell^{(a)}$ と定義し、$\ell^{(a)}=\overline{\ell}^{(a)}-1$ と再定義し、$\nu^{(a)}$ の特異なストリングの右端二ます（$\ell^{(a)}$ 列目と $\overline{\ell}^{(a)}$ 列目）を並んで削る。

4. もし $\overline{\ell}^{(a+1)}>\ell^{(a)}$ であれば、$\nu^{(a)}$ の特異なストリングで長さが $\overline{\ell}^{(a+1)}$ 以上であるもののうち最少のものの長さを $\overline{\ell}^{(a)}$ と定義する。もし $\overline{\ell}^{(a+1)}<\infty$ かつ $\nu^{(a)}$ に該当する特異なストリングが存在しなかった場合、$\overline{\ell}^{(a)}=\infty$ とおき、出力を $\overline{a+1}$ として停止する。一方 $\overline{\ell}^{(1)}<\infty$ ならば出力を $\overline{1}$ として停止する。

こうして選んだ特異なストリングから一ますないし二ますずつ削り、$\mu^{(1)}=(1^{L-1})$ とする。新しい rigging は、変化しなかったストリングに対しては元と同じものを、削除されたストリングに対しては新しい艤装配位において該当するストリングが特異になるように与えることは $A_n^{(1)}$ 型や $D_n^{(1)}$ 型の場合と同様である。

一般の $\bigotimes_{k=1}^{L} B^{r_k,s_k}$ の場合にも数値的に実験がなされている [ScScr14]。結果として、$B^{r,s}$ を表示するタブローは §8.8 で与えた $D_n^{(1)}$ 型の場合と全く同様である。より詳しく述べると、クリスタル $B^{r,s}$ を古典部分代数の表現として分解すると、$r\times s$ 型の長方形から横ドミノ □□ を任意個数取り除いたものが得られる。従って $D_n^{(1)}$ 型 KR タブローの定義で現れた奇数のステップ $k_c$ は存在せず、かえって単純なものになってしまう。もちろん他の KR タブローを得るために用いる柏原作用素 $\widetilde{f}_a$ には $C_n^{(1)}$ 型のものを用いる。また艤装配位上の柏原作用素の定義も $A_n^{(1)}$ 型や $D_n^{(1)}$ 型の場合と全く同じ定義で良いようである。

代数の種類が異なるのであるから、全ての種類に対して全く同じ結果が成り立つのでは各代数の個性を見落としてしまっている可能性が高い。一方代数の種類ごとに全く異なる構成をしているのでは普遍性の面で難がある。艤装配位の場合ここで見たように、各代数で共通するべき部分では共通の定義が成り立ち、一方各々の代数の個性が最も発揮されるであろう部分には適切な違いが発生するので、筆者には好ましく感じられる[31]。

**$B(\infty)$-クリスタル**　最近 Salisbury–Scrimshaw により KR クリスタルとは別の興味深いクリスタルのクラスである $B(\infty)$-クリスタルの構造が艤装配位上に導入された [SalScr14]。この場合 0 番目の頂点も含む全ての Dynkin 図の頂点に $\nu^{(a)}$ が配置される。柏原作用素の定義など代数構造の面だけに限れば一般の対称化可能な Kac–Moody リー代数の場合にも全く同様にして定義することができる [SalScr15b]。なお、この様な考え方では本稿で考えてきた KR クリスタルの場合は 0 番目の頂点には何も配置していないので有限型と呼ばれるが、KR クリスタルの場合にはアフィン代数による無限次元対称性が深くかかわっているので、ここでの議論とはかなり異なる話なのかもしれない。

もう少し詳細に彼らの構成を記述しておこう。本稿で述べてきた艤装配位と異なる点として、$B^{r,s}$ の形状を記述している $\mu^{(a)}$ は全て空集合として扱い、vacancy number からも $\mu^{(a)}$ の関わる項を消去する。次に柏原作用素 $\widetilde{e}_a$ と $\widetilde{f}_a$ であるが、基本的な定義は本稿で述べてきた KR クリスタルの場合と全く同様である。唯一の相違点は、$\widetilde{f}_a$ の定義において rigging が対応する vacancy number を超過した場合 $\widetilde{f}_a$ の作用を 0 と定めていたのを、その様な場合でも構わず（形式的に全く同じ手順で）そのまま計算を続けていってしまう点である。元々の Schilling [Sc06] による KR クリスタルの構造の導入では $\widetilde{f}_a$ の作用が 0 となる場合を含む定義である事を本質的に用いて Stembridge の方法に帰着させていたので $B(\infty)$-クリスタルの状況には適用できない議論であった。一方論文 [Sa14] による方法は関数 $\varepsilon_a$ を艤装配位上に導入する直接的な方法であるため $B(\infty)$-クリスタルの場合にもそのまま通用し、実際論文 [SalScr14] ではその方法で導入されている。

最終的に $B(\infty)$-クリスタルは、$\nu^{(a)}=\emptyset$ が全ての頂点 $a$ について成り立つような空の艤装配位から出発して、柏原作用素 $\widetilde{f}_a$ を任意の仕方で作用させたものの全体として定義される。

$B(\infty)$-クリスタルに対する艤装配位の持つ著しい性質として、Kashiwara involution と呼ばれる重要な双対性が非常に単純な形をとることがあげられる [SalScr16]。具体的には各ストリングの rigging と corigging を入れ替える操作

$$\left(\nu_i^{(a)}, J_i^{(a)}\right) \longmapsto \left(\nu_i^{(a)}, P_{\nu_i^{(a)}}(\nu) - J_i^{(a)}\right)$$

で移った先で柏原作用素 $\widetilde{e}_a$ と $\widetilde{f}_a$ を定義すればよい。実にすっきりとした状況である。同じ操作は 48 ページで Lusztig involution を実現した際にも用いたし、また Bethe 仮説方程式の解析でも 16 ページの式 (10) で登場して基本的であった。

$B(\infty)$-クリスタルは空の艤装配位から柏原作用素で得られるもの全体として定義されたが、最も単純な $A_n$ 型の場合にはより詳細な情報を得ることができる [SalScr15a, HL]。組み合わせ論的な記述には marginally large (reverse) tableaux と呼ばれる対象が登場する。

$B(\infty)$-クリスタルの場合の艤装配位の理論は現状では柏原作用素が定義され、表現が構成されました、という段階の様に筆者には感じられる。その意味では著者らも述べているようにまだ一つの「モデル」というべき段階であろう。モデルを超えてどのように本質的な意味が存在するのか解明するのは将来の重要な課題である。

[31] 老子 第 48 章 『無為にして為さざるは無し』

**組み合わせ論等との関連**　元々艤装配位の理論は Kostka–Foulkes 多項式やその一般化に関わる組み合わせ論的な恒等式の証明を目的として導入されたものだった。これはチャージ関数（38 ページ (39) 参照）の生成母関数が Kostka–Foulkes 多項式となることによる。この方面についても多くの研究がなされており、例えば [K01, HKOTT] を参照されたい。また Kostka–Macdonald 多項式との関連 [F95] や二重 Kostka 多項式との関連 [Liu] 等も議論されている。なおアフィン量子群の表現論との関連は、例えば [Her, dFK] を参照。

**おわりに**　艤装配位の研究は 1980 年代中ごろの Kirillov 氏らの研究を嚆矢として現在まで 30 年ほどの間様々な研究者により研究が進められてきた。本稿では、少なくとも筆者の視野に収まっている範囲でその全体像をお示ししたつもりである。ようやく広がりのある世界が見えてきたところだと思う。

艤装配位の研究は、単に表現を構成する、という地点を超えて、我々の扱っている代数の表す無限次元の対称性とは何なのか、という問題について考える際に深い洞察を与えてくれるようである。無限次元代数の表す対称性は有限次元代数における感覚からは想像もつかないほど大きいため、系に強い制限を課し適用可能な対象は限定されてしまう一方、驚くほど精緻で意外性のある現象が発生することが多い。従って有限次元代数の研究において伝統的に育まれてきた方法に加え、全く異なる発想からの研究も重要になると思われる。無限次元代数の表す対称性が本格的に研究されるようになったのは歴史的に見て決して長い期間ではなく、その特有の世界についてはまだまだ解明の余地が残っているように感じられる。その時本稿で述べた様な観点からの研究も一つの考え方を提供するであろう。

艤装配位は人間の頭の中で論理的に作られたようなものではなく、無限次元の対称性に由来した自然界に実在する数学的量であると言って良いであろう。今後とも艤装配位の研究が進み、様々な姿を見せてくれるのを楽しみにしている。

**謝辞**　本稿は筆者が 2016 年 9 月 7 日に第 61 回代数学シンポジウム（日本数学会、佐賀大学）において行った講演の議事録として執筆されました。講演の機会を与えて下さった鈴木武史様をはじめとする運営委員の皆様に深く感謝いたします。また本稿の一部は 2013 年度に東京大学大学院数理科学研究科で行った集中講義における講義ノートに基づいており、今回の執筆に当たっても講義を聴講してくださった方々からの鋭いコメントが大変役に立ちました。集中講義のホストをしてくださった寺田至様にもこの場を借りて感謝いたします。最後になりましたが、本稿の内容にかかわる共同研究を行い、色々なことを教えて下さった Anatol N. Kirillov 様にもお礼申し上げます。

# 参考文献

[A] V. I. Arnold, 『古典力学の数学的方法』（岩波書店, 1980 年）

[AV] L. V. Avdeev and A. A. Vladimirov, *Exceptional solutions to the Bethe ansatz equations.* Theor. Math. Physics **69** (1986) 1071–1079.

[BMSZ] N. Beisert, J. A. Minahan, M. Staudacher and K. Zarembo, *Stringing spins and spinning strings.* JHEP **09** (2003) 010 (27pp). arXiv:hep-th/0306139

[B] H. A. Bethe, *Zur theorie der metalle.* Zeit. für Physik **71** (1931) 205–226.

[BDS] E. Brattain, N. Do and A. Saenz, *The completeness of the Bethe ansatz for the periodic ASEP*. arXiv:1511.03762

[CJ] W. Cai and N. Jing, *Applications of a Laplace–Beltrami operator for Jack polynomials*. Eur. J. Comb. **33** (2012) 556–571. arXiv:1101.5544

[DG14] T. Deguchi and P. R. Giri, *Non self-conjugate strings, singular strings and rigged configurations in the Heisenberg model*, J. Stat. Mech: Theor. Exp. (2015) P02004. arXiv:1408.7030

[DS] L. Deka and A. Schilling: *New fermionic formula for unrestricted Kostka polynomials*, Journal of Combinatorial Theory, Series A **113** (2006) 1435–1461. See also the preprint version: arXiv:math/0509194

[dFK] P. Di Francesco and R. Kedem. *Proof of the combinatorial Kirillov–Reshetikhin conjecture*. Int. Math. Res. Not. IMRN, **2008** Art. ID rnn006, 57 pages. arXiv:0710.4415

[D] P. A. M. Dirac, *Quantum mechanics of many-electron systems*. Proc. Royal Soc. of London **A123** (1929) 714–733.

[DMN] B. A. Dubrovin, V. B. Matveev and S. P. Novikov, *Non-linear equations of Korteweg-de Vries type, finite-zone linear operators, and Abelian varieties*, Russian Math. Surveys **31** (1976) 59-146.

[EKS] F. H. L. Essler, V. E. Korepin and K. Schoutens, *Fine structure of the Bethe ansatz for the spin-$\frac{1}{2}$ Heisenberg $XXX$ model*, J. Phys. A: Math. Gen. **25** (1992) 4115–4126.

[F96] L. D. Faddeev, *How Algebraic Bethe Ansatz works for integrable model.* arXiv:hep-th/9605187

[F00] L. D. Faddeev, 『現代の数理物理はどのようなものであるべきか』数理科学 2000 年 11 月号 5–11.

[F95] S. Fishel, *Statistics for special $q,t$-Kostka polynomials*, Proc. Amer. Math. Soc. **123** (1995) 2961–2969.

[F95] E. Frenkel, *Affine algebras, Langlands duality and Bethe ansatz*, in Proceedings of the International Congress of Mathematical Physics, Paris, 1994, ed. D. Iagolnitzer, pp. 606–642, International Press, 1995). arXiv:q-alg/9506003

[FOY] K. Fukuda, M. Okado and Y. Yamada: *Energy functions in box ball systems*, Int. J. Modern Phys. **A15** (2000) 1379–1392. arXiv:math/9908116

[GHNS] A. M. Gainutdinov, W. Hao, R. I. Nepomechie and A. J. Sommese, *Counting solutions of the Bethe equations of the quantum group invariant open XXZ chain at roots of unity*, J. Phys. A **48** (2015) 494003 (38pp). arXiv:1505.02104

[GN] A. M. Gainutdinov and R. I. Nepomechie, *Algebraic Bethe ansatz for the quantum group invariant open XXZ chain at roots of unity*, Nucl. Phys **B909** (2016) 796–839. arXiv:1603.09249

[GD] P. R. Giri and T. Deguchi, *Heisenberg model and rigged configurations*, J. Stat. Mech: Theor. Exp. P07007 (2015). arXiv:1501.07801

[HNS13] W. Hao, R. I. Nepomechie and A. I. Sommese, *Completeness of solutions of Bethe's equations*, Phys. Rev. E **88** (2013) 052113 (8pp plus supplementary materials). arXiv:1308.4645

[HNS14] W. Hao, R. I. Nepomechie and A. J. Sommese, *Singular solutions, repeated roots and completeness for higher-spin chains*, J. Stat. Mech. **2014** (2014) P03024 (20pp). arXiv:1312.2982

[HHIKTT] G. Hatayama, K. Hikami, R. Inoue, A. Kuniba, T. Takagi and T. Tokihiro: *The* $A_M^{(1)}$ *automata related to crystals of symmetric tensors*, J. Math. Phys. **42** (2001) 274–308. arXiv:math/9912209

[HKOT] G. Hatayama, A. Kuniba, M. Okado and T. Takagi *Combinatorial R matrices for a family of crystals:* $B_n^{(1)}$, $D_n^{(1)}$, $A_{2n}^{(2)}$ *and* $D_{n+1}^{(2)}$ *cases*, J. Algebra **247** (2002) 577–615. arXiv:math/0012247

[HKOTT] G. Hatayama, A. Kuniba, M. Okado, T. Takagi and Z. Tsuboi, *Paths, Crystals and Fermionic Formulae*, Prog. Math. Phys. **23** (2002) 205–272. arXiv:math/0102113

[H] T. Hawkins, *Jacobi and the birth of Lie's theory of groups*. Archive for History of Exact Sciences **42** (1991) 187–278.

[H] W. Heisenberg, *Zur theorie des ferromagnetismus*. Zeit. Physik **49** (1928) 619–636. D. Delphenich 氏による英訳が以下の場所にある（項目は Physics: Electromagnetism)。 http://www.neo-classical-physics.info

[Her] D. Hernandez, *The Kirillov-Reshetikhin conjecture and solutions of T-systems*. J. für die reine und angewandte Math. **596** (2006) 63–87. arXiv:math/0501202

[HL] J. Hong and H. Lee, *Rigged Configuration Descriptions of the Crystals* $B(\infty)$ *and* $B(\lambda)$ *for Special Linear Lie Algebras*, arXiv:1604.04357v1

[J] C. G. J. Jacobi, *Lectures on dynamics* (1842–43). English translation by K. Balagangadharan (Hindustan Book Agency, 2009).

[JM] M. Jimbo and T. Miwa, *Solitons and infinite dimensional Lie algebras*, Publ. RIMS. Kyoto Univ. **19** (1983) 943–1001.

[K90] V. G. Kac, *Infinite Dimensional Lie Algebras*, third edition, Cambridge Univ. Press (1990).

[KR] V. G. Kac and A. K. Raina, *Bombay lectures on highest weight representations of infinite dimensional Lie algebras*. (World Scientific, 1987).

[$(KMN)^2$] S. J. Kang, M. Kashiwara, K. C. Misra, T. Miwa, T. Nakashima and A. Nakayashiki, *Affine crystals and vertex models*. Int. J. Modern Phys. **A7** (1992), 449–484.

[K91] M. Kashiwara, *On crystal bases of the q-analogue of universal enveloping algebras*. Duke Math. J. **63** (1991), 465–516.

[KN] M. Kashiwara and T. Nakashima, *Crystal graphs for representations of the q-analogue of classical Lie algebras*. J. Algebra **165** (1994), 295–345.

[KKR] S. V. Kerov, A. N. Kirillov and N. Y. Reshetikhin, *Combinatorics, Bethe Ansatz, and representations of the symmetric group*. J. Sov. Math. **41** (1988) 916–924.

[K01] A. N Kirillov, *Bijective correspondences for rigged configurations*, (Russian) Algebra i Analiz **12** (2000) 204–240, translation in St. Petersburg Math. J. **12** (2001) 161–190.

[KR] A. N. Kirillov and N. Yu. Reshetikhin: *The Bethe ansatz and the combinatorics of Young tableaux*, J. Sov. Math. **41** (1988) 925–955.

[KS09] A. N. Kirillov and R. Sakamoto, *Relationships Between Two Approaches: Rigged Configurations and 10-Eliminations*, Lett. Math. Phys. **89** (2009) 51–65. arXiv:0902.2286

[KS14a] A. N. Kirillov and R. Sakamoto, *Singular solutions to the Bethe ansatz equations and rigged configurations*, J. Phys. A: Math. Theor. **47** (2014) 205207 (20pp). arXiv1402.0651

[KS14b] A. N. Kirillov and R. Sakamoto, *Some remarks on Nepomechie–Wang eigenstates for spin 1/2 XXX model*, Moscow Math. J. **15** (2015) 337–352. arXiv:1406.1958

[KS15] A. N. Kirillov and R. Sakamoto, *Bethe's Quantum Numbers And Rigged Configurations*, Nucl. Phys. **B905** (2016) 359–372. arXiv:1509.02305

[KSS] A. N. Kirillov, A. Schilling and M. Shimozono, *A bijection between Littlewood–Richardson tableaux and rigged configurations*, Selecta Math. (N.S.) **8** (2002) 67–135. arXiv:math/9901037

[K06] A. Kitaev, *Anyons in an exactly solved model and beyond*. Ann. Phys. **321** (2006) 2–111. arXiv:cond-mat/0506438

[K70] D. E. Knuth, *Permutations, matrices, and generalized Young tableaux*, Pacific J. Math. **34** (1970) 709–727.

[KOSTY] A. Kuniba, M. Okado, R. Sakamoto, T. Takagi and Y. Yamada, *Crystal interpretation of Kerov–Kirillov–Reshetikhin bijection*, Nucl. Phys. **B740** (2006) 299–327. arXiv:math/0601630

[KS06] A. Kuniba and R. Sakamoto, *The Bethe ansatz in a periodic box-ball system and the ultradiscrete Riemann theta function*, J. Stat. Mech. (2006) P09005 (12pp). arXiv:math/0606208

[KSY] A. Kuniba, R. Sakamoto and Y. Yamada, *Tau functions in combinatorial Bethe ansatz*, Nucl. Phys. **B786** (2007) 207–266. arXiv:math/0610505

[KTT] A. Kuniba, T. Takagi and A. Takenouchi *Bethe ansatz and inverse scattering transform in a periodic box-ball system*, Nucl. Phys. **B747** (2006) 354–397. arXiv:math/0602481

[L] T. Lam, *Loop symmetric functions and factorizing matrix polynomials*. Fifth International Congress of Chinese Mathematicians. Part 1, 2, 609–627, AMS/IP Stud. Adv. Math., 51, pt. 1, 2, Amer. Math. Soc., Providence, RI, 2012. arXiv:1012.1262

[LP] T. Lam and P. Pylyavskyy, *Intrinsic energy is a loop Schur function*. J. Comb. 4 (2013) **4** 387–401. arXiv:1003.3948

[LPS1] T. Lam, P. Pylyavskyy and R. Sakamoto, *Box-Basket-Ball Systems*, Rev. Math. Phys. **24** (2012) 1250019 [23 pages]. arXiv:1011.5930

[LPS2] T. Lam, P. Pylyavskyy and R. Sakamoto, *Rigged Configurations and Cylindric Loop Schur Functions*, Ann. Inst. Henri Poincaré (to appear). arXiv:1410.4455

[LS94] R. P. Langlands and Y. Saint-Aubin, *Algebro-geometric aspects of the Bethe equations*. in Proceedings of the Gursey Memorial Conference I–Strings and Symmetries, ed. G. Atkas, C. Saclioglu, M. Serdaroglu, Lecture Notes in Physics, vol. 447, Springer Verlag (1994), pp40–53.

[LS97] R. P. Langlands and Y. Saint-Aubin, *Aspects combinatoires des équations de Bethe*. in Advances in Mathematical Sciences: CRM 's 25 years, ed. Luc Vinet, CRM Proceedings and Lecture Notes vol.11 Am. Math. Soc. (1997) pp231–301.

[LS05] C. Lecouvey and M. Shimozono, *Lusztig's q-analogue of weight multiplicity and one-dimensional sums for affine root systems*, Adv. Math. **208** (2007) 438–466. arXiv:math/0508511

[LOS] C. Lecouvey, M. Okado and M. Shimozono, *Affine crystals, one-dimensional sums and parabolic Lusztig q-analogues* Math. Zeit. **271** (2012) 819–865. arXiv:1002.3715

[LR] S.-S. Lin and S.-S. Roan, *Bethe Ansatz for Heisenberg XXX Model*, arXiv:cond-mat/9509183

[Liu] S. Liu, *Fermionic formula for double Kostka polynomials*. arXiv:1602.08792

[L] J. Lützen, *Joseph Liouville 1809–1882. Master of pure and applied mathematics*. (Springer, 1990).

[MIT] J. Mada, M. Idzumi and T. Tokihiro, *On the initial value problem of a periodic box-ball system*, J. Phys. A: Math. Gen. **39** L617–L623. arXiv:nlin/0608037

[MM] D. N. Mermin and H. Wagner, *Absence of ferromagnetism or antiferromagnetism in one-or two-dimensional isotropic Heisenberg models*. Phys. Rev. Lett. **17** (1966) 1133–1136.

[MY] K. Mimachi and Y. Yamada, *Singular vectors of the Virasoro algebra in terms of Jack symmetric polynomials*. Comm. Math. Phys. **174** (1995) 447–455.

[MW] K. C. Misra and E. A. Wilson, *Soliton cellular automaton associated with $D_n^{(1)}$-crystal $B^{2,s}$*, J. Math. Phys. **54** (2013) 043301 (33pp)

[M] N. F. Mott, 『初歩の量子力学』(共立出版, 1977).

[NW] R. I. Nepomechie and C. Wang, *Algebraic Bethe ansatz for singular solutions*, J. Phys. A: Math. Theor. **46** (2013) 325002 (8pp). arXiv:1304.7978

[OS] M. Okado and R. Sakamoto, *Stable Rigged Configurations for Quantum Affine Algebras of Nonexceptional Types*, Adv. Math. **228** (2011) 1262–1293. arXiv:1008.0460

[OSS11] M. Okado, R. Sakamoto and A. Schilling: *Affine crystal structure on rigged configurations of type $D_n^{(1)}$*, J. Algebraic Comb. **37** (2013) 571–599. arXiv:1109.3523

[OSSS] M. Okado, R. Sakamoto, A. Schilling and T. Scrimshaw, *Type $D_n^{(1)}$ rigged configuration bijection*, J. Algebraic Comb. (to appear). arXiv:1603.08121

[OSano] M. Okado and N. Sano, *KKR type bijection for the exceptional affine algebra $E_6^{(1)}$*, Contemp. Math. **565** (2012) 227–242. arXiv:1105.1636

[OSS03] M. Okado, A. Schilling and M. Shimozono: *A crystal to rigged configuration bijection for nonexceptional affine algebras*, in "Algebraic Combinatorics and Quantum Groups" (Ed. N. Jing), World Scientific (2003), 85–124. arXiv:math/0203163

[PS] V. Pasquier and H. Saleur, *Common structures between finite systems and conformal field theories through quantum groups*. Nucl. Phys. **B330** (1990) 523–556.

[Sa06] R. Sakamoto, *Crystal Interpretation of Kerov-Kirillov-Reshetikhin Bijection II. Proof for $sl_n$ Case*, J. Algebraic Comb. **27** (2008) 55–98. arXiv:math/0601697

[Sa07] R. Sakamoto, *Kirillov–Schilling–Shimozono bijection as energy functions of crystals*, Int. Math. Res. Notices (2009) 2009: 579–614. arXiv:0711.4185

[Sa12a] R. Sakamoto, 超離散ソリトン系と組み合わせ的表現論、数理科学 2012 年 1 月号

[Sa12b] R. Sakamoto, *Ultradiscrete Soliton Systems and Combinatorial Representation Theory*, RIMS Kokyuroku 1913 (2014) 141–158. arXiv:1212.2774

[Sa14] R. Sakamoto, *Rigged Configurations and Kashiwara Operators*, SIGMA 10 (2014), 028, 88 pages. arXiv:1302.4562 (式番号などはリンク先の v2 に基づく)

[Sa15] R. Sakamoto, *Rigged configurations approach for the spin-1/2 isotropic Heisenberg model*. J. Phys. A: Math. Theor. **48** (2015) 165201. arXiv:1501.00532

[SSAFR] R. Sakamoto, J. Shiraishi, D. Arnaudon, L. Frappat and E. Ragoucy, *Correspondence between conformal field theory and Calogero–Sutherland model*. Nucl. Phys. **B704** (2005) 490–509. arXiv:hep-th/0407267

[SalScr14] B. Salisbury and T. Scrimshaw, *A rigged configuration model for $B(\infty)$*, J. Combin. Theory Ser. A 133 (2015), 29–57. arXiv:1404.6539

[SalScr15a] B. Salisbury and T. Scrimshaw, *Connecting marginally large tableaux and rigged configurations via crystals.* Algebr. Represent. Theory **19** (2016), 523–546. arXiv:1505.07040

[SalScr15b] B. Salisbury and T. Scrimshaw, *Rigged configurations for all symmetrizable types.* arXiv:1509.07833

[SalScr16] B. Salisbury and T. Scrimshaw, *Rigged configurations and the $*$-involution.* arXiv:1601.06137

[Sc04] A. Schilling: *A bijection between type $D_n^{(1)}$ crystals and rigged configurations*, J. Algebra **285** (2005) 292–334. arXiv:math/0406248

[Sc06] A. Schilling, *Crystal structure on rigged configurations*, Int. Math. Res. Notices **2006** (2006) Article ID 97376, Pages 1–27. arXiv:math/0508107

[Sc08] A. Schilling: *Combinatorial structure of Kirillov–Reshetikhin crystals of type $D_n^{(1)}$, $B_n^{(1)}$, $A_{2n-1}^{(2)}$*, J. Algebra **319** (2008) 2938–2962. arXiv:0704.2046

[ScScr14] A. Schilling and T. Scrimshaw, *Crystal structure on rigged configurations and the filling map.* Electronic J. Comb. **22** (2015) #P1.73 (54 pages) arXiv:1409.2920

[ScSh] A. Schilling and M. Shimozono: *$X = M$ for symmetric powers*, J. Algebra **295** (2006) 562–610. arXiv:math/0412376

[SRFB] U. Schollwöck, J. Richter, D. J. J. Farnell and R. F. Bishop (Eds.), *Quantum magnetism.* Lecture Notes in Physics 645 (Springer, 2004).

[Scr16] T. Scrimshaw, *Rigged configurations as tropicalizations of loop Schur functions.* arXiv:1607.03232

[Scr15] T. Scrimshaw, *A crystal to rigged configuration bijection and the filling map for type $D_4^{(3)}$*, J. Algebra **448** (2016) 294–349. arXiv:1505.05910

[Sh98] M. Shimozono, *Affine type A crystal structure on tensor products of rectangles, Demazure characters, and nilpotent varieties*, J. Algebraic Combin. **15** (2002) 151–187. arXiv:math/9804039

[Sh05] M. Shimozono, *On the $X = M = K$ conjecture*, arXiv:math/0501353

[SZ] M. Shimozono and M. Zabrocki, *Deformed universal characters for classical and affine algebras*, J. of Algebra **299** (2006) 33–61. arXiv:math/0404288

[Su71] B. Sutherland, *Exact results for a quantum many-body problem in one dimension.* Phys. Rev. **A4** (1971) 2019–2021.

[Su72] B. Sutherland, *Exact results for a quantum many-body problem in one dimension II.* Phys. Rev. **A5** (1972) 1372–1376.

[T] M. Takahashi, *One-dimensional Heisenberg model at finite temperature*, Prog. Theor. Phys. **46** (1971) 401–415.

[TS] 武部 尚志, 関谷 信寛, 『可解格子模型と共形場理論の話題から』上智大学数学講究録 47 (2006).

[YY] C. N. Yang and C. P. Yang, *One-dimensional chain of anisotropic spin-spin interactions. I. Proof of Bethe's hypothesis for ground state in a finite system.* Phys. Rev. **150** (1966) 321–327.

[Y00] Y. Yamada, *A birational representation of Weyl group, combinatorial R-matrix and discrete Toda equation*, in Physics and Combinatorics (Nagoya, 2000), Editors A. N. Kirillov and N. Liskova, World Sci. Publ., River Edge, NJ (2001), 305–319.

[Y13] F. Yura, 箱とバスケットと玉の系におけるソリトン解、九州大学応用力学研究所研究集会報告 No.24AO-S3, (2013) 156–161.

[YT] F. Yura and T. Tokihiro, *On a Periodic Soliton Cellular Automaton*, J. Phys. A: Math. Gen. **35** (2002) 3787–3801. arXiv:nlin/0112041

[V] A. A. Vladimirov, *Proof of the invariance of the Bethe-ansatz solutions under complex conjugation*, Theor. Math. Phys. **66** (1986) 102–105.

**Invitation to the Bethe ansatz** (by Reiho Sakamoto)
Proceedings for the 61st Symposium on Algebra (Mathematical Society of Japan, 2016)

Abstract: We review the algebraic Bethe ansatz for the Heisenberg model. The exposition includes some of recent advancements with emphasis on a relation with the rigged configurations. We also provide somewhat thorough review of the crystal bases and the rigged configurations.

In particular, we provide the inverse scattering transform for the type $D_n^{(1)}$ box-ball systems, see (54) at page 51. The result is analogous to the type $A_n^{(1)}$ case and is a consequence of [OSSS]. We compare the result with Misra–Wilson's box-ball system [MW].

We also provide a reformulation of a result of [Sa07] at page 33. More precisely, we provide a crystal interpretation of the algorithm of the type $A_n^{(1)}$ rigged configuration bijection by the simplified quantity (35), instead of $\varepsilon^{(a)}_{l,j,k}$ of [Sa07]. Let the columns of a tableau $b_k \in B^{r_k,s_k}$ be $c_1, c_2, \ldots, c_{s_k}$ from left to right. Then we introduce the intermediate tensor product $C_j$ as in (34). Here we use the Kashiwara convention for the tensor products of crystals. Then we compute the sum of energy functions as in (32) to derive (35). As for the definitions used in the sum (32), see (31) where $u^{a,l} \in B^{a,l}$ is the classical highest weight element and $\overset{R}{\simeq}$ stands for the isomorphism under the combinatorial $R$-matrix.

Contents:

*Date: Dec. 14, 2016*